\documentclass{itogi2017}

 \currentyear{}
 \currentvolume{}

\usepackage[russian]{babel}
\usepackage[cp1251]{inputenc}
\usepackage{graphicx}

\begin{document}
\centerline
{\bf Holomorphically Homogeneous Real Hypersurfaces in $\mathbb C^3$
}

\

\centerline{Loboda A.V.}

\

{\bf Abstract:}
We give a complete description and classification of locally homogeneous real hypersurfaces in $\mathbb C^3$. Various groups of mathematicians have been studying this problem in the last 25 years, and several significant classes of  hypersurfaces under consideration have been studied and classified. The final results in the classification problem presented in this paper are obtained by using the classification of abstract 5-dimensional real Lie algebras, and by studying their representations by algebras of holomorphic vector fields in complex 3-space. The complete list of pairwise inequivalent hypersurfaces that we obtain contains 47 types of homogeneous hypersurfaces; some of the types are 1- or 2-parametric families, and each of the others is single hypersurface or a finite family of hypersurface.

This paper is to appear in the Proceeding of the Moscow Mathematical Society ("Trudy Moskovskogo Matematicheskogo Obshchestva"); it was submitted to the journal in March 2020.

The current version is the original Russian one, and the English version is to appear soon. The complete classification of locally homogeneous hypersurfaces in $\mathbb C^3$ is given on pages 53-54.

\

УДК {517.55, 515.172.2, 512.816}

\

\centerline{\bf
Голоморфно-однородные вещественные гиперповерхности в $\Bbb C^3 $
}

\centerline{Лобода А.В.
}

* Работа выполнена при поддержке РФФИ (проекты NN 17-01-00592, 20-01-00497).

\

{\bf Аннотация.} Обсуждается полное решение задачи описания и классификации (с локальной точки зрения) голоморфно однородных вещественных гиперповерхностей пространства $\mathbb{C}^3$. В рамках различных подходов к этой задаче за 25 лет ее исследования несколькими
группами математиков были изучены большие семейства таких многообразий.
Итоговые результаты в этой задаче получены автором настоящей статьи с использованием классификации абстрактных 5-мерных алгебр
Ли и техники их голоморфных представлений в трехмерном комплексном пространстве.

  Полный список попарно не эквивалентных голоморфно однородных гиперповерхностей в полученной классификации содержит 47 типов таких многообразий, включающих отдельные поверхности, а также одно-параметрические и двух-параметрические семейства поверхностей.

\

{\bf Ключевые слова:} однородное многообразие, вещественная гиперповерхность, нормальная форма, голоморфное преобразование, векторное поле,
алгебра Ли.

\

\

\centerline{\bf Содержание}

Введение

1. Основные определения и некоторые подходы к задаче

2. Общие формулировки теорем об однородности в $ \Bbb C^3 $

3. Схема изучения <<просто однородных>> гиперповерхностей в $ \Bbb C^3 $

\qquad 3.1. Классификация Мубаракзянова 5-мерных алгебр Ли

\qquad 3.2. Совместное упрощение вида векторных полей

\qquad 3.3. Изученные блоки 5-мерных алгебр Ли

4. Однородные поверхности для блока из 11 пятимерных алгебр

\qquad 4.1. Первый случай: выпрямление базиса идеала $ \mathfrak h_2 = <e_1, e_2, e_4> $

\qquad 4.2. Второй случай: выпрямление базиса идеала $ \mathfrak h_2 = <e_1, e_3, e_4> $

5. Однородные поверхности для блока из 7 пятимерных алгебр

\qquad 5.1. Орбиты алгебр $ \mathfrak g_{5,38}$ и $ \mathfrak g_{5,39} $

\qquad 5.2. Семейства алгебр $ \mathfrak g_{5,33}, \mathfrak g_{5,34}, \mathfrak g_{5,35}$: первый случай

\qquad 5.3. Семейства алгебр $ \mathfrak g_{5,33}, \mathfrak g_{5,34}, \mathfrak g_{5,35}$: второй случай

\qquad 5.4. Орбиты алгебры $ \mathfrak g_{5,36}$

\qquad 5.5. Упрощение базисов голоморфных реализаций алгебры $ \mathfrak g_{5,37} $

6. Интегрирование <<трудных>> алгебр и итоговые выводы

\qquad 6.1. Сферичность Леви-невырожденных орбит алгебры $ \mathfrak g_{5,37} $

\qquad 6.2. Описание орбит семейства алгебр $ \mathfrak g_{5,35} $: второй случай

\qquad 6.3. Проверка гипотезы о новизне орбит алгебр $ \mathfrak g_{5,35} $

\qquad 6.4. Итоговые выводы о <<просто однородных>> гиперповерхностях

Приложение. Полный список однородных гиперповерхностей в в $ \Bbb C^3 $

Список литературы

\

{\bf Введение}

  Изучение свойств границ многомерных комплексных областей и,
более широко, вещественных гиперповерхностей многомерных комплексных пространств, является одной из традиционных задач многомерного комплексного анализа.
Так, работа Пуанкаре [1], в которой сравнивались бидиск и шар в $ \Bbb C^2 $
с точки зрения их границ, показала наличие голоморфно различных свойств у этих границ, в силу чего существенно различаются и
голоморфные свойства самих этих областей.

  Начало изучению голоморфно однородных вещественных гиперповерхностей многомерных комплексных пространств было положено работой
Э. Картана [2]. В ней дано описание (с локальной и глобальной точек зрения) голоморфно однородных вещественных гиперповерхностей 2-мерных комплексных пространств.
   Идея Картана о связи решаемой им задачи с соответствующими изучаемым поверхностям алгебрами Ли (голоморфных векторных полей)
позволила относительно легко получить полный список однородных гиперповерхностей в $\Bbb C^2 $: список всех 3-мерных алгебр Ли
в начале 20-го века был уже известен, и список этот достаточно краток.

   Эта же идея Картана была, в конце концов, реализована и в случае 3-мерного комплексного пространства. Настоящая статья подводит
завершающие итоги 25-летнего изучения задачи описания (локально) голоморфно однородных вещественных гиперповерхностей в $ \Bbb C^3 $.

   При обсуждении ее полного решения естественно выделить три части, на которые оказалась разделена
полученными результатами совокупность изучаемых
однородных многообразий. Одну из таких частей составляют вырожденные по Леви однородные гиперповерхности, описанные в работе [3].
В отличие от локальной сводимости любой вырожденной однородной гиперповерхности в $ \Bbb C^2 $ к вещественной гиперплоскости,
количество в $ \Bbb C^3 $ аналогичных однородных многообразий оказалось весьма представительным.

  Еще более обширной оказалась вторая часть, содержащая Леви-невырожденные однородные гиперповерхности в $ \Bbb C^3 $ с
<<богатыми>> алгебрами симметрий.
Примеры таких многообразий, от обычной сферы $ S^5 \subset \Bbb C^3 $  и поверхности Винкельманна
$
    v = z_1 \bar z_2 + z_2 \bar z_1 + |z_1|^4
$
до более <<экзотических>> объектов, постепенно появлялись на протяжении всего периода изучения задачи. Полный же их список был
предъявлен в работе [4] лишь три года назад.

Эта часть семейства однородных гиперповерхностей, самая большая из трех по числу представителей,
содержит, в частности, обобщения (на 3-мерный случай) проективно однородных поверхностей Картана и большое количество
трубчатых многообразий с аффинно однородными основаниями. В определенном смысле наиболее интересным здесь является семейство
из 8 типов голоморфно однородных трубчатых гиперповерхностей, основания которых не являются аффинно однородными. В пространстве $ \Bbb C^2 $ такая ситуация невозможна, что объясняется строгой псевдо-выпуклостью всех невырожденных по Леви поверхностей в двумерном случае. Существование же такого семейства в $ \Bbb C^3 $ характеризует связь свойств голоморфной и аффинной однородности в многомерных комплексных пространствах
как более тонкую по сравнению с представлявшейся на начальном этапе исследований.

  Третья, последняя, часть всего семейства, содержащая <<просто однородные>> гиперповерхности пространства $ \Bbb C^3 $,
предполагалась изначально самой обширной. Однако исследование возможных 5-мерных орбит в $ \Bbb C^3 $ всех 5-мерных алгебр Ли,
опирающееся на известный обширный список [5] (содержащий 67 типов алгебр), показало обратную картину. <<Просто однородных>>,
т.е. имеющих дискретные стабилизаторы, строго псевдо-выпуклых гиперповерхностей в $ \Bbb C^3 $ вообще не существует.
А <<просто однородных>> гиперповерхностей с индефинитной формой Леви в этом пространстве имеется лишь три. С точностью до биголоморфной эквивалентности это поверхности с уравнениями
$$
    v\, ( 1 \pm y_2 x_2) = y_1 y_2 \quad \mbox{и}
\quad  (v - x_2 y_1)^2 + y_1^2 y_2^2 = y_1.
$$

  Изучение третьего подсемейства голоморфно однородных гиперповерхностей в $ \Bbb C^3 $ потребовало,
при всей краткости его итогового описания, большой технической работы.
При этом
  5-мерные орбиты в $ \Bbb C^3 $, отвечающие большей части упомянутого списка 5-мерных алгебр [5],
были рассмотрены в серии совместных работ [6-10] автора данной статьи. В то же время два блока из семи, на которые весь список [5] был разделен самим автором, обсуждаются в последних разделах настоящей  статьи.
Новых однородных поверхностей в этих блоках в итоге не оказалось, но идея Картана, перенесенная на 3-мерный случай,
реализована полностью.

   Отметим одну особенность изучения третьего подсемейства однородных поверхностей.
При обсуждении орбит 5-мерных алгебр часто удается понять, не выписывая конкретных уравнений этих орбит, что
многие из них сводятся либо к поверхностям из двух первых (уже описанных) частей изучаемой совокупности, либо к
трубчатым поверхностям с аффинно-однородными основаниями в $ \Bbb R^3 $, список которых также известен (см. [11]).
Заинтересованность только в новых однородных объектах позволяет в итоге существенно сократить обсуждения.

  Уточним, что основной технический прием, позволивший осуществить перенос идеи Картана на 3-мерный случай, был разработан в [12].
Этот прием связан с совместным выпрямлением (нормализацией) нескольких голоморфных векторных полей в пространстве $ \Bbb C^3 $.
  На заключительных этапах
исследования получаемых орбит и проверки их возможной эквивалентности известным однородным поверхностям часто оказывается полезной техника нормальных форм уравнений самих изучаемых поверхностей. Первые результаты в обсуждаемой задаче были получены автором статьи
именно с применением этой техники.

   Ниже обсуждаются (с большей или меньшей степенью подробности) все три названные выше подсемейства однородных
гиперповерхностей в $ \Bbb C^3 $. При этом желание автора подробно изложить завершающий фрагмент решения объемной задачи
не позволило прокомментировать более детально два первых подсемейства. Их описания имеются и доступны в оригинальных текстах,
что, возможно, извиняет чрезмерный акцент автора настоящей статьи на собственные результаты.

  Отметим еще два момента, связанные с общей характеристикой обсуждаемых в работе результатов:

1. размерность пространства модулей отдельных компонент (стратов) общей совокупности голоморфно однородных гиперповерхностей 3-мерных
комплексных пространств не превышает 2;

2. все однородные гиперповерхности в изучаемой ситуации допускают реализации, описываемые элементарными функциями.

\

   Автор выражает глубокую благодарность участникам семинара Витушкина (Белошапке В.К., Немировскому С.Ю., Чирке Е.М.,
Сергееву А.Г., Кружилину Н.Г.) за постоянное внимание к его научной активности и, в частности, за интерес к
подготовке этой статьи. Отдельные благодарности автор выражает Винбергу Э.Б., Горбацевичу В.В., Циху А.К., Шмальцу Г., Лопатину И. А.
Одним из первых заинтересовался исследованиями автора в задаче об однородности и проявлял к ним постоянный интерес Исаев А.В.,
недавно ушедший из жизни. Этот интерес также поддерживал автора в многолетней работе над задачей.

\

$ ^* $ От рецензента статьи автор получил (помимо множества замечаний) интересный вопрос
 о возможной сводимости голоморфной однородности в обсуждаемой ситуации к однородности бирациональной. Проверка этой гипотезы предполагает
приложение существенных дополнительных усилий, связанных с исследованием локальных групп преобразований всех 47 типов однородных гиперповерхностей из итогового списка настоящей работы. Автор глубоко признателен рецензенту за высказанные им замечания, способствовавшие существенному улучшению текста статьи, и за этот вопрос, но в настоящее время не имеет возможностей дополнить статью таким исследованием.

\

{\bf 1.  Основные определения и некоторые подходы к задаче
}

\

   Везде ниже мы опираемся на
традиционное определение однородного многообразия (см., например, [13, с. 444]) и используем его модификации, учитывающие
локальный характер голоморфных функций и отображений и естественность перехода от групп Ли к соответствующим алгебрам Ли.

\textbf{\it Определение 1.} Многообразие $ M $ {\it однородно
относительно некоторой группы Ли} (преобразований) $ G $, если эта
группа транзитивно действует на $ M $, то есть любую точку из $M$
можно перевести в любую другую точку этого многообразия
подходящим преобразованием из группы $ G $.

  Переходя от групповых преобразований к векторным полям (инфинитезимальным преобразованиям), можно
переформулировать первое определение. В случае локальных голоморфных преобразований ростков подмногообразий
(в частности, гиперповерхностей) многомерных комплексных пространств так получается более естественное

\textbf{\it Определение 2.}
Росток $ M_{\xi} $ вещественно-аналитической гиперповерхности $ M \subset \Bbb C^n $ (с центром в точке $ \xi \in M $) мы называем голоморфно однородным, если алгебра Ли $ g(M) $ касательных к $ M $ ростков голоморфных векторных полей (определенных в точке $ {\xi} $) накрывает своими значениями всю касательную плоскость $ T_{\xi}M $.

{\bf Замечание 1.} В работе [14] обсуждаются несколько эквивалентных определений однородных CR-многообразий. В интересующем нас случае
вещественных гиперповерхностей одно из них совпадает с определением 2.

{\bf Замечание 2.} Введенную в определении 2 алгебру $ g(M) $ ВСЕХ касательных к $ M $ ростков голоморфных векторных полей мы будем называть далее {\it алгеброй симметрий} (ростка) однородной гиперповерхности $ M $.

{\bf Замечание 3.} Ниже мы будем использовать вариант определения 2, связанный с наличием на ростке однородной
гиперповерхности $ M \subset \Bbb C^3 $ некоторой алгебры Ли $ g $ ростков касательных к $ M $ голоморфных векторных полей, имеющей (полный) ранг 5 в центре ростка $ p $. При этом алгебра $ g $ может быть лишь подалгеброй алгебры симметрий $ g(M) $,
не обязательно совпадающей с $ g(M) $.

  Отметим, что этот вариант определения 2 можно перенести
на случай алгебр Ли произвольной природы (например, аффинных)
и на произвольные размерности обсуждаемого многообразия $ M $ и объемлющего пространства. При этом смысл условия на ранг алгебры
(полнота ранга) в том, что значения векторных полей (соответствующей природы) в каждой
точке $ p\in M $ накрывают всю касательную плоскость $ T_p(M) $.

  В изучаемой нами многомерной задаче о голоморфной однородности результаты об однородности аффинной
(причем, мало-размерной и вещественной) имеют важное значение. Это связано с содержанием
и результатом основополагающей статьи Э. Картана [2] 1932 г., первой работы о голоморфной (псевдо-конформной) однородности
вещественных гиперповерхностей многомерных комплексных пространств. Полное описание таких гиперповерхностей в $ \Bbb C^2 $,
полученное Картаном, существенно использует описание аффинно-однородных плоских кривых.
Приведем здесь этот результат.

{\bf ТЕОРЕМА 1.1 (E. Cartan).} {\it Произвольная голоморфно-однородная
ве\-щест\-венная гиперповерхность 2-мерного комплексного пространства
голоморфно эквивалентна вблизи любой своей точки
\newline
либо 1) трубке над одной из аффинно-различных кривых:
$$
   a)\ y = x^s \ (-1 \le s < 1),
\qquad
   b) \ y = \ln x,
\qquad
   c) \ y = x \ln x,
\eqno (1.1)
$$
$$
  d) \ r = e^{a\varphi} \
   (r-\mbox{полярный радиус},\varphi - \mbox{полярный угол},  a \ge 0),
$$
либо
2) одной из проективно-однородных поверхностей следующего семейства:
$$
  a) \ 1 + |z|^2 + |w|^2 = a | 1 + z^2 + w^2 | \ \ ( a > 1 ),
\eqno (1.2)
$$
$$
  b) \ 1 + |z|^2 - |w|^2 = a | 1 + z^2 - w^2 | \ \ ( a > 1 ),
$$
$$
\qquad \
  c) \  |z|^2 + |w|^2 - 1 = a | z^2 + w^2 - 1 | \ \ ( 0 < | a | < 1 ).
$$
}

{\bf Замечание.} Набор кривых 1a) --- 1d) из формулировки теоремы представляет собой полный (с точностью до аффинной эквивалентности) список плоских аффинно-однородных кривых.

  Отталкиваясь от этой теоремы Э. Картана естественно ставить вопросы о голо\-морфно-однородных вещественных
гиперповерхностях и вещественных подмногообразиях произвольной коразмерности в
комплексных пространствах больших размерностей.
В рамках
развивающегося многомерного комплексного анализа, а именно в теории CR-многообразий и CR-отображений,
такие вопросы ставились и публиковались соответствующие работы.

    Упомянем здесь, например, работу [15]
об однородных абстрактных CR-мно\-гообразиях. Отметим, однако, что, в отличие
от статьи Э. Картана, ни одного примера конкретного однородного многообразия авторы этой статьи не приводят, ограничиваясь изучением
схем Дынкина алгебр Ли, допускающих реализацию в виде алгебр векторных полей на абстрактных CR-многообразиях.

   Укажем также на <<параллельное>> с комплексным анализом изучение вопросов об однородности в рамках вещественной геометрии.
Так, описание аффинно-однородных кривых на плоскости (см., например, [16]) было получено, как принято считать, школой Бляшке на рубеже 19-го и 20-го веков.
Поверхности 3-мерного вещественного пространства,
однородные относительно различных подгрупп
аффинной группы,
были описаны в середине 20-го века. Тогда же эти описания были
включены в учебники по дифференциальной геометрии. Например, в
книге [17] приведен <<полный>> (оказавшийся все-таки не полным [18]) список эквиаффинно-однородных
поверхностей пространства $ \Bbb R^3 $.

   Ясно, что даже не полные списки аффинно-однородных гиперповерхностей пространства $ \Bbb R^n $ позволяют получать
достаточно большое количество примеров голоморфно однородных гиперповерхностей пространства $ \Bbb C^n $ (с тем же $ n $)
за счет рассмотрения т.н. трубчатых многообразий (трубок) вида
$$
   M = \Gamma + i \Bbb R^n, \ \Gamma \in \Bbb R^n.
\eqno (1.3)
$$

  С учетом получения в 1995-96 годах в работе [11] полного списка аффинно-однородных гиперповерхностей в $ \Bbb R^3 $
представлялась относительно простой задача описания голоморфно-однородных гиперповерхностей в $ \Bbb C^3 $ (например, в духе Э. Картана).
  Как часто бывает, вопрос оказался <<несколько>>
сложнее начальных представлений о нем. Тем не менее, по
прошествии 25 лет усилиями целого ряда математиков, включая и автора данной
статьи, поставленная задача оказалась решена. Именно ее решению посвящено все нижеследующее изложение.

  Уточним, что задача была сформулирована в рамках семинара по комплексному анализу на мех-мате МГУ,
которым в 80-90-е годы прошлого века руководили Шабат Б.В., Гончар А.А. и Витушкин А.Г. В то время
популярными темами семинара были голоморфные отображения вещественных подмногообразий многомерных комплексных пространств
(или, в современной терминологии, вопросы CR-геометрии) и, в частности, голоморфные свойства вещественных гиперповерхностей
([19-22]).
  Значительный интерес участников семинара вызвала статья [23] и возможные приложения ее результатов в
близких задачах комплексного анализа (см. [24-27]).

  Автор настоящей статьи занялся обсуждаемой в ней задачей приблизительно за год до выхода статьи [11]. К этому времени
сложилось системное осмысление результата Картана, представляемое, например, в следующей краткой  форме:

   среди голоморфно однородных вещественных гиперповерхностей в двумерном комплексном пространстве имеются две
(с точностью до локальной голоморфной эквивалентности) <<исключительных>> поверхности: вещественная гиперплоскость и сфера.
Их алгебры голоморфных симметрий являются, соответственно, бесконечномерной и 8-мерной, что было известно еще в начале 20-го века
(см. [1]). Для всех остальных однородных поверхностей из списка Картана алгебры голоморфных симметрий имеют размерность 3,
совпадающую с размерностью самих этих поверхностей.

   При этом появляющиеся примеры и отдельные факты о голоморфно однородных вещественных гиперповерхностях 3-мерных комплексных пространств
наглядно отличались от предложенного выше описания ситуации с однородностью в двумерном случае.

  Во-первых, в трехмерном комплексном пространстве были известны (см., например, [23]) две однородные вещественные гиперповерхности
$$
   Im\,z_3 = |z_1|^2 + |z_2|^2 \, \quad \mbox{и} \quad Im\,z_3 = |z_1|^2 - |z_2|^2,
\eqno (1.4)
$$
на которых алгебры голоморфных векторных полей имеют одинаковые размерности 15. У других поверхностей,
обладающих свойством невырожденности по Леви  (в т.ч., голоморфно однородных), такие алгебры имеют меньшую размерность.

  <<Естественной>> размерностью таких алгебр для других однородных гиперповерхностей в $ \Bbb C^3 $ могла бы быть
размерность самих этих поверхностей, т.е. 5.
Однако в 1994 г. появился пример (см. [28]) алгебраической гиперповерхности 4-й степени (т.н. поверхности Винкельманна)
$$
    Im\, z_3 = (z_1 \bar z_2 + z_2 \bar z_1) + |z_1|^4
\eqno (1.5)
$$
в пространстве трех комплексных переменных. Эта поверхность
является голоморфно однородной, и группа ее голоморфных преобразований имеет размерность 8, не совпадающую ни с 5, ни с 15 !

   Второе, наиболее существенное, отличие ситуации с гиперповерхностями двумерных комплексных пространстве от пространств
большей размерности (и, в частности, трехмерных) проявляется в такой
голоморфной характеристике гиперповерхностей, как {\it форма Леви}.

  Напомним, что любую вещественно-аналитическую гиперповерхность в многомерном комплексном пространстве с координатами
$ z_1, ..., z_n $ можно задать
вблизи неособой точки уравнением, разрешенным относительно одной из вещественных переменных
$$
    Im\, z_n = F(z_1, ... , z_{n-1}).
\eqno (1.6)
$$

  При этом можно считать, что в правой части этого уравнения нет линейных слагаемых, так что
комплексная касательная плоскость к обсуждаемой поверхности $ M $ имеет уравнение $ z_n = 0 $.
Формой Леви поверхности $ M $ (в начале координат) в такой ситуации является ([29]) эрмитова форма
$$
    \sum_{k,j=1}^{n-1} \frac{\partial^2 F }{\partial z_k \partial{\bar z_j}}(0) z_k \bar{z_j}.
\eqno (1.7)
$$

   Эта форма может быть вырожденной или невырожденной, а в случае размерности $ n $ обсуждаемого пространства, большей чем 2,
возникают различные возможности для
сигнатуры этой формы. В трехмерном пространстве $ \Bbb C^3 $ форма Леви гиперповерхности зависит от двух переменных,
а потому возможен случай ее знакоопределенности (положительный и отрицательный случаи в локальной голоморфной геометрии не различаются),
знаконеопределенности (индефинитный случай) и вырождения.
    В случае двумерного пространства возможны лишь (положительная) определенность или вырожденность эрмитовой формы от
одной комплексной переменной.

   Гиперповерхности с положительно определенной формой Леви называются {\it строго псевдо-выпуклыми} (СПВ-случай); {\it
индефинитными} естественно называть гиперповерхности, форма Леви которых является знаконеопределенной (и невырожденной).

Отметим, что единичная сфера
$ |z_1|^2 + ... + |z_n|^2 = 1 $
в комплексном пространстве любой размерности является СПВ-гиперповерхностью в каждой своей точке. В
соответствии с представлением (1.13), ее уравнение можно привести (дробно-линейным преобразованием) к виду
$$
   Im\,z_n = |z_1|^2 + ... + |z_{n-1}|^2.
$$

  Тем самым, первая из двух поверхностей (1.4) есть (локальный) голоморфный образ сферы;
она является СПВ-гиперповерхностью. Обе поверхности второй степени (1.4)
часто называют {\it квадриками Мозера} или {\it сферическими}, а поверхности, не сводимые голоморфными преобразованиями ни к
одной из них --- несферическими.
   Поверхность Винкельманна является индефинитной несферической гиперповерхностью;  вторая из
поверхностей (1.4) --- индефинитная сферическая гиперповерхность.

  Отличия в аналитических описаниях вещественно-аналитических гиперповерхностей двумерного и трехмерного
(<<существенно>> многомерного) пространств достаточно наглядно иллюстрирует понятие {\it нормальной формы} Мозера ([23]).
Смысл этого понятия состоит в упрощении за счет голоморфных преобразований уравнения вида (1.6) для произвольной аналитической гиперповерхности комплексного $ n $-мерного пространства, имеющей невырожденную форму Леви.
В частности, для интересующего нас случая $ n = 3 $ уравнение (1.6) Леви-невырожденной поверхности можно
довести до состояния
$$
    v = \, < z,z > + \sum_{k,l \ge 2, m \ge 0 } N_{klm}(z,\bar z) u^m.
\eqno (1.8)
$$

   Здесь и в дальнейшем переменные в пространстве $ \Bbb C^3 $ обозначаются через $ z_1, z_2, w $; у переменной $ w $
выделены вещественная и мнимая часть, т. что $ u = Re\,w, v = Im\,w $; через $ z $ обозначен комплексный вектор $ (z_1,z_2) $,
и, соответственно, $ \bar z = (\bar z_1, \bar z_2) $.

   Через $ < z,z > $ в (1.15) обозначена форма Леви поверхности, а $ N_{klm}(z,\bar z) $ --- однородный многочлен степени
 $ k $ по переменной $ z $ и степени $ l $ по $ \bar z $. Многочлены $ N_{22m} $, $ N_{32m} $, $ N_{33m} $ удовлетворяют некоторым дополнительным ограничениям (т.н. tr-условиям). В силу таких условий самый младший из таких многочленов, т.е. $ N_{220} $, определяется пятью (вместо формально возможных девяти) вещественными коэффициентами. В случае знакоположительной формы Леви
$
   < z,z > = |z_1|^2 +|z_2|^2
$
этот многочлен является гармоническим относительно переменных $ z_1, z_2 $, т.е. $ \Delta N_{220} = 0 $.
А в случае знаконеопределенной формы Леви
$
   < z,z > = z_1  \bar z_2 + z_2  \bar z_1
$
он имеет вид (см. [30])
$$
   N_{220}(z, \bar z)  = \lambda_1 |z_1|^4 + \lambda_2 \left(4|z_1|^2 |z_2|^2 - (z_1^2 \bar z_2^2 + z_2^2 \bar z_1^2)\right) +
\lambda_3 |z_2|^4 +
\eqno (1.9)
$$$$
\qquad \qquad
+ i \mu_1 (z_1^2 \bar z_1 \bar z_2 - z_1 z_2 \bar z_1^2 ) + i \mu_2 (z_1 z_2 \bar z_2^2 - z_2^2 \bar z_1 \bar z_2 ),
$$
где $ \lambda_1, \lambda_2, \lambda_3, \mu_1, \mu_2 \in \Bbb R $.

  При переходе от одного нормального уравнения поверхности к другому многочлен $ N_{220}(z,\bar z) $ преобразуется в
$ \rho N_{220}(Uz, \overline{Uz}) $, где
$ \rho $ ---некоторое ненулевое вещественное число, $ U $ --- унитарное (или псевдо-унитарное) линейное преобразование.

  Если (в фиксированной точке поверхности) этот многочлен является нулевым, то поверхность называется {\it омбилической} (в этой точке).
Омбиличность голоморфно однородной поверхности хотя бы в одной точке означает, что поверхность является сферической (т.е. голоморфно эквивалентна одной из поверхностей (1.4)).

  Многочлены $ N_{klm} $ из уравнения (1.15) суммарных весов $ k+l+2m > 4 $ также важны при обсуждении свойства однородности, однако они
могут существенно изменяться под действием голоморфных преобразований, сохраняющих (в целом) нормальный вид уравнения поверхности. Для фиксированной
гиперповерхности $ M $ семейство таких преобразований (и, тем самым, алгебра голоморфных векторных полей на $ M $) зависит не более чем от 15 вещественных параметров.

   Упомянем для сравнения с описанной ситуацией нормальные уравнения Мозера в пространстве двух комплексных переменных $ z, w $. Здесь
форму Леви можно всегда считать равной $ |z|^2 $, а первым голоморфным инвариантом
поверхности является многочлен шестого веса $ N_{420} $ (а не 4-го, как в трехмерном случае).

 В работе [31] показано, что голоморфно однородная невырожденная по Леви вещественная гиперповерхность двумерного
комплексного пространства однозначно определяется тройкой вещественных коэффициентов своего {\it специального} нормального уравнения. Вся
эта тройка удовлетворяет двум уравнениям второй степени, а потому семейство (локально) голоморфно однородных Леви-невырожденных (а следовательно, псевдо-выпуклых) вещественных гиперповерхностей в $ \Bbb C^2 $ может быть описано посредством одного вещественного  параметра.

Установленное в [31] полное соответствие такого описания однородных поверхностей с работой Э. Картана укрепило надежды на эффективность применения аналогичных коэффициентных подходов к задаче об однородности гиперповерхностей в $ \Bbb C^3 $. А начальным результатом в этой задаче, полученным с использованием нормальных форм Мозера, можно считать опубликованную годом
ранее работу [32]. В ней были получены оценки размерностей групп изотропии (или стабилизаторов) для
вещественных Леви-невырожденных гиперповерхностей в $ \Bbb C^3 $.

   Следствием этих оценок является утверждение о том, что поверхность Винкельманна --- единственная
(с точностью до голоморфной эквивалентности)  голоморфно однородная Леви-невырожденная вещественная
гиперповерхность трехмерного
комплексного пространства, имеющая 8-мерную алгебру симметрий. Для всех остальных однородных Леви-невырожденных
вещественных гиперповерхностей в $ \Bbb C^3 $ эта размерность не превышает 7.

Напомним, что
оценка снизу числом 5 (т.е. размерностью самих таких многообразий) для размерности алгебры голоморфных векторных
полей на голоморфно однородных гиперповерхностях в $ \Bbb C^3 $  является очевидной. Отметим еще, что
однородные гиперповерхности пространства $ \Bbb C^3 $, (максимальные) алгебры симметрий которых имеют размерность, равную именно 5,
естественно называть {\it просто однородными} ([33]) в отличие от {\it кратно транзитивных} ([4]) однородных многообразий.

  В последующие несколько лет на основе коэффициентного подхода автором настоящей статьи были получены списки голоморфно
однородных гиперповерхностей в $ \Bbb C^3 $, имеющих в точности 7-мерные алгебры симметрий (СПВ-случай -- [34],
индефинитный случай -- [30]). Также с использованием метода нормальных форм был получен список однородных СПВ-гиперповерхностей
 $ \Bbb C^3 $ с 6-мерными алгебрами ([35]), дополненный в работе [4] одной <<потерянной>> поверхностью.

   Попытки исследования индефинитных однородных гиперповерхностей в $ \Bbb C^3 $ в рамках того же коэффициентного подхода
оказались менее эффективными. Они
показали, во-первых, что количество таких поверхностей, имеющих, например, 6-мерные алгебры симметрий, значительно превосходит
аналогичное количество СПВ-поверхностей. Кроме того, важный в таких обсуждениях многочлен $ N_{220} $ из нормального уравнения
невырожденной вещественно-аналитической гиперповерхности имеет в индефинитном случае гораздо больше возможных
(не сводимых друг к другу) форм, чем в СПВ-случае. Приведем здесь конкретный результат (см. [30] или [36]) о таких многочленах,
который потребуется нам в последних разделах статьи.

{\bf Предложение 1.2.} {\it Пару многочленов $ (<z,z>, N_{220}) $ из нормального уравнения неомбилической индефинитной гиперповерхности
$ M \subset \Bbb C^3 $ можно преобразовать голоморфной заменой координат так, что
$
   < z,z > = z_1  \bar z_2 + z_2  \bar z_1,
$
а набор коэффициентов
$ \lambda_1, \lambda_2, \lambda_3, \mu_1, \mu_2 $
многочлена (1.9) перейдет в одно из семи следующих состояний:

$
    1) \ (1, A, 1, 0, 0), |A| > 1/3;
$

$
   2) \ (1,A,-1, 0, 0), A\in \Bbb R;
$

$
  3) \ (1, \pm 1, 0, 0, 0);
$

$
   4) \ (1, 1/3, 1, 0, 0) $ --- Картанов тип;

$
  5) \ (1, -1/3, 1, 0, 0)  $ --- псевдо-Картанов тип;

$
   6) \ (0, 0, 0, 1, 0) $ --- <<антитрубчатый>> тип;

$
  7) \ (1, 0, 0, 0, 0) $ --- Винкельманнов тип.
}

{\bf Замечание.} Изучение вопроса
о существовании однородных индефинитных поверхностей в каждом из этих семи классов привело в свое
время лишь к отдельным примерам.
Из полного описания
голоморфно однородных поверхностей в $ \Bbb C^3 $, приводимого в следующих разделах статьи,
ответы на этот вопрос выводятся достаточно легко.

  Описанные здесь результаты автора являются лишь частичными продвижениями в задаче описания голоморфно однородных гиперповерхностей
в $ \Bbb C^3 $. В силу их частного характера сами полученные в названных работах списки однородных
поверхностей здесь лишь упоминаются. В следующих разделах статьи такие списки будут приведены в полном объеме и представлены с
точки зрения других подходов к задаче.

   Основным из них можно считать подход, основанный на детальном изучении алгебраических
структур, связанных с однородными гиперповерхностями. Его эффективность более чем наглядно иллюстрируется, например, работами
[11] и [4], обсуждаемыми в следующих разделах настоящей статьи. Результаты этих работ составляют значительную часть
представляемого полного описания однородных гиперповерхностей в $ \Bbb C^3 $.

   Важно отметить, что в работе [4] изучаются лишь невырожденные по Леви поверхности. В связи с этим
отдельное описание вырожденных по Леви однородных гиперповерхностей можно считать
естественным направлением в изучении общей задачи об однородности.
Такое описание Леви-вырожденных голоморфно однородных гиперповерхностей 3-мерных комплексных пространств (или в более современных терминах, классификация вырожденных однородных 5-мерных CR-многообразий)
было получено в работе [3]. В основе этой работы также лежит использование достаточно тонких алгебраических свойств
вырожденных поверхностей и соответствующих им алгебр Ли.

   Еще одним из подходов к задаче об однородности можно назвать исследование строго псевдо-выпуклых поверхностей.
В этом более узком классе поверхностей имеются свои упрощения, геометрия СПВ-многообразий обеспечивает их дополнительные, более жесткие
свойства по сравнению с общей ситуацией. Однако важные эффекты типа существования поверхности Винкельманна (со всеми ее экстремальными
свойствами) при таком <<одностороннем>> подходе могут быть не замечены. Поэтому естественно сочетать в решении поставленной задачи разные  подходы: понятие однородности позволяет использовать при ее исследовании как упомянутые выше алгебраические и аналитические методы, так и
инструментарий дифференциальной геометрии.

   В качестве иллюстрации использования геометрических свойств СПВ-гиперпо\-верхностей в обсуждаемой задаче можно привести уже упомянутые
работы [34-35], а также работу [33]. Полученное в последней работе утверждение о голоморфной сводимости
<<просто однородных>> СПВ-гиперповерхностей к трубчатым многообразиям не является общим для всех однородных поверхностей, но показывает очередной интересный эффект, связанный с геометрическим свойством строгой псевдо-выпуклости.

  Говоря о возможных подходах к задаче об однородности в широком смысле слова, имеет смысл упомянуть и топологические ее аспекты,
развиваемые, например, в работе [37].

\

Отметим еще несколько достаточно тонких моментов обсуждаемой задачи,
которые необходимо учитывать при ее полном исследовании в рамках любых подходов.

 1) В первую очередь упомянем вопрос о возможной эквивалентности получаемых (например, разными методами) однородных поверхностей.
Так, в списке аффинно однородных плоских кривых имеются две разные кривые $ y = x^2 $ и $ y = e^x $; соответственно,
трубки в $ \Bbb C^2 $ над ними имеются в списке Картана. Однако в работе [2] имеется уточнение о {\it локальной} голоморфной эквивалентности обеих этих трубок сфере, а значит, и друг другу.
  Подобные уточнения необходимо иметь в виду и при построении списков голоморфно однородных поверхностей в $ \Bbb C^3 $,
многие из которых также сводятся к трубкам над поверхностями из $ \Bbb R^3 $, а некоторые из них <<неожиданно>> оказываются сферическими.

 2)  При исследовании однородных поверхностей с точки зрения их алгебр Ли, необходимо помнить, что отдельная алгебра может иметь
голоморфно неэквивалентные орбиты в обсуждаемом многомерном комплексном пространстве. Пример такой неоднозначности доставляют проективно-однородные поверхности Э. Картана в $ \Bbb C^2 $; имеются такие примеры и в 3-мерном пространстве, как будет показано ниже.

   С другой стороны, могут быть голоморфно эквивалентными орбиты
двух различных (в смысле алгебраического изоморфизма) алгебр Ли ростков голоморфных векторных полей в $ \Bbb C^3 $,
имеющих ранг 5.
 Утверждение о различии орбит в такой ситуации
может быть справедливым при некоторых дополнительных допущениях (типа <<простой однородности>>, трубчатости и др.) относительно обсуждаемых поверхностей. Но, вообще говоря, однородная поверхность может быть орбитой {\it разных} алгебр Ли.

В трехмерном комплексном пространстве такой эффект часто оказывается связан с наличием на 5-мерной однородной поверхности $ M $
(максимальной) алгебры $ g(M) $
голоморфных векторных полей, имеющей размерность, бОльшую чем 5, но неизвестную заранее. В
подобной ситуации возможно существование у алгебры $ g(M) $, например, 5-мерных подалгебр (одной или нескольких), каждая из которых по-отдельности
обеспечивает однородность обсуждаемой поверхности. При получении завершающих классификационных выводов об отдельных
однородных многообразиях необходим учет такого обстоятельства.

   3)  В решаемой нами задаче составления полного списка однородных гиперповерхностей приходится использовать уже
известные большие блоки таких поверхностей, а также (не менее обширные) описания, например, 5-мерных алгебр Ли.
В таких <<вспомогательных>> результатах имеются как опечатки, так и фактические ошибки (на одну из таких ошибок,
допущенную автором данной статьи в работе [35], справедливо указано в [4]; в книге [38] и статье [7] отмечена одна
достаточно значимая неточность в статье [5]; еще одна аналогичная неточность той же работы [5], не зафиксированная
в книге [38], обсуждается в разделе 6 настоящей статьи). Разумеется, в работе автора над настоящей статьей, как и в
любой объемной работе, также естественно ожидать некоторое количество погрешностей.

   В такой ситуации автор настоящей статьи выражает надежду, что
количество неточностей в предлагаемом итоговом <<полном>> списке однородных поверхностей окажется не очень значительным, а сам такой
список может оказаться полезным другим исследователям.

Подводя итоги этого раздела, подчеркнем еще раз, что в нем излагаются лишь
предварительные результаты, связанные с задачей описания голоморфно
однородных гиперповерхностей в $ \Bbb C^3 $.
  Вместе с тем, эти результаты по сути задают схему разбиения обсуждаемой сложной задачи на несколько
менее объемных блоков. Полное решение основной задачи будет описано в следующих разделах статьи именно в рамках
такой схемы.

    Отметим еще, что при
описании голоморфно-однородных гиперповерхностей в комплексных пространствах
следующих размерностей, несомненно, потребуются другие идеи и средства. Количество таких многообразий, как и соответствующих им алгебр Ли,
растет весьма быстро, как показывает сравнение двумерной картановской классификации с итоговым списком однородных гиперповерхностей настоящей работы.

В то же время отдельные фрагменты и идеи данной статьи могут оказаться продуктивными в более высоких размерностях. Так,
в настоящее время готовится публикация некоторых результатов об орбитах нильпотентных 7-мерных вещественных алгебр Ли в пространстве $ \Bbb C^4 $, полученных (с участием автора настоящей статьи) на основе техники,  описываемой ниже для 3-мерного случая.

\

{\bf 2. Общие формулировки теорем об однородности в $ \Bbb C^3 $
}

\

   Изложение предыдущего раздела и результаты, полученные в задаче об однородности в последние 15 лет, позволяют
считать естественной схему изучения этой задачи в рамках трех ее блоков:

1. Вырожденные по Леви однородные гиперповерхности в $ \Bbb C^3 $.

2. Леви-невырожденные однородные гиперповерхности в $ \Bbb C^3 $ с <<богатыми>> алгебрами симметрий (multiply-transitive hypersurfaces).

3. <<Просто однородные>> Леви-невырожденные гиперповерхности в $ \Bbb C^3 $.

   Напомним, что первый блок из этого набора описан в работе [3] 2008 г.; итоговому описанию поверхностей второго
блока посвящена работа [4], опубликованная в 2017 г. Наконец, завершающий третий блок, как и объемная схема его изучения,
будут описаны ниже.

   Отметим, что предлагаемые схемы используют еще несколько работ, которые нельзя не упомянуть в контексте изложения
достаточно продолжительной истории построения полного списка голоморфно-однородных гиперповерхностей в $ \Bbb C^3 $.
   Помимо процитированной выше работы Картана, важную роль в решении задачи о голоморфной однородности в трехмерном
комплексном случае сыграла техника голоморфной реализации (голоморфных представлений) 5-мерных алгебр Ли в $ \Bbb C^3 $,
разработанная в статье [12] и фактически переносящая идеи Картана на эти размерности.

  Работой же, которую можно первой отнести к изучаемой задаче, является
статья Дуброва-Комракова-Рабиновича [11]. Этими авторами в 1995 г. был представлен полный список аффинно-однородных поверхностей 3-мерного вещественного пространства. Несмотря на вещественный характер решенной в ней задачи, именно трубки над
такими поверхностями составляют значительную часть общего списка
голоморфно-однородных гиперповерхностей пространства $ \Bbb C^3 $.

  В частности, трубчатые гиперповерхности (трубки) над аффинно-однородными основаниями присутствуют во всех трех блоках
представляемой в этой работе классификации. Приведем здесь результат [11].

{\bf ТЕОРЕМА 2.1 ([11]).} {\it
Всякая локально однородная поверхность в 3-мерной вещественной
аффинной геометрии является открытым подмножеством либо
некоторой поверхности второго порядка, либо
цилиндра над одной из однородных плоских кривых из теоремы 0.1 Картана, либо
аффинно эквивалентна открытому подмножеству одной из следующих поверхностей
($\alpha, \beta $ - вещественные параметры):

\
1) $ z = x^{\alpha} y^{\beta} $,
\qquad
\qquad
\qquad
\qquad
\qquad
2) $ z = (x^2 + y^2)^{\alpha} e^{\beta arg (x+iy)} $,

\
3) $ z = \ln x + \alpha \ln y $,
\qquad
\qquad
\qquad
\quad
4) $ z =  arg (x + i y ) + \beta \ln (x^2 + y^2) $,

\
5) $ z =  \ln (x^2 + y^2) $,
\qquad
\qquad
\qquad
\quad
6) $ z = x( \alpha \ln x + \ln y ) $,

7) $ (z - xy + x^3/3 )^2 = \alpha (y - x^2/2)^3, $
\
8) $  z = y^2 \pm e^x $,

\
9) $  z = y^2 \pm x^{\alpha} $,
\qquad
\qquad
\qquad
\qquad  \
10) $  z = y^2 \pm \ln x $,

\
11) $  z = y^2 \pm x \ln x $,
\qquad
\qquad
\qquad
\quad
12) $  z = x y + e^x $,

\
13) $  z = x y +  x^{\alpha} $,
\qquad
\qquad
\qquad
\qquad
14) $  z = x y + \ln x $,

\
15) $  z = x y + x \ln x $,
\qquad
\qquad
\qquad
\ \
16) $  z = x y + x^2 \ln x $,

\
17) $ x z = y^2 \pm x^{\alpha} $,
\qquad
\qquad
\qquad
\quad
\quad
18) $ x z = y^2 \pm x \ln x $,

\
19) $ x z = y^2 \pm x^2 \ln x $.
}

{\bf Замечание.} Имеется лишь 5 типов поверхностей 2-го порядка, являющихся строго выпуклыми или седловидными (только над такими поверхностями трубки являются невырожденными по Леви):

1. Эллипсоид (сфера),

2-3. Два гиперболоида (однополостный и двухполостный),

4-5. Два параболоида (эллиптический и гиперболический).

  При этом оба параболоида имеются в общем списке теоремы 1.1: гиперболический параболоид с уравнением $ z = x y $ входит в п.1
при $ \alpha = \beta = 1 $, а эллиптический $ z= x^2 + y^2 $  --- в п. 2 при $ \alpha = 1, \beta = 0 $.

  Можно отметить, что идея получения списка [11], базирующаяся на использовании алгебр Ли,
отвечающих аффинно однородным поверхностям, достаточно проста. Сначала из геометрических соображений делается оценка
сверху для возможных размерностей алгебр аффинных векторных полей, отвечающих однородным поверхностям в $ \Bbb R^3 $.
После этого перебор всех таких возможных алгебр (количество которых оказывается вполне обозримым) и их интегрирование
приводят к итоговому списку однородных поверхностей.
   При этом полученный в [11] окончательный результат поставил точку в достаточно продолжительной активности
в задаче об аффинно однородных поверхностях в $\Bbb R^3 $ ряда других математиков, и, в частности, представителей бельгийской школы
дифференциальной геометрии (см., например, [39]).

   Отметим также, что в списке теоремы 1.1 имеются поверхности, трубки над которыми в пространстве $ \Bbb C^3 $
вырождены по Леви. К ним, в частности, относятся конус и цилиндры над плоскими кривыми a), b), c), d) из п.1 Теоремы 1.1 Картана.
   Полное описание вырожденных по Леви (локально) голоморфно однородных гиперповерхностей в $ \Bbb C^3 $,
полученное Фелсом и Каупом и составляющее первый блок обсуждаемой полной классификации, приводится ниже.

{\bf ТЕОРЕМА 2.2 ([3]).}
{\it Каждая локально однородная относительно голоморфных преобразований
вещественная гиперповерхность 3-мерного комплексного пространства, имеющая
в фиксированной точке вырожденную форму Леви,
голоморфно эквивалентна вблизи этой точки одному из следующих
многообразий:

$
  1) \ M = F_k + i\, \Bbb R^3,
$
где $ F_k\, (k = 1,.., 5) $ - одна из перечисленных ниже поверхностей пространства $ \Bbb R^3_{(x,y,z)} $
(с какими-либо значениями параметров);

$
 \quad 1.1. \
 x^2 + y^2 = z^2 \ (z > 0),
$

$
\quad
  1.2. \ z = |x+iy| e^{ \omega \arg (x+iy)}, \
$

$\quad
  1.3. \ z = x (\ln y - \ln x), \
$

$
\quad  1.4. \ z = x^{1 - \theta} y^{\theta} \ (x>0, y>0),
$

$\quad
  1.5. \ (z - 3 xy + 2 x^3)^2 = 4(x^2 - y)^3.
$

$
  2) \ M = \Bbb C \times M',
$
где $ M' $ - одна из невырожденных по Леви голоморфно однородных
вещественных гиперповерхностей из списка Картана (Теорема 1.1);

$
  3) \ M = \Bbb C^2 \times \Bbb R .
$
}

{\bf Замечание 1.} В оригинальной работе [3] поверхности 1.2 -- 1.5 заданы в параметрической форме. Здесь мы предпочли их координатные описания, наглядно показывающие пересечения списков из теорем 2.1 и 2.2
и выделяющие в классификационной теореме 1.1 поверхности в $ \Bbb R^3 $ с нулевой гауссовой кривизной. Вырожденными по Леви являются
трубки в $ \Bbb C^3 $ именно над такими поверхностями.

   Описание следующего, второго, блока из обсуждаемой полной классификации однородных гиперповерхностей в $ \Bbb C^3 $
является более объемным.
  Приведем краткую формулировку соответствующей теоремы, а затем развернутые комментарии к ней.

{\bf ТЕОРЕМА 2.3 ([4]).}
{\it
 Любая <<кратно-транзитивная>> Леви-невы\-рожденная гиперповерхность в $ \Bbb C^3 $
локально биголоморфно эквивалентна одному из многообразий следующего списка:

(1) максимально симметричные гиперповерхности
$$
   v = |z_1|^2 \pm  |z_2|^2
$$
с 15-мерными алгебрами симметрий;

(2) трубчатые гиперповерхности с размерностями алгебр симметрий 6, 7 или 8;

(3) Картановы гиперповерхности и кватернионные модели, алгебры симметрий которых являются вещественными формами алгебры
$ \mathfrak{so}(4,\Bbb C ) \cong \mathfrak{sl}(2,\Bbb C) + \mathfrak{sl}(2,\Bbb C) $;

(4) гиперповерхности Винкельманнова типа, имеющие 6-мерные алгебры симметрий.
}

  В качестве первого комментария к этой теореме уточним, что все трубчатые поверхности, заявленные в п. (2) теоремы, сведены авторами в две
таблицы (Таблицы 7 и 8).

   При этом в первой таблице, содержащей 11 типов голоморфно однородных трубчатых гиперповерхностей, имеется
6 типов, совпадающих (при замене переменной $ u $ на $ z $) с пп. 12 -16 списка из теоремы 1.1 (трубка над поверхностью
Викельманна $ v = xy+x^4 $, имеющая 8-мерную алгебру симметрий, и другие трубки $ v = xy + x^{\alpha} $ при $ \alpha \in \Bbb R \setminus \{0,1,2,3,4 \} $ занимают в таблице 7 из [4] две строки). Основания таких трубок являются аффинно однородными поверхностями в $ \Bbb R^3 $.

  Пять новых (по сравнению с [11]) типов голоморфно однородных трубчатых гиперповерхностей в $ \Bbb C^3 $,
основания которых НЕ являются аффинно однородными, приведены ниже.

$
  1) \  v = (1 + \exp{2x_1}) (x_2 + \ln(1 + \exp{2x_1})),
$

$
  2) \  v = x_1 x_2 + x_1^3 \ln x_1,
$

$
  3) \ v = x_2\exp(x_1) + \exp(\alpha x_1),
$

$
  4) \ v \cos x_1 + x_2 \sin x_1 = \exp(\alpha x_1),
$

$
  5) \ v = x_1 e^{x_2} + x_1^2.
$

  Аналогично, в таблице 8, содержащей 10 типов однородных трубчатых поверхностей новыми (по сравнению с [11]) являются лишь три:

$
   6) \  v = \alpha \ln(1 + \exp{2x_1}) + \ln x_2,
$

$
   7) \ v = \alpha \ln(1 + \exp{2x_1}) + \ln(1 + \exp{2x_2}),
$

$
   8) \ v = x_2 + \varepsilon \ln(1 + \exp{2x_1}).
$

   Остальные трубки из таблицы 8 совпадают (с точностью до незначительных переобозначений параметров) с шестью типами поверхностей
из пп. 3), 4), 9), 10), 11), 18) теоремы 1.1.

{\bf Замечание.} У части из этих 8 типов алгебра симметрий является 7-мерной, а у другой части --- 6-мерной.

 В качестве еще одного общего типа поверхностей в таблице 8 выделены трубки над поверхностями
второго порядка
$$
    u^2 + \varepsilon_1 x^2 + \varepsilon_2 y^2 = 1 \ (\varepsilon_1, \varepsilon_2 = \pm 1),
\eqno (2.1)
$$
обладающие 6-мерными алгебрами симметрий и
не имеющие персональных номеров в перечне теоремы 1.1.

\

  Второй комментарий касается обобщений проективно однородных поверхностей Картана из пространства  $ \Bbb C^2 $ или
поверхностей Картанова типа. В однородных координатах $ (\xi_0 : \xi_1 : \xi_2 : \xi_3) $ пространства  $ \Bbb{CP}^3 $
любая из этих поверхностей может быть задана уравнением
$$
   \varepsilon_0 |\xi_0|^2 + \varepsilon_1 |\xi_1|^2 + \varepsilon_2 |\xi_2|^2 + \varepsilon_3 |\xi_3|^2  =
     \beta |\varepsilon_0 \xi_0^2 + \varepsilon_1 \xi_1^2 + \varepsilon_2 \xi_2^2 + \varepsilon_3 \xi_3^2|,
\eqno (2.2)
$$
где $ \varepsilon_k = \pm 1, \ \beta\in \Bbb R $.

   В зависимости от набора знаков параметров $ \varepsilon_k $ на таких поверхностях транзитивно действуют 6-мерные группы $ SO(4) $
или $ SO(p,q) \ (p + q = 4) $. Соответствующие этим группам алгебры Ли являются вещественными формами алгебры
$ \mathfrak{so}(4,\Bbb C ) \cong \mathfrak{sl}(2,\Bbb C) + \mathfrak{sl}(2,\Bbb C) $. В [4] выписаны естественные ограничения на параметр $ \beta $, связанные с сигнатурами набора
$ \varepsilon_k $: например, из положительности всей четверки этих коэффициентов очевидно следует, что $ \beta > 1 $.

  Помимо естественных обобщений поверхностей Картана в п. (3) теоремы 2.3 работы [4] вводится еще одно семейство поверхностей, т.н.
однородные <<кватернионные  модели>>, алгебра симметрий для которых имеет структуру
$$
     \mathfrak{so}^*(4) \cong \mathfrak{sl}(2,\Bbb R) + \mathfrak{su}(2).
$$

    Уравнения этих новых однородных поверхностей в (однородных) координатах имеют достаточно сложный вид
$$
   Im( \bar \xi_1 \xi_3 + \bar \xi_2 \xi_4) = \gamma \sqrt{ |Re(\bar \xi_1 \xi_3 + \bar \xi_2 \xi_4)|^2 + |\xi_1 \xi_4 - \xi_2 \xi_3|^2},
\quad \gamma \in \Bbb R\setminus\{0\}.
\eqno (2.3)
$$

   Наконец, комментируя последний п. (4) теоремы 2.3, рассмотрим поверхности Винкельманнова типа. Отметим, что исходя из
своего определения этого типа, авторы работы [4] выписывают здесь одно семейство поверхностей, зависящее от двух
вещественных параметров и две одиночных поверхности:

$ (i) \ Im(w + \bar z_1 z_2) = (z_1)^{\alpha}\overline{(z_1)^{\alpha}} $, где $ \alpha \in \Bbb C \setminus \{1, 0, 1, 2\} $;

$ (ii) \ Im(w + \bar z_1 z_2) = \exp(z_1 + \bar z_1) $;

$ (iii) \ Im(w + \bar z_1 z_2) = \ln(z_1) \ln(\bar z_1) $.

 При этом вторая и третья одиночные поверхности, как констатируют сами авторы, голоморфно эквивалентны, соответственно,
трубчатым гиперповерхностям
$$
     v = x_1 x_2 + e^{x_1}
\quad \mbox{и} \quad
     v = x_1 e^{x_2} + x_1^2,
$$
уже упоминавшимся выше. Все три выписанные типа поверхностей имеют 6-мерные алгебры симметрий.

\

   Завершая описание семейства голоморфно однородных гиперповерхностей пространства $ \Bbb C^3 $,
сформулируем итоговое утверждение, доказательству которого посвящена остальная часть настоящей статьи. Это утверждение
относится к однородным поверхностям, алгебры голоморфных векторных полей на которых являются в точности 5-мерными.
Оно показывает значение теорем 1.1 -- 1.3 в заявленной схеме полного описания однородных гиперповерхностей в $ \Bbb C^3 $.

{\bf ТЕОРЕМА 2.4.} {\it Пусть $ M $ -- невырожденная по Леви голоморфно однородная (вблизи некоторой точки $ p \in M $)
гиперповерхность в $ \Bbb C^3 $, а $ g(M) $ -- 5-мерная алгебра голоморфных векторных полей на $ M $, имеющая (вблизи $ p $)
ранг, равный 5. Тогда $ M $ голоморфно эквивалентна (вблизи $ p $) либо одной из известных поверхностей, упомянутых
в теоремах 1.1 -- 1.3, либо одной из следующих новых <<просто однородных>> поверхностей:
$$
    v\, ( 1 \pm y_2 x_2) = y_1 y_2.
\eqno (2.4)
$$
$$
  (v - x_2 y_1)^2 + y_1^2 y_2^2 = y_1.
\eqno (2.5)
$$
}

{\bf Замечание 1.} Однородная поверхность $ M $в $ \Bbb C^3 $, на которой имеется 5-мерная алгебра Ли $ g(M) $
голоморфных векторных полей полного ранга,
может иметь и бОльшую по размерности алгебру касательных к ней полей. Такая ситуация является достаточно
распространенной, а ответ на вопрос о
новизне орбиты 5-мерной алгебры $ g(M) $ в обсуждаемой задаче описания всех однородных
гиперповерхностей часто не является прозрачным.

{\bf Замечание 2.} Отметим, что в работе [33] решена задача описания всех голоморфно однородных СПВ-гипер\-поверхностей
пространства $ \Bbb C^3 $. Основной результат этой работы состоит в том, что не существует <<просто однородных>>
СПВ-гиперповерх\-ностей, отличных от трубок с аффинно-однородными основаниями. В силу теоремы 2.4, ситуация с новыми индефинитными
<<просто однородными>> гиперповерхностями в $ \Bbb C^3 $, не сводимыми к трубкам, оказалась интереснее
СПВ-случая.

  Используя теоремы 2.2, 2.3 и удаляя из списка теоремы 2.1 поверхности, трубки над которыми вырождены по Леви или обладают
<<богатыми>> алгебрами симметрий, можно получить следующий вывод из теоремы 2.4.

{\bf ТЕОРЕМА 2.5.(Полный список <<просто однородных>> Леви-невырож\-денных гиперповерхностей в $ \Bbb C^3 $)} {\it Любая <<просто однородная>>
невырожденная по Леви гиперповерхность в $ \Bbb C^3 $ локально голоморфно эквивалентна одному из следующих попарно не эквивалентных многообразий:

1) $ \ v= x_1^{\alpha} x_2^{\beta}, \ 0 < |\alpha|\le |\beta|\le 1, \alpha + \beta \ne 1, $

2) $ \  v = (x_1^2 + x_2^2)^{\alpha} e^{\beta arg (x_1 + ix_2)},  \ \alpha \ne 1/2, \beta \ge 0, (\alpha,\beta) \ne (1,0) $,

3) $  \ v = x_1( \alpha \ln x_1 + \ln x_2 ), \ \alpha \notin \{ -1, 0 \} $,

4) \ $ ( v - x_1 x_2 + x_1^3/3 )^2 = \alpha (x_2 - x_1^2/2)^3, \ \notin \{-8/9, 0 \} $,

5) $ \ x_1 v = x_2^2 \pm x_1^{\alpha}, \ \alpha\notin \{0,1,2\} $,

6) $ \  x_1 v = x_2^2 \pm x_1^2 \ln x_1 $,

$
   7) \  v\, ( 1 \pm y_2 x_2) = y_1 y_2, $

$  8) \ (v - x_2 y_1)^2 + y_1^2 y_2^2 = y_1.
$
}

{\bf Замечание 1.} Семейства поверхностей 1)--- 6) из этой теоремы содержат как строго псевдо-выпуклые, так и индефинитные однородные
поверхности.

{\bf Замечание 2.} Выписанные в этой теореме ограничения на параметры, будут прокомментированы в заключительном разделе статьи после содержательных доказательств теоремы 2.5 и сопутствующих ей утверждений. При этом
 все следующие разделы статьи фактически посвящены доказательству Теорем 2.4 и 2.5, так как (достаточно объемные)
 доказательства теорем 2.1 -- 2.3 уже опубликованы (см. [11],[3],[4]).

\

{\bf 3. Схема изучения <<просто однородных>> гиперповерхностей в $ \Bbb C^3 $}

\

   С учетом предыдущих обсуждений для полного решения поставленной задачи
необходимо изучить случай невырожденных по Леви однородных гиперповерхностей в $ \Bbb C^3 $,
алгебры Ли голоморфных векторных полей на которых являются в точности 5-мерными.
   Предлагаемая для этого техническая идея заключается в
описании однородных поверхностей,
ассоциированных с 5-мерными алгебрами Ли, т.е. таких гиперповерхностей,
на которых имеются 5-мерные алгебры Ли голоморфных векторных полей (не обязательно совпадающие с их полными алгебрами симметрий). Затем, с учетом
обсуждений предыдущих разделов, нужно будет исключить из полученного списка кратно транзитивные поверхности.

  Описанию однородных поверхностей, ассоциированных с одной из 5-мерных алгебр,
посвящена статья [12]. Подчеркнем, что в ней изучена лишь отдельная алгебра из всего семейства 5-мерных алгебр Ли.
Совокупность же орбит алгебр голоморфных векторных полей в $ \Bbb C^3 $, имеющих одну и ту же фиксированную алгебраическую
структуру, содержит однопараметрическое семейство голоморфно различных однородных гиперповерхностей.

   При этом, как известно, все семейство 5-мерных алгебр Ли является весьма обширным. Например, классификация Мубаракзянова
включает в себя 67 типов (см. [5,40]) абстрактных алгебр (и содержит, в частности, двух- и трех-параметрические подсемейства алгебр).
При сопоставлении таких результатов работ [5] и [12] может создаться
впечатление о фактической безнадежности применения техники статьи [12] и обозначенной выше идеи в задаче описания всех однородных гиперповерхностей в $ \Bbb C^3 $, а также о необходимости других подходов к этой задаче.

  Однако в серии совместных публикаций [6 - 10] автора настоящей статьи, Акопян Р.С., Атанова А.В. и Коссовского И.Г.
методы работы [12] были применены к описанию большого подсемейства 5-мерных алгебр Ли. Более точно, из семи блоков алгебр,
представляющих полную классификацию Мубаракзянова [5] абстрактных 5-мерных алгебр Ли, были полностью изучены пять блоков.

Описанию однородных поверхностей, ассоциированных с алгебрами из двух оставшихся блоков, посвящены следующие разделы настоящей статьи.
Здесь же мы кратко опишем классификацию Мубаракзянова 5-мерных алгебр Ли, а также некоторые идеи и фрагменты схемы, позволившей изучить пять упомянутых блоков из этой классификации и получить описания соответствующих однородных гиперповерхностей.

\

{\bf 3.1. Классификация Мубаракзянова 5-мерных алгебр Ли}

\

   В работах [5,40] выделены следующие семь блоков 5-мерных алгебр Ли, составляющих их полную классификацию:

I. Разложимые алгебры Ли (27 типов алгебр);

II. Нильпотентные неразложимые алгебры Ли (6 типов алгебр);

III. Разрешимые алгебры Ли с одним ненильпотентным базисным элементом,
содержащие абелев идеал $ 4 g_1 $ (12 типов алгебр);

IV. Алгебры Ли с одним ненильпотентным базисным элементом,
содержащие идеал $ g_{3,1} + g_1 $ (11 типов алгебр);

V. Разрешимые алгебры Ли с одним ненильпотентным базисным элементом,
содержащие идеал $ g_{4,1} $ (3 типа алгебр);

VI. Разрешимые алгебры Ли с двумя линейно ниль-независимыми элементами (7 типов алгебр);

VII. Неразрешимая неразложимая алгебра Ли $ g_5 $.

  Как отмечалось выше, общее число различных типов абстрактных алгебр, содержащихся в обозначенных блоках, составляет 67.

  Каждая из алгебр Ли описывается, как известно, своими коммутационными соотношениями. Ниже приведено описание
такими соотношениями алгебр из блока V.

$$
\begin{array}{|c|c|c|c|c|c|c|c|c|c|c|c|}
\hline
{\mbox{ Алгебры }} & [e_1,e_5] & [e_2,e_3] & [e_2,e_4] & [e_2,e_5] & [e_3,e_4] & [e_3,e_5] & [e_4,e_5] \\ \hline
\mathfrak{g}_{5,30}&    (2 +h)e_1 & \, &  e_1 & ( 1+h) e_2 & e_2  & h e_3 & e_4 \\ \hline
\mathfrak{g}_{5,31}&    3 e_1 &  \,& e_1  &  2 e_2 &  e_2  & e_3  & e_3 + e_4 \\ \hline
\mathfrak{g}_{5,32}&    e_1 & \,  & e_1  & e_2 & e_2 & h e_1 + e_3 &  \, \\ \hline
\end{array}
$$

{\bf Таблица 1.} Коммутационные соотношения для алгебр из блока V

{\bf Замечание.}  Параметр $ h $ в семействе алгебр $ \mathfrak{g}_{5,30} $ может принимать любое
  вещественное значение. На обозначенный тем же символом параметр семейства алгебр $ \mathfrak{g}_{5,32} $
  никакие условия в [5] не оговариваются. Однако
  несложно убедиться, что (формально) однопараметрическое семейство таких алгебр сводится по существу к трем
отдельным алгебрам, отвечающим значениям параметра $ h $ из множества $ \{-1,0,1\} $. На аналогичное уточнение
информации о семействе алгебр $ \mathfrak{g}_{5,30} $, содержащееся в книге [38], автору указал Горбацевич В.В.

  Напомним, что в 5-мерной алгебре имеется, вообще говоря, 10 различных коммутационных соотношений. В приведенную таблицу 1 не включены
(c целью ее упрощенного представления) тривиальные соотношения
$$
        [e_1,e_2] = [e_1,e_3] = [e_1,e_4] = 0,
$$
выполняющиеся для всех трех алгебр из блока V.

  Укажем также на следующее важное свойство обсуждаемой классификации
и используемых в ней базисов: в каждом из 67 типов алгебр имеется не более семи нетривиальных коммутационных соотношений. Максимальное
их количество (равное именно 7) достигается для неразрешимой неразложимой алгебры Ли $ \mathfrak g_5 $; <<таблица умножения>> для этой алгебры приведена ниже.
$$
\begin{array}{|c|c|c|c|c|c|c|c|c|c|c|c|}
\hline
 [e_1,e_2] & [e_1,e_3] & [e_1,e_4] & [e_1,e_5] & [e_2,e_3] & [e_2,e_4] & [e_2,e_5] & [e_3,e_4] & [e_3,e_5] & [e_4,e_5] \\ \hline
 2e_1 & -e_2 &  e_5 & \, & 2 e_3 & e_4 & -e_5 & \, & e_4 & \, \\ \hline
\end{array}
$$

{\bf Таблица 2.} Описание неразрешимой неразложимой 5-мерной алгебры Ли $ \mathfrak g_5 $

  Выписывать здесь коммутационные соотношения для всех остальных алгебр из списка [5] мы не будем. В следующих разделах статьи
соответствующие таблицы будут приведены лишь для обсуждаемых подробно блоков IV и VI классификации Мубаракзянова.

\

{\bf 3.2. Совместное упрощение вида векторных полей}

\

При рассмотрении реализаций произвольной 5-мерной (вещественной) алгебры Ли с базисом $e_1$, $e_2$, $e_3$, $e_4$, $e_5$ в виде алгебры ростков голоморфных векторных полей в $\mathbb{C}^3$ будем записывать каждый элемент этого базиса в виде
$$
\label{Bazis}
e_k = f_k(z_1,z_2,w)\frac{\partial}{\partial z_1} + g_k(z_1,z_2,w)\frac{\partial}{\partial z_2} + h_k(z_1,z_2,w)\frac{\partial}{\partial w} \quad (k = 1,...,5)
\eqno (3.1)
$$
с голоморфными (вблизи обсуждаемой точки поверхности) функциональными коэффициентами $f_k$, $g_k$, $h_k$.

  Отметим, что здесь и везде ниже вещественные векторные поля
$
    2 Re\,(f \frac{\partial}{\partial z_1} + g \frac{\partial}{\partial z_2} + h\frac{\partial}{\partial w})
$
представляются своими (1,0)-компонентами.
  Мы будем использовать также записи вида
$$
e_k = (f_k, g_k, h_k)
$$
для сокращения формулы (3.1).

Через $z_1$, $z_2$, $w$ здесь и далее обозначаются координаты в пространстве $\mathbb{C}^3$. Их вещественные и мнимые части будем обозначать через $x_k = Re\, z_k$, $y_k = Im\, z_k$ $(k = 1, 2)$, $u = Re\,w$, $v = Im\, w$. Пару координат $(z_1, z_2)$ мы будем объединять
в двумерный комплексный вектор~$z$.

  Ниже приводятся несколько лемм, сводящих обсуждение произвольных 5-мерных алгебр голоморфных векторных полей в $ \Bbb C^3 $
к рассмотрению алгебр, базисы которых имеют достаточно простой вид. Именно эти леммы позволили получить описания однородных
гиперповерхностей в случае принадлежности соответствующих им алгебр Ли упомянутым пяти блокам из классификации Мубаракзянова.
Эти же леммы являются основой при рассмотрении двух последних блоков в следующих разделах настоящей статьи. Доказательства этих лемм
можно найти в [10].

  Одним из основных вопросов при обсуждении абстрактной 5-мерной алгебры Ли в связи с задачей об однородности можно считать вопрос о
существовании однородных невырожденных по Леви поверхностей, ассоциированных с такой алгеброй. Важным утверждением, часто
дающим ответ на этот вопрос без детальных исследований голоморфных векторных полей, является

ЛЕММА 1 ({Коссовский}).
Пусть в 5-мерной алгебре $g(M)$ голоморфных векторных полей на вещественной гиперповерхности $M \subset \mathbb{C}^3$ имеется 4-мерная абелева подалгебра. Тогда $M$ вырождена по Леви.

  Для алгебр Ли, не удовлетворяющих условиям леммы 1, приходится разрабатывать более детальную программу их исследования.

  Одна из промежуточных задач, решаемых в процессе голоморфной реализации конкретной алгебры Ли, --- получить за счет голоморфных
преобразований, по возможности, <<простой>> вид функциональных коэффициентов $f_k$, $g_k$, $h_k$ для всей базисной пятерки полей.
Для алгебр векторных полей на однородных поверхностях оказываются справедливыми следующие свойства, связанные с одновременным
упрощением или выпрямлением нескольких базисных полей.

ЛЕММА 2.
Если на невырожденной по Леви гиперповерхности $M \subset \mathbb{C}^3$ имеется пара коммутирующих голоморфных векторных полей $e_1$ и $e_2$, линейно независимых над $\mathbb{R}$, то голоморфной заменой координат эта пара может быть выпрямлена, т.е. приведена (вблизи некоторой точки поверхности) к виду
$$
e_1 = (0,0,1), \, e_2 = (1,0,0).
$$

ЛЕММА 3.
Пусть $M \subset \mathbb{C}^3$ --- вещественно-аналитическая гиперповерхность, невырожденная по Леви, $g(M)$ --- 5-мерная алгебра касательных к $M$ голоморфных векторных полей, имеющая полный ранг (вблизи некоторой точки поверхности $M$), $g^*(M) = \langle e_i, e_j, e_k \rangle$ --- 3-мерная абелева подалгебра в $g(M)$. Тогда голоморфным преобразованием базис этой подалгебры можно привести (вблизи некоторой точки поверхности $M$) к одному из двух видов:
$$
1) \,  e_i =  (1,0,0), \, e_j =  (0,1,0), \, e_k =  (0,0,1),
$$
либо
$$
2) \, e_i =  (1,0,0), \, e_j =  \left(f_j(z_2),0,  h_j(z_2)\right), \, e_k = (0,0,1).
$$

{\bf Замечание.} Далее в статье часто будут появляться голоморфные функции, зависящие лишь от одной переменной $ z_2 $ трехмерного комплексного пространства. Будем выделять такое свойство функции дополнительным символом \ $ \hat{} $\,. Например, во втором случае
предыдущей леммы поле $ e_j $ имеет в краткой записи вид
$
  e_j =  (\hat f_j,0,  \hat h_j).
$

ЛЕММА 4 (о четырех нулях).
Если четверка базисных голоморфных полей 5-мерной алгебры $g(M)$ полного ранга имеет вблизи некоторой точки $M$  вид
$$
\label{4_nulya}
\begin{array}{l}
  e_1 = \left(f_1(z_1,z_2,z_3), 0, h_1(z_1,z_2,z_3)\right),\\
  e_2 = \left(f_2(z_1,z_2,z_3), 0, h_2(z_1,z_2,z_3)\right), \\
  e_3 = \left(f_3(z_1,z_2,z_3), 0, h_3(z_1,z_2,z_3)\right), \\
  e_4 = \left(f_4(z_1,z_2,z_3), 0, h_4(z_1,z_2,z_3)\right),
\end{array}
$$
то поверхность $M$ является вырожденной по Леви (вблизи обсуждаемой точки).

  Использование этих лемм позволяет привести базис любой из обсуждаемых алгебр к виду, удобному для интегрирования.
Многие из интегральных поверхностей, возникающих при таком интегрировании, оказываются вырожденными по Леви. Еще одна
значительная их часть совпадает с уже известными однородными поверхностями, допускающими более чем 5-мерные алгебры симметрий.
  И лишь малая часть таких интегральных поверхностей претендует на <<новый статус>> в задаче об однородности.

Поиск новых однородных объектов является
  главной задачей, решаемой в рамках проблемы полного описания однородных гиперповерхностей в $ \Bbb C^3 $. Вместе с тем,
нужны в этой задаче и методы проверки возможной голоморфной эквивалентности <<разных>> однородных поверхностей.
В частности, является важным факт сводимости обсуждаемых (однородных) поверхностей
к трубчатым многообразиям. Приведем здесь одно утверждение такого рода.

ЛЕММА 5. Пусть 5-мерная алгебра Ли $ \mathfrak g $ имеет две 3-мерных абелевых подалгебры  $ \mathfrak h_1, \mathfrak h_2 $ таких, что
$ \dim(\mathfrak h_1 \cap  \mathfrak h_2) = 2 $, а
голоморфная реализация в $ \Bbb C^3 $ алгебры $ \mathfrak g $ имеет полный ранг. Тогда
 любая невырожденная по Леви орбита этой реализации голоморфно эквивалентна трубчатой гиперповерхности.

   Доказательство.
   По лемме 3 тройку базисных полей $ e_1, e_2, e_4 $ абелевой подалгебры $ \mathfrak h_1 $ можно упростить до одного из двух состояний:

$$
\begin{array}{ll}
  \mbox{Случай 1:} & e_1 = \left(1, 0, 0\right), \\
  \, & e_2 = \left(0, 1, 0\right), \\
  \, & e_3 = \left(f_3(z_1, z_2, w), g_3(z_1, z_2, w), h_3(z_1, z_2, w)\right), \\
  \, & e_4 = \left(0, 0, 1\right), \\
  \, & e_5 = \left(f_5(z_1, z_2, w), g_5(z_1, z_2, w), h_5(z_1, z_2, w)\right).
\end{array}
\eqno (3.2)
$$
$$
\begin{array}{ll}
  \mbox{Случай 2:} & e_1 = \left(1, 0, 0\right), \\
  \, & e_2 = \left(\hat f_2(z_2), 0, \hat h_2(z_2)\right), \\
  \, & e_3 = \left(f_3(z_1, z_2, w), g_3(z_1, z_2, w), h_3(z_1, z_2, w)\right), \\
  \, & e_4 = \left(0, 0, 1\right), \\
  \, & e_5 = \left(f_5(z_1, z_2, w), g_5(z_1, z_2, w), h_5(z_1, z_2, w)\right).
\end{array}
\eqno (3.3)
$$

  В первом случае наличие трех полей $ e_1, e_2, e_4 $, порождающих сдвиги, автоматически обеспечивает трубчатость любой орбиты алгебры
$ \mathfrak g $.

  Во втором случае рассмотрим четверку полей $ e_1, e_2, e_4, e_3 $, три первых из которых образуют базис
первой абелевой подалгебры $ \mathfrak h_1 $, а базисом второй абелевой подалгебры $ \mathfrak h_2 $ является набор $ e_1, e_3, e_4 $.

  Так как поле $ e_3 $ коммутирует с $ e_1 = (1,0,0) $ и $ e_4 = (0,0,1) $, то в координатной записи оно имеет вид
$ e_3 = (\hat f_3, \hat g_3, \hat h_3) $ с компонентами, зависящими только от переменной $ z_2 $.

   При этом случай тождественно нулевой (на невырожденной орбите алгебры $ \mathfrak g $) функции $ \hat g_3(z_2) $
невозможен по лемме 4 о четырех нулях. Если же
$ \hat g_3(z_2) $ отлично от нуля в некоторой точке поверхности, то голоморфной заменой переменной $ z_2 $ этот
функциональный коэффициент можно превратить в тождественную единицу в окрестности исходной точки. После этого
производится еще одна голоморфная замена
$$
    z_1^* = z_1 + \varphi(z_2), \quad z_2^* = z_2, \quad w^* = w + \psi(z_2).
$$
с функциями $ \varphi(z_2),\psi(z_2) $, удовлетворяющими условиям
$$
     \hat f_3 + \varphi'(z_2) = 0, \quad \hat h_3 + \psi'(z_2) = 0.
$$
Она позволяет (см. [12] или [10]) сохранить выпрямленный вид полей $ e_1 = (1,0,0), \  e_4 = (0,0,1) $ и привести поле $ e_3 $ к виду
$
   e_3 = (0,1,0).
$

  Наличие трех выпрямленных полей $ e_1, e_3, e_4 $ в голоморфной реализации алгебры $ \mathfrak g $ снова
гарантирует трубчатость (в некоторых координатах) любой невырожденной орбиты такой алгебры.
  Лемма 5 доказана.
\hfill $\Box$

\

{\bf 3.3. Изученные блоки 5-мерных алгебр Ли}

\

{\bf 3.3.1}   Первым из семи выписанных блоков 5-мерных алгебр Ли в задаче об однородности был изучен
в [6] блок II нильпотентных
неразложимых алгебр (дополненный тремя разложимыми нильпотентными алгебрами Ли). Из всех элементов этого блока
лишь две алгебры допускают
реализацию в виде алгебр голоморфных векторных полей на однородных невырожденных по Леви гиперповерхностях.
С точностью до локальной голоморфной эквивалентности таких поверхностей имеется,
согласно [6], ровно две: это известные квадрики (1.4), т.е
$$
    v = |z_1|^2 + |z_2|^2, \quad
    v = |z_1|^2 - |z_2|^2,
$$
имеющие в действительности 15-мерные алгебры симметрий. Рассмотренные 5-мерные нильпотентные алгебры Ли реализуются в
$ \Bbb C^3 $ в виде подалгебр ранга 5 таких 15-мерных алгебр, обеспечивающх однородность этих поверхностей <<самостоятельно>>,
без апелляции к объемлющим алгебрам.

\

{\bf 3.3.2.} Проще всего в рамках изучения однородных гиперповерхностей проводится рассмотрение блока III алгебр, содержащих
4-мерный абелев идеал. Согласно лемме 1 все интегральные поверхности таких алгебр (если они существуют) могут быть только вырожденными по Леви. Никаких новых примеров однородности здесь нет, т.к. все вырожденные однородные гиперповерхности в $ \Bbb C^3 $ описаны в работе
[3].

\

{\bf 3.3.3.} Блок V, содержащий всего три типа алгебр Ли с коммутационными соотношениями, описанными в Таблице 1,
является более интересным.
Во-первых, именно алгебра из семейства $ g_{5,30} $ при $ h = 1 $ и ассоциированные с этой алгеброй однородные поверхности,
изученные в работе [7], послужили стартовой точкой при разработке описанной выше схемы изучения однородности.

Во-вторых, в блоке V имеются однородные поверхности, отличные от известных трубок над аффинно-однородными основаниями. Более точно,
справедлива

{\bf ТЕОРЕМА 3.1 ([7]).}
{\it Любая невырожденная по Леви голоморфно однородная гиперповерхность, ассоциированная с одной из алгебр блока V,
локально голоморфно эквивалентна либо трубчатой поверхности над известным аффинно-однородным основанием,
либо трубке
$$
  v = x_1 x_2 + x_1^3 \ln x_1,
\eqno (3.4)
$$
основание которой не является аффинно-однородным, либо одной из двух поверхностей
$$
    v (1\pm x_2 y_2) = y_1 y_2.
\eqno (3.5)
$$
}

{\bf Замечание 1.} Поверхность (3.4) имеется в списке [4] и, тем самым, не является новой. Однако обнаружение (в рамках
предложенной схемы изучения однородности) трубки, основание которой НЕ является аффинно однородным, подтвердило (на соответствующем этапе исследования задачи об однородности) обоснованность использования и эффективность именно такой схемы исследования.

{\bf Замечание 2.} Поверхности же (3.5), являющиеся орбитами голоморфной реализации алгебры $ \mathfrak{g}_{5,32} $ при $ h = 0 $,
оказались новыми (см. [7]) в обсуждаемой задаче.

\

{\bf 3.3.4.} Блок VII, содержащий всего одну алгебру $ g_5 $, также представляет интерес в задаче об однородности.

{\bf ТЕОРЕМА 3.2 ([8]).}
{\it Базис (любой) голоморфной реализации алгебры $ \mathfrak{g}_5  $, допускающей невырожденные по Леви орбиты,
можно привести локальным голоморфным преобразованием к виду
$$
 \begin{array}{ccccccc}
  e_1 :   (\ \ \ \ 1\ , & 0\ , &\ 0 \ \ \ \ \ ) \\
  e_2 : ( \ 2 z_1 \ , &  - z_2 \ , & w  ) \\
  e_3: ( \-z_1^2 \ , &  z_1 z_2 - w \ , & - z_1 w \ ) \\
  e_4 : ( \ 0 \ , & 1 \ , & z_1 \ ) \\
  e_5: ( \ 0 \ , &  0 \ , & 1\ ).
\end{array}
\eqno (3.6)
$$
Любая невырожденная по Леви орбита этой алгебры локально голоморфно эквивалентна поверхности
$$
        ( v - x_2 y_1)^2 + y_1^2 y_2^2 =  y_1.
\eqno (3.7)
$$
}

   Поверхность (3.7) также является новой в задаче описания голоморфно однородных вещественных гиперповерхностей
   пространства $ \Bbb C^3 $.
В целом, именно поверхности (3.5) и (3.7), отвечающие алгебрам из блоков V и VII, и только эти поверхности
оказались новыми в обсуждаемой задаче! Доказательство отсутствия других новых поверхностей, связанное с
изучением блоков IV и VI, приводится в следующих разделах настоящей статьи.

\

{\bf 3.3.5.} Обсудим, наконец, блок I, самый объемный по числу представителей.
Напомним, каждая 5-мерная алгебра Ли из этого блока раскладывается в прямую сумму нескольких своих подалгебр
(меньших размерностей). За счет
этого многие из коммутационных соотношений, отвечающих таким алгебрам, оказываются тривиальными.
Легко устанавливается (см. [9-10]) следующий факт.

{\bf Предложение 3.3 ([10]).} {\it Тринадцать из 27 типов алгебр Ли, входящих в блок I, содержат 4-мерный абелев идеал.
}

   В силу леммы 1 орбиты таких алгебр могут быть только вырожденными по Леви. Оставшиеся 14 алгебр приходится разбирать более подробно.
Отметим при этом, что в их число входят 10 разрешимых и 4 неразрешимые алгебры.

{\bf Предложение 3.4 ([10]).}
{\it Любая невырожденная по Леви однородная гиперповерхность, ассоциированная с одной из разрешимых разложимых алгебр из блока I,
локально голоморфно эквивалентна трубчатой поверхности либо над известными аффинно-однородными основаниями,
либо над одной из поверхностей
$$
  v \sin y_2 + y_1\cos y_2 = e^{p y_2} \ (p \in \mathbb{R}),
\eqno (3.8)
$$
не являющихся аффинно-однородными.
}

{\bf Замечание.}
Любую поверхность (3.8) можно привести посредством (локального) голоморфного преобразования в аффинно-однородное состояние,
но с потерей свойства трубчатости.
В работе [4], посвящённой изучению голоморфно однородных поверхностей с богатыми алгебрами симметрии, построено несколько семейств
голоморфно однородных невырожденных по Леви <<странных>> трубок. Их основания,
подобно поверхностям (3.4) и (3.8), не являются аффинно-однородными.
Отметим, что в $ \Bbb C^2 $ такая ситуация невозможна (см. [41]).

  Разложимые неразрешимые алгебры из блока I
представляют собой прямые суммы классических трехмерных алгебр Ли с тривиальной и нетривиальной двумерными алгебрами:
$$
   \mathfrak{sl}(2)+ 2\mathfrak{g}_1, \ \mathfrak{su}(2)+ 2\mathfrak{g}_1, \ \mathfrak{sl}(2)+\mathfrak{g}_2, \ \mathfrak{su}(2)+\mathfrak{g}_2.
\eqno (3.9)
$$
 Именно с двумя последними алгебрами из (3.9) связаны наиболее интересные утверждения (см. [9]) об однородных поверхностях,
ассоциированных с 27 представителями блока 1.1 разложимых алгебр Ли.

{\bf Предложение 3.5 [9].} {\it Голоморфно однородные поверхности, отвечающие двум первым алгебрам из набора (3.9), могут быть только вырожденными по Леви.
Невырожденными по Леви однородными поверхностями, ассоциированными с алгебрами
$
   \mathfrak{sl}(2)+\mathfrak{g}_2 $ и $ \mathfrak{su}(2)+\mathfrak{g}_2
$
и не эквивалентными трубчатым поверхностям с аффинно однородными основаняими, являются (с точностью до локальной голоморфной эквивалентности) лишь:
$$
  vx_1 + x_2 = \alpha|z_2|, \ 0 < \alpha \ne 1 \ (\mbox{алгебра } \mathfrak{sl}(2)+\mathfrak{g}_2),
\eqno (3.10)
$$
и
$$
   v \ch x_1 - x_2 \sh x_1 = \alpha|z_2|, \ \alpha > 0 \ (\mbox{алгебра } \mathfrak{su}(2)+\mathfrak{g}_2).
\eqno (3.11)
$$
}

  Визуально поверхности семейств (3.10) и (3.11) отличаются от всех известных ранее однородных гиперповерхностей. В связи с этим
в течение некоторого времени они представлялись новыми автору настоящей статьи. Вместе с тем имеют место следующие два утверждения.

{\bf Предложение 3.6.} {\it При $ 0 < \alpha \ne 1 $ поверхность (3.10), т.е.
$
    v y_1 + x_2 = \alpha |z_2|,
$
отвечающая 5-мерной алгебре $ sl(2) + g_2$, голоморфно эквивалентна поверхности
$$
   (1 - |z_1|^2 - |z_2|^2 + |w|^2) = \alpha|1-z_1^2 - z_2^2 + w^2|
\eqno (3.12)
$$
с тем же значением параметра $ \alpha $.
}

 Для доказательства предложения 3.6 достаточно проверить непосредственными вычислениями, что голоморфное отображение
$$
    z_1^* = \frac{1+z_2}{w-z_1},\quad z_2^* = \frac{1-z_1^2 - z_2^2 + w^2}{(w-z_1)^2}, \quad w^* = \frac{2(1 - z_2)}{w-z_1}
\eqno (3.13)
$$
переводит обсуждаемую поверхность (3.10) в поверхность (3.12) с тем же $ \alpha $.
\hfill $\Box$

{\bf Замечание 1.} Известное (см. [23]) семейство (3.12) обобщает однородные поверхности Картана [2] из пространства $ \Bbb C^2 $.
Каждая поверхность этого семейства при $ 0 < \alpha \ne 1 $ имеет
6-мерную алгебру симметрий.

{\bf Замечание 2.} Еще один интересный момент связан с (локальной) голоморфной эквивалентностью каждой из Картановых поверхностей (3.12),
являющихся алгебраическими и имеющими четвертую степень,
некоторой поверхности второй степени. Такая эквивалентность является следствием предложения 3.6. Квадратичное отображение
$ z_2^*= z_2^2 $ превращает уравнение
$
 v x_1 + x_2 = \alpha |z_2|
$
в
$$
    v x_1 + x_2^2 - y_2^2 = \alpha |z_2|^2.
$$

  Отметим, что отображение (3.13) можно получить, учитывая вложение алгебры
$  \mathfrak{sl}(2)+\mathfrak{g}_2 $ в алгебру
$
   \mathfrak{sl}(2)+ \mathfrak{sl}(2) \cong \mathfrak{so}(2,2),
$
ассоциированную с однородными поверхностями (3.12). Аналогичное вложение алгебры
$  \mathfrak{su}(2)+\mathfrak{g}_2 $ в алгебру
$
   \mathfrak{su}(2)+ \mathfrak{sl}(2) \cong \mathfrak{so}^*(4)
$
позволяет получить

{\bf Предложение 3.7.} {\it Интегральные поверхности (3.11) неразрешимой алгебры $ m_{17} =\mathfrak{su}(2) + \mathfrak g_2 $
имеют 6-мерные алгебры касательных векторных полей.
}

Доказательство.
   Пятерку базисных полей
алгебры  $ m_{17} =\mathfrak{su}(2) + \mathfrak g_2 $ можно преобразовать за счет голоморфных замен координат к виду
$$
 \begin{array}{ccccccc}
  e_1 :   ( \ i \cos z_1  , & -i z_2 \sin z_1 \ , &\  z_2 \cos z_1 ) \\
  e_2 : (\ - i \sin z_1  , &  -i z_2 \cos z_1 \ , &\  z_2 \sin z_1) ) \\
  e_3: ( \ \  1 \ , & 0, & 0 \ ) \\
  e_4 : ( \ \ \ 0\ , & 0\ , &\ 1 \ \  ) \\
  e_5: ( 0 \ , &  z_2 \ , & w\ )
\end{array}.
\eqno (3.14)
$$

   Поверхности (3.11) получены в [9] интегрированием алгебр с такими базисами.
При этом пятерка полей (3.14) легко дополняется полем
$$
    e_6 = (2i z_2,\  2z_2 w,\  z_2^2 + w^2),
$$
также касательным к поверхностям (3.11). А вещественная линейная оболочка расширенного набора полей $ e_1,..., e_6 $ образует 6-мерную алгебру векторных полей в $ \Bbb C^3 $,
изоморфную $ \mathfrak{su}(2) + \mathfrak{sl}(2) $.
\hfill $\Box $

{\it Следствие.} Поверхности (3.11) НЕ являются новыми в задаче об однородности.

{\bf Замечание 1.} В заметке [9] поверхности (3.11) названы новыми в обсуждаемой задаче. Предложение 3.7 показывает ошибочность такого утверждения. К этой ошибке привели недостаточно выверенные компьютерные вычисления коэффициентов нормальных уравнений поверхностей (3.11).
Более тщательная проверка таких коэффициентов подтверждает вывод предложения 3.7.

{\bf Замечание 2.} Компьютерные вычисления с использованием пакетов символьной математики находят все большее применение в задачах
<<чистой>> математики и, в частности, в близких к обсуждаемому рассмотрениях (см., например, [42-43]). К сожалению, ошибки в этом
относительно новом виде математической активности пока случаются достаточно часто. Но одновременно растет количество математических утверждений, которые можно легко подтвердить или опровергнуть в рамках такой активности.

  Представляет интерес
вопрос об известных однородных гиперповерхностях, которым эквивалентны поверхности (3.11). Однако он подобен большой серии
аналогичных вопросов обо всех полученных в статье орбитах 5-мерных алгебр Ли, которые мы рассматриваем лишь <<по модулю>> возможной
их новизны. А потому детально мы здесь на нем не останавливаемся.

   Отметим лишь, что
максимальные алгебры симметрий для всех поверхностей семейства (3.11) или отдельных его представителей могут иметь
размерность 6 или 7. 6-мерные алгебры Ли $ \mathfrak{su}(2)+ \mathfrak{sl}(2) $ имеют кватернионные модели (2.3) из раздела 2,
и следовательно, возможна эквивалентность поверхностей семейства (3.11) именно этим однородным моделям.

   В то же время
7-мерная алгебра Ли, отвечающая семейству трубок
$$
    u = \alpha \ln(1+ e^{2 x_1})+ \ln x_2, \ \alpha \in \Bbb R \setminus\{-1,0\}, \ \alpha \sim 1/\alpha,
\eqno (3.15)
$$
имеет структуру $ \mathfrak{su}(2)+ \mathfrak{sl}(2) + \Bbb R $ (см. [4]) и содержит в качестве подалгебры
сумму двух первых своих слагаемых. Тем самым, требует проверки возможность голоморфной эквивалентности поверхностей (3.11)
и трубок (3.15).

\

{\bf 4. Однородные поверхности для блока IV из 11 пятимерных алгебр
}

\

   Этот раздел начнем с таблицы коммутационных соотношений для алгебр блока IV из классификации Мубаракзянова.

\begin{equation*}
\begin{array}{|c|c|c|c|c|c|c|c|c|c|c|c|}
\hline
\rotatebox[origin=c]{90}{\mbox{ Алгебры }} & [e_1,e_j], j=2,3,4 & [e_1,e_5] & [e_2,e_3] & [e_2,e_4] & [e_2,e_5] & [e_3,e_4] & [e_3,e_5] & [e_4,e_5] \\ \hline
\mathfrak{g}_{5,19} &  \, & (1+\alpha)e_1 & e_1 & \, &  e_2 & \, & \alpha e_3 & \beta e_4 \\ \hline
\mathfrak{g}_{5,20} &  \, & (1+\alpha)e_1 & e_1 & \, & e_2 & \, & \alpha e_3  &e_1 + (1+ \alpha) e_4 \\ \hline
\mathfrak{g}_{5,21} &  \, & 2 e_1 & e_1 & \, & e_2 + e_3 & \, & e_3 + e_4 & e_4 \\ \hline
\mathfrak{g}_{5,22} &  \, & \, &  e_1 & \, & e_3 & \, & \, &  e_4 \\ \hline
\mathfrak{g}_{5,23} &  \, & 2e_1 & e_1 & \, & e_2 + e_3 & \, & e_3 & \beta e_4 \\ \hline
\mathfrak{g}_{5,24} &  \, & 2e_1 & e_1 & \, & e_2 + e_3 & \, & e_3 & \varepsilon e_1 + 2 e_4 \\ \hline
\mathfrak{g}_{5,25} &  \, & 2pe_1 & e_1 & \, & p e_2 + e_3 & \, & -e_2 +p e_3 & \beta e_4 \\ \hline
\mathfrak{g}_{5,26} &  \, & 2pe_1 & e_1 & \, & p e_2 + e_3 & \, & -e_2 +p e_3 & \varepsilon e_1 + 2 p e_4 \\ \hline
\mathfrak{g}_{5,27} &  \, & e_1 & e_1 & \, & \, & \, & e_3 + e_4 & e_1 + e_4 \\ \hline
\mathfrak{g}_{5,28} &  \, & (1+\alpha)e_1 & e_1 & \, & \alpha e_2 & \, & e_3 + e_4 &  e_4 \\ \hline
\mathfrak{g}_{5,29} &  \, & e_1 & e_1 & \, & e_2 & \, & e_4 & \, \\ \hline
\end{array}
\end{equation*}

{\bf Таблица 3.}
Алгебры Ли с одним ненильпотентным базисным элементом $ e_5 $,
содержащие идеал $ g_{3,1} + g_1 $.

{\bf Замечание 1.} Параметры семейств алгебр, содержащихся в этой таблице, удовлетворяют, согласно [5] некоторым ограничениям. Так,
параметр $ \beta \ne 0 $
для всех трех типов алгебр $ g_{5,19}, g_{5,23}, g_{5,25} $ (в описаниях которых он присутствует), а для алгебр $ g_{5,24}, g_{5,26} $ имеются в виду ограничения
$ \varepsilon = \pm 1 $. Вещественные значения параметров $ \alpha $ и $ p $ никак не ограничиваются.

  Основным результатом этого раздела является следующее утверждение.

{\bf ТЕОРЕМА 4.1.} {\it Алгебрам Ли $ g_{5,19} - g_{5,29} $ из этого блока соответствуют
лишь вырожденные по Леви однородные гиперповерхности и голоморфные образы невырожденных по Леви трубок с аффинно-однородными основаниями.
}

{\bf Замечание.} В задачах об однородности важное место занимают сферические поверхности (строго псевдо-выпуклые или имеющие индефинитную форму Леви). Встречаются они в качестве орбит и у алгебр Ли из блоков IV и VI. Всякая такая поверхность голоморфно эквивалентна одной из двух трубок $ v = x_1^2 \pm x_2^2 $, основания которых аффинно однородны. Поэтому в формулировке теоремы сферические поверхности отдельно не оговариваются.

\

Для доказательства теоремы 4.1 заметим (см. таблицу 3), что
у каждой из 11 алгебр этого блока имеется два абелевых идеала
$ \mathfrak h_1 = <e_1, e_2, e_4> $ и $ \mathfrak h_2 = < e_1, e_3, e_4 > $. По лемме 5 это означает, что все
Леви-невырожденные орбиты таких алгебр являются трубчатыми гиперповерхностями. Утверждение же об аффинной однородности
оснований всех таких трубок мы докажем в двух случаях, возникающих в той же лемме 5 в связи с выпрямлением базисов двух этих идеалов.

\

{\bf 4.1. Первый случай: выпрямление базиса идеала $ \mathfrak h_1 = <e_1, e_2, e_4> $}

\

  В соответствии с леммой 5 базис любой алгебры из нашего блока в первом случае можно представить в виде
$$
\begin{array}{ll}
  \, &  e_1 = \left(1, 0, 0\right), \\
  \, & e_2 = \left(0, 1, 0\right), \\
  \, & e_3 = \left(f_3(z_1, z_2, w), g_3(z_1, z_2, w), h_3(z_1, z_2, w)\right), \\
  \, & e_4 = \left(0, 0, 1\right), \\
  \, & e_5 = \left(f_5(z_1, z_2, w), g_5(z_1, z_2, w), h_5(z_1, z_2, w)\right).
\end{array}
\eqno (4.1)
$$

  Доказательство аффинной однородности оснований трубчатых поверхностей-орбит алгебр с такими базисами требует детальных обсуждений двух оставшихся   базисных полей $ e_3, e_5 $, дополнительных к подалгебре $  \mathfrak h_1 $.

  Для уточнения вида поля $ e_3 $ воспользуемся коммутационными соотношениями
$$
   [e_1,e_3] = 0, \ [e_2,e_3] = e_1, \ [e_3,e_4] = 0,
\eqno (4.2)
$$
имеющими один и тот же вид для всех алгебр блока и связывающими базис выделенной подалгебры с полем $ e_3 $. Из них
получаем упрощенный вид этого поля
$$
   e_3 = (z_2+ A_3, B_3, C_3)
\eqno (4.3)
$$
с некоторыми комплексными константами $ A_3, B_3, C_3 $.

  После этого в каждой алгебре останется рассмотреть по 4 коммутационных соотношения. Запишем их в обобщенном виде
$$
    [e_1,e_5] = T e_1, \quad [e_2,e_5] = N e_2 + Q e_3, \quad [e_4,e_5] =  R e_1 + S e_4
\eqno (4.4)
$$
и
$$
   [e_3, e_5] = r_2 e_2 + r_3 e_3 + r_4 e_4,
\eqno (4.5)
$$
имея в виду, что значения вещественных коэффициентов $ T,N,Q,R,S, r_2, r_3, r_4 $ для каждой из алгебр определяются из таблицы 3.

  Из трех соотношений (4.4) можно получить упрощенный (хотя и достаточно сложный) вид поля $ e_5 $:
$$
   e_5 = \left(T z_1 + R w + \frac12 Q (z_2+ A_3)^2 + A_5, (N + Q B_3)z_2 + B_5, S w + Q C_3 z_2 + C_5\right)
\eqno (4.6)
$$
с некоторыми комплексными константами $ A_5, B_5, C_5 $.

  С учетом возможной нелинейной зависимости формулы (4.6) от переменной $ z_2 $ дальнейшее обсуждение удобно
разбить еще на два подслучая, связанные с коэффициентом $ Q $. Отметим, что согласно таблице 3,
для всех алгебр из изучаемого блока  $ g_{5,19} - g_{5,29} $ коэффициент $ Q $ может равняться либо нулю, либо единице.

\

Сначала рассмотрим
{\bf подслучай 1.1: $ Q = 0 $},
в котором формулы (4.4) не зависят от поля $ e_3 $.
  Это рассмотрение охватывает 5 из 11 типов алгебр обсуждаемого блока, а именно:
$$
    g_{5,19}, g_{5,20}, g_{5,27}, g_{5,28}, g_{5,29}.
\eqno (4.7)
$$

{\bf Предложение 4.1.} Теорема 4.1 верна для пяти типов алгебр (4.7) в подслучае 1.1.

  Доказательство.
  Коэффициенты из соотношений (4.4) и (4.5) для этих алгебр вынесем в отдельную таблицу, имея в виду, что
для всего набора (4.7) выполняется равенство $ r_2 = 0 $.

$$
\begin{array}{|c|c|c|c|c|c|c|c|c|c|c|c|}
\hline
\rotatebox[origin=c]{90}{\mbox{ Алгебры }} & T & N & R & S & r_3 & r_4  \\ \hline
\mathfrak{g}_{5,19} & 1+ \alpha & 1 & 0 & \beta & \alpha & 0 \\ \hline
\mathfrak{g}_{5,20} & 1+ \alpha & 1 & 1 & 1+ \alpha & \alpha & 0  \\ \hline
\mathfrak{g}_{5,27} & 1 & 0 & 1 & 1 &1  & 1   \\ \hline
\mathfrak{g}_{5,28} &1+ \alpha & \alpha & 0 & 1 & 1  & 1   \\ \hline
\mathfrak{g}_{5,29} & 1 & 1 & 0 & 0 &0  & 1   \\ \hline
\end{array}
$$

{\bf Таблица 4.} Коэффициенты из соотношений (4.4) и (4.5) для пяти типов алгебр (4.7)

\

  В этом подслучае поле
$$
    e_5 = (T z_1 + R w + A_5, N z_2 + B_5, S w + C_5)
\eqno (4.8)
$$
является аффинным (линейным), как и весь набор из пяти базисных полей, принимающий вид
$$
\begin{array}{ll}
\label{m3sl1}
  \  & e_1 = \left(1, 0, 0\right), \\
  \, & e_2 = \left(0, 1, 0\right), \\
  \, & e_3 = \left( z_2 + A_3, B_3, C_3)\right), \\
  \, & e_4 = \left(0, 0, 1\right), \\
  \, & e_5 = \left(T z_1 + R w + A_5, N z_2 + B_5, Sw + C_5\right).
\end{array}
\eqno (4.9)
$$

   При этом за счет сдвига переменной $ z_2 + A_3 \rightarrow z_2 $ и замены полей $ e_3, e_5 $ на их линейные комбинации
$$
   e_3^* = e_3 + \mu_1 e_1 + \mu_2 e_2 + \mu_3 e_3, \quad
   e_5^* = e_5 + \nu_1 e_1 + \nu_2 e_2 + \nu_3 e_3
\eqno (4.10)
$$
с подходящими вещественными коэффициентами $ \mu_k, \nu_k $
можно считать, что в формулах (4.9) $ A_3 = 0 $, а коэффициенты
$$
    B_3, C_3, A_5, B_5, C_5
$$
являются чисто мнимыми.

  Так как коммутант (или производная алгебра) $ [g,g] $) любой из алгебр Ли (4.7), не содержит элемента $ e_5 $,
после замены базиса (4.10) соотношение (4.5) изменяется до вида
$$
   [e_3^*,e_5^*] = r_3 e_3^*  + (r_1^*e_1 + r_2^*e_2 +r_4^*e_4)
\eqno (4.11)
$$
с тем же самым $ r_3 $ и с некоторыми вещественными коэффициентами $ r_k^* $.

  При этом, в силу вещественного характера коэффициентов $ T,R,N,S $ в действительности
$$
   r_1^* =  r_2^* = r_4^* = 0.
\eqno (4.12)
$$

В самом деле, вычисление коммутатора полей (4.3) и (4.8) приводит к формуле
$$
  [e_3^*, e_5^*] = ((T-N) z_2 + (C_3 R - B_5), B_3 N, C_3 S).
$$

   Из линейной по $ z_2 $ части первой компоненты этого поля определяется коэффициент $ r_3 = T-N $ в разложении (4.11).
А из чисто мнимого характера сдвиговых частей $ (C_3 R - B_5), B_3 N, C_3 S  $ коммутатора
$
  [e_3^*, e_5^*]
$
получаем желаемый вывод (4.12).

  Напомним теперь, что однородная поверхность $ M $ с 5-мерной алгеброй симметрий и определяющей функцией $ \Phi $ удовлетворяет системе
  соотношений
$$
    Re\, \left(e_k (\Phi)_{|M}\right) \equiv 0, \ k = 1,...,5.
\eqno (4.13)
$$

  Тривиальный вид тройки полей $ e_1, e_2, e_4 $ означает в силу (4.12), что функция $ \Phi $ в этом
  случае не зависит от трех вещественных
переменных $ x_1 = Re\,z_1, \ x_2 = Re\, z_2, \ u = Re\, w $. Обсуждаемая однородная поверхность может быть задана
в такой ситуации трубчатым уравнением
$$
    \Phi(y_1, y_2, v) = 0,
\eqno (4.14)
$$
удовлетворяющим системе двух тождеств (4.12) при $ k = 3,5 $ (напомним, что $ y_1 = Im\,z_1, y_2 =Im\, z_2, v=Im\,w $).

   Записывая два этих тождества в вещественных координатах, получим с учетом вещественного характера коэффициентов $ M,N,R,S $ и
упрощенных обозначений
$$
     B_3 = i b_3, \ C_3= i c_3, \ A_5 = i a_5, \ B_5 = i b_5, \ C_5 = i c_5
$$
следующую систему:
$$
    y_2 \frac{\partial \Phi}{\partial y_1} +  b_3 \frac{\partial \Phi}{\partial y_2} +  c_3 \frac{\partial \Phi}{\partial v} = 0,
\eqno (4.15)
$$
$$
   (M y_1 + R v + a_5) \frac{\partial \Phi}{\partial y_1} +  (N y_2 + b_5) \frac{\partial \Phi}{\partial y_2} +
           (Sv + c_5) \frac{\partial \Phi}{\partial v} = 0.
$$

  Если какая-либо однородная поверхность имеет 5-мерную алгебру Ли голоморфных векторных полей со структурой одной из алгебр (4.7),
то после описанных преобразований координат определяющая функция этой поверхности удовлетворяет системе (4.15).

   Отметим, что при произвольных значениях коэффициентов в (4.15) эта система может не иметь решений. Но для всякой однородной поверхности
соответствующая ей система такого вида гарантированно имеет решения. Остается заметить, что
система (4.15) для такой поверхности является развернутой записью факта касания
основания трубки (4.14) в пространстве трех вещественных переменных $ y_1,y_2,v $
   двумя вещественными векторными полями
$$
   \hat e_3 =y _2 \frac{\partial }{\partial y_1} +  b_3 \frac{\partial }{\partial y_2} +  c_3 \frac{\partial }{\partial v},
\eqno (4.16)
$$$$
   \hat e_5 = (M y_1 + R v + a_5) \frac{\partial }{\partial y_1} +  (N y_2 + b_5) \frac{\partial }{\partial y_2} +
           (Sv + c_5) \frac{\partial }{\partial v}.
$$

   При этом, в силу равенств (4.11) и (4.12) и вытекающего из них следствия
$
   [\hat e_3, \hat e_5] = r_3 \hat e_3,
$
поля (4.16) образуют алгебру Ли в $ \Bbb R^3_{(y_1,y_2,v)} $. Это означает, что у любой трубчатой гиперповерхности
(4.14) в $ \Bbb C^3 $, голоморфная однородность которой обеспечивается алгеброй голоморфных векторных полей с одной из структур (4.7),
ее основание является орбитой двумерной алгебры аффинных векторных полей в $ \Bbb R^3 $, а потому аффинно однородно.

   Тем самым, предложение 4.1 доказано.
\hfill{$\Box$}

\

  Обсудим теперь {\bf подслучай 1.2: $ Q=1 $}, включающий шесть типов алгебр
$$
    g_{5,21}, g_{5,22}, g_{5,23}, g_{5,24}, g_{5,25}, g_{5,26}.
\eqno (4.17)
$$

{\bf Предложение 4.2.} Теорема 4.1 верна для шести типов алгебр (4.17) в подслучае 1.2.

Доказательство.
Алгебра векторных полей в $ \Bbb C^3 $, реализующая какую-либо алгебру из
набора (4.17), должна содержать в этом случае квадратичное поле $ e_5 $, определяемое формулой (4.6).

  Здесь мы для начала вычислим коммутатор двух полей
$ e_3 $ и $ e_5 $, используя формулы (4.3) и (4.6). Имеем здесь
$$
   [e_3, e_5] =  (z_2 + A_3)(T,0,0) + B_3 (z_2+A_3,N+B_3, C_3) + C_3(R,0,S) - ((N+B_3)z_2 + B_5)(1,0,0).
\eqno (4.18)
$$

  Сравнение этого равенства с формулой (4.4) в сдвиговых частях трех компонент коммутатора приводит к
следующей системе трех уравнений на коэффициенты полей $ e_3 $ и $ e_5 $:
$$
   B_3 A_3 + M A_3 +C_3 R - B_5  = r_3 A_3,
\quad
  B_3 (N+B_3) = r_2 + r_3 B_3,
\quad
  B_3 C_3 + C_3 S = r_3 C_3 + r_4.
\eqno (4.19)
$$

   Для произвольной алгебры из набора (4.17) первое из этих уравнений, очевидно, позволяет выразить
коэффициент $ B_5 $ через остальные коэффициенты полей $ e_3, e_5 $ и параметры самой алгебры:
$$
    B_5 = A_3(B_3 + M - r_3) + C_3 R.
$$

  А два последних равенства из (4.19) перепишем в виде
$$
   B_3( N + B_3 - r_3) = r_2,
\quad
   C_3( S + B_3 - r_3) = r_4
\eqno (4.20)
$$
и рассмотрим их для каждого из шести обсуждаемых типов алгебр с учетом следующей таблицы.

$$
\begin{array}{|c|c|c|c|c|c|c|c|c|c|c|c|}
\hline
\rotatebox[origin=c]{90}{\mbox{ Алгебры }} & M & N & R & S & r_2 & r_3 & r_4 & B_3( N + B_3 - r_3) = r_2 & C_3( S + B_3 - r_3) = r_4 \\ \hline
\mathfrak{g}_{5,21} & 2 & 1 & 0 & 1 & 0 & 1 & 1 & B_3^2 = 0 & B_3 C_3 = 1 \\ \hline
\mathfrak{g}_{5,22} & 0 & 0 & 0 & 1 & 0 & 0 & 0  & B_3^2 = 0 & C_3=0 \\ \hline
\mathfrak{g}_{5,23} & 2 & 1 & 0 & \beta & 0 & 1 & 0  & B_3^2 = 0 & C_3(\beta-1) = 0 \\ \hline
\mathfrak{g}_{5,24} & 2 & 1 & \varepsilon & 2 &  0 & 1 & 0  & B_3^2 = 0 &  C_3=0 \\ \hline
\mathfrak{g}_{5,25} & 2p & p & 0 & \beta  & -1  & p & 0  & B_3^2 = -1 &  C_3(\beta-p\pm i)=0 \\ \hline
\mathfrak{g}_{5,26} & 2p & p & \varepsilon & 2p & -1  & p & 0  & B_3^2 = -1 & C_3(p \pm i)=0 \\ \hline
\end{array}
$$
\

{\bf Таблица 5.} Коэффициенты из соотношений (4.20) для шести типов алгебр (4.17)

\

  Из двух последних столбцов этой таблицы можно получить простые выводы о трех типах алгебр из шести. Так, для алгебры
$
  g_{5,21}
$
два равенства (4.20) очевидно противоречат друг другу, следовательно реализовать эту абстрактную алгебру в виде алгебры
векторных полей на (Леви-невырожденной) однородной вещественной гиперповерхности в $ \Bbb C^3 $ невозможно.

   Еще для двух алгебр
$ g_{5,22} $ и $ g_{5,24} $ из таблицы следует вывод
$$
     B_3 = C_3 = 0.
\eqno (4.21)
$$
Так как в теореме 4.1 основной интерес представляют невырожденные по Леви однородные поверхности, пара
$$
    e_1 = (1,0,0) \mbox{ и } \ e_3 = ( z_2 + A_3, B_3, C_3)
$$
базисных векторных полей алгебры $ g(M) $, касательных к любой
такой поверхности $ M $, должна быть линейно независима над $ \Bbb C $. В случае же выполнения равенств (4.21) эти поля оказываются зависимыми, а однородная гиперповерхность, являющаяся орбитой такой алгебры может быть только вырожденной.

   Аналогично, для семейства алгебр
$
  g_{5,23},
$
зависящего от одного вещественного параметра $ \beta $, лишь при $ \beta = 1 $ можно обсуждать существование голоморфных реализаций (с невырожденными орбитами). При остальных значениях параметра для алгебр из этого семейства выполняются условия (4.21), означающие отсутствие невырожденных орбит.

   Таким образом остается разобраться с возможными голоморфными реализациями семейств алгебр
$
  g_{5,25}, g_{5,26}
$
и отдельной алгебры
$ g_{5,23} $ при $ \beta = 1, C_3 \ne 0 $. Из всего обсуждаемого набора (4.17) только в этих случаях пока остается возможность существования невырожденных по Леви орбит. Доказательство предложения 4.2 завершает следующее
утверждение.

{\bf Предложение 4.3.} {\it Все невырожденные по Леви орбиты голоморфных реализаций любой алгебры из совокупности
$$
    g_{5,23} (\mbox{при }  \beta = 1, C_3 \ne 0), \ g_{5,25}, \ g_{5,26}
\eqno (4.22)
$$
голоморфно эквивалентны трубчатым гиперповерхностям с аффинно однородными основаниями.
}

{\bf Замечание.} Все однородные невырожденные по Леви поверхности из этого предложения были известны ранее. Все они имеют
<<богатые>> алгебры симметрий, размерность которых строго больше 5.
Обсуждаемые в этом разделе пятимерные алгебры, орбитами которых являются эти поверхности, содержатся в качестве подалгебр в некоторых бОльших алгебрах голоморфных векторных полей в $ \Bbb C^3 $, касательных к тем же поверхностям.

  Доказательство предложения 4.3 проведем по-отдельности для всех трех заявленных в нем типов алгебр.

  Во-первых, с учетом сказанного выше, базис голоморфной реализации алгебры $ g_{5,23} $ (при $ \beta = 1, C_3 \ne 0 $) имеет следующий вид
($ B_5 = A_3 $):
$$
\begin{array}{ll}
\label{m3sl1}
  \  & e_1 = \left(1, 0, 0\right), \\
  \, & e_2 = \left(0, 1, 0\right), \\
  \, & e_3 = \left( z_2 + A_3, 0, C_3)\right), \\
  \, & e_4 = \left(0, 0, 1\right), \\
  \, & e_5 = \left(2 z_1 + (1/2)(z_2+A_3)^2 + A_5,  z_2 + A_3, w + C_3 z_2 + C_5\right).
\end{array}
\eqno (4.23)
$$

   Сдвигая переменные
$$
   z_1 + \frac{A_5}2 \rightarrow z_1, \quad z_2 + A_3 \rightarrow z_2, \quad w + C_5 - {A_3}{C_3} \rightarrow w,
$$
освободимся от части свободных констант в формулах (4.23), т. что два <<сложных>> базисных поля в (4.23) примут вид:
$$
    e_3 = (z_2 , 0, C_3), \quad
    e_5 = (2 z_1 + \frac12 z_2^2, z_2, w + C_3 z_2)
\eqno (4.24)
$$

    Еще одной заменой $ w \rightarrow C_3 w $, сохраняющей выражение
$ w \frac{\partial}{\partial w} $, константа $ C_3 $ в этих формулах превращается в единицу.

   Учтем теперь, что производную $ \partial\Phi(y_1, y_2, v)/\partial v $ определяющей функции любой обсуждаемой невырожденной по Леви трубки можно считать отличной от нуля (иначе из тождественного равенства нулю еще и производной $ \partial \Phi/{\partial u} $ следует вырожденность
такой поверхности). Тогда можно задавать искомую поверхность явным уравнением вида
$
  \Phi = - v + F(y_1, y_2) = 0,
$
т. что, соответственно,
$$
    \frac{\partial \Phi}{\partial v} = -1, \quad \frac{\partial \Phi}{\partial y_k} = \frac{\partial F}{\partial y_k}, \ k = 1,2.
$$

 Обозначая еще
$
  Re\, C_3 = c_{31}, \ Im\, C_3 = c_{32},
$
выпишем систему вещественных уравнений в частных производных, отвечающую двум полям (4.24):
$$
    y_2 \frac{\partial F}{\partial y_1} = c_{32},
\qquad
    (2 y_1 + x_2 y_2) \frac{\partial F}{\partial y_1} + y_2 \frac{\partial F}{\partial y_2} = F + c_{31} y_2 + c_{32} x_2.
\eqno (4.25)
$$

   Два слагаемых, содержащие <<лишнюю>> переменную $ x_2 $, можно удалить из второго уравнения за счет вычитания из него первого уравнения,
умноженного на $ x_2 $. Упрощенная система примет вид
$$
    y_2 \frac{\partial F}{\partial y_1} = c_{32},
\qquad
    2 y_1 \frac{\partial F}{\partial y_1} + y_2 \frac{\partial F}{\partial y_2} = F + c_{31} y_2.
\eqno (4.26)
$$

   Подобно обсуждениям системы (4.15), для каждой пары коэффициентов $ (c_{31}, c_{32}) $, входящей в (4.26), легко строится двумерная алгебра аффинных векторных полей в $ \Bbb R^3_{y_1,y_2,v} $, орбитами которой являются поверхности $ v = F(y_1,y_2) $.
Эти аффинно-однородные поверхности и являются основаниями обсуждаемых голоморфно однородных трубок, отвечающих
алгебре $ g_{5,23} $ (при $ \beta = 1, C_3 \ne 0 $).

 Не разворачивая здесь всех вычислений и манипуляций с константами, связанных с
решением (стандартными методами) системы (4.26), зафиксируем лишь итоговые выводы об однородных гиперповерхностях для этой алгебры :

\

  при $ c_{32} = 0, c_{31} \ne 0 $ всем решениям системы (4.26) соответствуют вырожденные по Леви поверхности в $ \Bbb C^3 $;

  при  $ c_{31} = 0, c_{32} \ne 0 $ все орбиты исходной алгебры голоморфно эквивалентны индефинитной квадрике $ v = y_1^2 - y_2^2 $;

  при  $ c_{31} \ne 0, c_{32} \ne 0 $ все орбиты исходной алгебры голоморфно эквивалентны (известной) аффинно-однородной индефинитной
поверхности
  $ v = y_1 y_2 + y_1^2 \ln y_1 $.

   Следовательно, орбиты этой алгебры удовлетворяют утверждению предложения 4.1.

\

  Рассмотрим далее орбиты семейства алгебр $ g_{5,25} $. Базис голоморфной реализации любой такой алгебры в обсуждаемом подслучае 1.1 можно привести к виду
$$
\begin{array}{ll}
\label{m3sl1}
  \  & e_1 = \left(1, 0, 0\right), \\
  \, & e_2 = \left(0, 1, 0\right), \\
  \, & e_3 = \left( z_2 + A_3, i \varepsilon_1, 0\right), \
\varepsilon_1 = \pm 1,\\
  \, & e_4 = \left(0, 0, 1\right), \\
  \, & e_5 = \left(2 p z_1 + (1/2)(z_2 + A_3)^2 + A_5, (p+i\varepsilon_1) (z_2 + A_3), \beta w  + C_5\right).
\end{array}
\eqno (4.27)
$$

  Напомним, что, согласно замечанию к таблице 3, параметр $ \beta $ для этого семейства отличен от нуля.

  Преобразованиями, аналогичными описанным выше, два уравнения в частных производных, отвечающие полям $ e_3, e_5 $,
приводятся к виду ($ a_5 = Im\, A_5 $)
$$
    y_2 \frac{\partial F}{\partial y_1} + \varepsilon_1\frac{\partial F}{\partial y_2} = 0  ,
\qquad
    (2 p y_1 + x_2 y_1 + a_5) \frac{\partial F}{\partial y_1} +(p y_2 + \varepsilon_1 x_2) \frac{\partial F}{\partial y_2} = \beta F.
$$

  Освобождаясь далее от вхождения переменной $ x_2 $ во второе уравнение (так же, как это сделано при
  обсуждении предыдущей алгебры), останется разобрать все решения упрощенной системы
$$
    y_2 \frac{\partial F}{\partial y_1} + \varepsilon_1\frac{\partial F}{\partial y_2} = 0  ,
\qquad
    (2 p y_1 + a_5) \frac{\partial F}{\partial y_1} + p y_2 \frac{\partial F}{\partial y_2} = \beta F.
\eqno (4.28)
$$

  При произвольных значениях параметров $ \beta\ne 0, p, a_5 $ эта система также допускает интерпретацию в терминах двумерной алгебры аффинных векторных полей в $ \Bbb R^3 $. Поэтому утверждение предложения 4.1 верно и в случае семейства алгебр $ g_{5,25} $.

   Приведем конкретные выводы
  относительно аффинно-однородных оснований трубок в этом семействе голоморфно-однородных поверхностей:

 1) при $ a_5 = 0, p = 0 $ трубки над такими основаниями вырождены по Леви;

 2) при $ a_5 \ne 0, p = 0 $ основания этих трубок аффинно эквивалентны поверхностям
$$
     v = y_2^2 \pm \ln y_1;
$$

 3) наконец, при $ p \ne 0 $ получаем поверхности в $ \Bbb R^3 $ с уравнениями
$$
    v = y_2^2 \pm y_1^{\lambda}, \ \lambda = \frac{2p}{\beta}.
$$

\

  Последний фрагмент доказательства Предложения 4.1 связан с орбитами семейства алгебр
$ g_{5,26} $. В соответствии с начальными обсуждениями подслучая 1.1, базис голоморфной реализации любой алгебры
из этого семейства считаем имеющим вид
$$
\begin{array}{ll}
\label{m3sl1}
  \  & e_1 = \left(1, 0, 0\right), \\
  \, & e_2 = \left(0, 1, 0\right), \\
  \, & e_3 = \left( z_2 + A_3, i \varepsilon_1, 0\right), \
\varepsilon_1 = \pm 1, \varepsilon_2 = \pm 1, \\
  \, & e_4 = \left(0, 0, 1\right), \\
  \, & e_5 = \left(2 p z_1 +\varepsilon_2 w + (1/2)(z_2 + A_3)^2 + A_5, (p+i\varepsilon_1) (z_2 + A_3), 2 p w  + C_5\right).
\end{array}
\eqno (4.29)
$$

  Система уравнений в частных производных, отвечающая такому набору базисных полей, как и ранее, описывает трубчатую поверхность. Основание
этой трубки $ v = F(y_1, y_2) $ удовлетворяет системе двух уравнений, в одно из которых, как и в предыдущих рассмотрениях, входит
<<лишняя>> переменная $ x_2 $. Ее удаление описанным выше приемом приводит к системе
$$
    y_2 \frac{\partial F}{\partial y_1} + \varepsilon_1\frac{\partial F}{\partial y_2} = 0  ,
\qquad
    (2 p y_1 + \varepsilon_2 F + a_5) \frac{\partial F}{\partial y_1} + p y_2 \frac{\partial F}{\partial y_2} = 2 p F + c_5,
\eqno (4.30)
$$
где
$ a_5 = Im\, A_5, \ c_5 = Im\, C_5 $.

  В случае $ a_5 = c_5 = p = 0 $ решения этой системы, описывающие основания интересующих нас трубок,
могут приводить лишь к вырожденным по Леви поверхностям в $ \Bbb C^3 $.

    Если $ p =0, c_5 \ne 0 $, то все решения системы (4.30) задают поверхности в  $ \Bbb R^3 $, аффинно эквивалентные двум квадрикам
$$
    v = y_1^2 \pm y_2^2.
$$
     Если же $ p \ne 0 $, то независимо от значений коэффициентов $ a_5, c_5 $ все решения системы сводятся к двум (известным)
аффинно-однородным поверхностям
$$
     v = y_1^2 \pm y_2 \ln y_2.
$$

  Тем самым, предложение 4.3 (а значит, и предложение 4.2) доказано полностью.
\hfill{$\Box$}

  Вместе с этим завершено рассмотрение случая 1 для всего блока алгебр
$ g_{5,19} - g_{5,29} $.

\

{\bf 4.2. Второй случай: выпрямление базиса идеала $ \mathfrak h_2 = <e_1, e_3, e_4> $}

\

   Согласно рассуждениям из доказательства леммы 5, базис голоморфной реализации любой из рассматриваемых алгебр
в этом случае можно привести к виду
$$
\begin{array}{ll}
  \, & e_1 = \left(1, 0, 0\right), \\
  \, & e_2 = (\hat f_2(z_2) , 0, \hat h_2(z_2)), \\
  \, & e_3 = \left(0, 1, 0\right), \\
  \, & e_4 = \left(0, 0, 1\right), \\
  \, & e_5 = \left(f_5(z_1, z_2, w), g_5(z_1, z_2, w), h_5(z_1, z_2, w)\right).
\end{array}
\eqno (4.31)
$$

  Пользуясь коммутационным соотношением
$
   [e_2,e_3] = e_1,
$
общим для всех алгебр из блока IV, получим уточненный вид поля
$$
   e_2 = (-z_2 + A_2, 0, C_2)
\eqno (4.32)
$$
с некоторыми комплексными константами $ A_2, C_2 $.

  Вычислим теперь остальные коммутаторы базисных полей обсуждаемой алгебры Ли. По формулам, уже использованным выше
при обсуждении случая 1, имеем:
$$
    [e_1,e_5] = T e_1, \ [e_3,e_5] = r_2 e_2 + r_3 e_3 + r_4 e_4, \ [e_4,e_5] =  R e_1 + S e_4, \ [e_2,e_5] =  N e_2 + Q e_3.
\eqno (4.33)
$$

  Из трех первых соотношений (4.33) получаем следующий вид поля $ e_5 $:
$$
   e_5 = \left(T z_1+ R w - \frac12\, r_2\, (z_2 - A_2)^2 + A_5, \ r_3 z_2 + B_5, \ (r_2 C_2 + r_4)z_2 +  Sw + C_5 \right).
$$

  Тогда
$$
  [e_2, e_5] =(-T(z_2 -A_2) + r_3 z_2 - B_5, \ 0, \ C_2 S).
$$

  Разворачивая покомпонентно последнее равенство из (4.33) и учитывая формулы для базисных полей изучаемых алгебр, получаем
 систему следующих трех уравнений:
$$
   (r_3 - T ) z_2 + (T  A_2 - B_5) = N(A_2 - z_2), \quad 0 = Q, \quad C_2 S = N C_2.
\eqno (4.34)
$$

   Второе из трех равенств (4.34) существенно сужает наши рассмотрения, т.к. только для 5 типов алгебр из 11, составляющих блок
$ g_{5,19} - g_{5,29} $, имеет
место равенство $ Q = 0 $, а именно, для алгебр, рассмотренных в подслучае 1.1.
  Напомним, что для всех таких алгебр коэффициент $ r_2 $ равен нулю, а потому поле $ e_5 $ оказывается аффинным.

  Третье равенство из (4.34), т.е. $ C_2 (S - N) = 0 $ еще больше сужает возможности для требуемых голоморфных реализаций алгебр из
  обсуждаемого блока. В самом деле, в силу условия невырожденности по Леви обсуждаемых однородных поверхностей, коэффициент $ C_2 $
не должен быть нулевым (иначе векторные поля $ e_1 $ и $ e_2 $ окажутся линейно зависимыми).

Но равенство $ S= N $,
обеспечивающее выполнение третьего условия (4.34), НЕ выполняется для алгебр $ g_{5,27} $ и $ g_{5,29} $. А для трех оставшихся семейств
алгебр оно выполняется только при дополнительных ограничениях на параметры алгебр. Все однородные гиперповерхности, возможные при этих
ограничениях, описываются в следующем утверждении.

{\bf Предложение 4.4.}
{\it Все невырожденные по Леви орбиты голоморфных реализаций любой алгебры из совокупности
$$
    g_{5,19} \ (\mbox{при }  \beta = 1), \ g_{5,20} \ (\mbox{при }  \alpha = 0), \ g_{5,28} \ (\mbox{при } \alpha = 1)
\eqno (4.35)
$$
голоморфно эквивалентны трубчатым гиперповерхностям с аффинно однородными основаниями.
}

  Доказательство.

Трубчатый характер орбит алгебр (4.35) следует из наличия трех
выпрямленных полей $ e_1, e_3, e_4 $ в базисе любой из этих алгебр, т. что необходимо доказать лишь аффинную однородность
оснований всех таких трубок. Уравнения этих оснований определяются двумя оставшимися полями $ e_2, e_5 $ в каждой алгебре. При этом поле
$ e_2 $ имеет вид (4.32), общий для всех обсуждаемых алгебр. А последнее
поле $ e_5 $ выписано ниже для всех трех интересующих нас типов алгебр с учетом полученных уточнений на коэффициенты:
$$
   g_{5,19} \ : \ e_5 = ((1+\alpha) z_1 + A_5, \alpha z_2 + B_5, w + C_5);
$$
$$
   g_{5,20} \ : \ e_5 = ( z_1 + w + A_5, B_5, w + C_5);
$$
$$
   g_{5,28} \ : \ e_5 = (2 z_1 + A_5, z_2 + B_5, z_2 + w + C_5).
$$

     Подчеркнем, что поля $ e_2, e_5 $ являются аффинными, все коэффициенты $ A_2, C_2, A_5, B_5, C_5 $ здесь можно считать чисто мнимыми,
а коммутационное соотношение, связывающее два этих поля в исходных 5-мерных алгебрах имеет очень простой вид
$$
    [e_2,e_5] = N e_2.
$$

   Это означает, что вещественные уравнения в частных производных, связанные с действиями двух этих полей на определяющую функцию
голоморфно однородной гиперповерхности в $ \Bbb C^3 $, описывают орбиты двумерных алгебр аффинных векторных полей в $ \Bbb R^3 $.
Эти орбиты и являются интересующими нас аффинно-однородными основаниями изучаемых трубок.

     Тем самым, предложение 4.4 доказано.
\hfill{$\Box$}

{\bf Замечание.} Конкретный вид возникающих в предложении 4.4 аффинно-однородных поверхностей в ходе доказательства не уточнялся.
В частности, некоторые из таких поверхностей могут оказаться вырожденными по Леви, что не запрещалось приведенным доказательством.

   Для полноты картины дадим краткую информацию об орбитах полученных двумерных алгебр.

 I. У семейства $ g_{5,19} $ основания трубок, возникающих в случае 2, могут быть аффинно эквивалентны (в зависимости от значений параметров, определяющих алгебру):

1) индефинитной квадрике $ v = y_1 y_2 $,

2) поверхностям семейства $ v = y_1 y_2 + y_2^{\alpha}, \ \alpha \in \Bbb R $,

3) поверхности $ v = y_1 y_2 + \ln y_2 $.

II. Все орбиты алгебры $ g_{5,20} $ при $ \alpha = 0 $ аффинно эквивалентны поверхности
$ v = y_1 y_2 + e^{y_2} $.

III. Все орбиты алгебры  $ g_{5,28} $ при $ \alpha = 1 $ аффинно эквивалентны поверхности
$ v = y_1 y_2 + y_2^2 \ln y_2 $.

\

  Завершая этот раздел, кратко прокомментируем аналогичные <<неявные>> результаты Теоремы 4.1, связанные с
невырожденными орбитами алгебр (4.7) в подслучае 1.1.

  Можно показать, что алгебры типов $ \mathfrak g_{5,27}, \mathfrak g_{5,28}, \mathfrak g_{5,29} $
  не реализуются в рамках этого подслучая.

Из всего двухпараметрического семейства алгебр $ \mathfrak g_{5,19} $ с параметрами $ \alpha, \beta $ здесь
реализуются лишь два однопараметрических подсемейства при $ \alpha =1 $ и $ \alpha = \beta $.
  Орбитами этих подсемейств являются трубки, эквивалентные, соответственно, поверхностям
$$
   v = y_1^2 \pm y_2^{2/\beta} \ \quad \mbox{и} \quad
$$
$$
   v = y_1 y_2 + y_1^{\alpha} \ (\alpha \ne -1), \ v = y_1 y_2 + \ln y_1 \ (\alpha =-1)
$$

  Трубчатые орбиты семейства алгебр $ \mathfrak g_{5,20} $ могут быть невырожденными по Леви только при $ \alpha = 1 $. Любая из них
аффинно эквивалентна одной из двух поверхностей
$$
     v =  y_2^2 \pm y_1 \ln y_1.
$$

  Все эти конкретные поверхности получаются за счет интегрирования стандартными методами систем двух квазилинейных
уравнений в частных производных.
  Отметим еще, что все выписанные здесь трубки имеют 6-мерные алгебры голоморфных симметрий.

\

{\bf 5. Однородные поверхности для блока VI из 7 пятимерных алгебр
}

\

  Выпишем, прежде всего, таблицу коммутационных соотношений для семи алгебр из блока VI.

\begin{equation*}
\begin{array}{|c|c|c|c|c|c|c|c|c|c|c|c|}
\hline
\rotatebox[origin=c]{90}{\mbox{ Алгебры }} & [e_1,e_2] & [e_1,e_3] & [e_1,e_4] & [e_1,e_5] & [e_2,e_3] & [e_2,e_4] & [e_2,e_5] & [e_3,e_4] & [e_3,e_5] & [e_4,e_5] \\ \hline
\mathfrak{g}_{5,33} & \, & \, & e_1  & \, & \, & \, & e_2  & \beta e_3 & \gamma e_4 & \, \\ \hline
\mathfrak{g}_{5,34} & \, & \, & \alpha e_1  & e_1& \, & e_2 & \, & e_3 & e_2  & \, \\ \hline
\mathfrak{g}_{5,35} & \, & \, & \beta e_1  & \alpha e_1 & \, & e_2 & - e_3 & e_3  &  e_2 & \, \\ \hline
\mathfrak{g}_{5,36} & \, & \, & e_1  & \, & e_1 & e_2 & -e_2 & \, & e_3 &  \, \\ \hline
\mathfrak{g}_{5,37} & \, & \, & 2 e_1  & \, & e_1 &e_2 & -e_3 & e_3 & e_2 & \, \\ \hline
\mathfrak{g}_{5,38} & \, & \, & e_1  & \, & \, & \, & e_2 & \, & \, & e_3 \\ \hline
\mathfrak{g}_{5,39} & \, & \, & e_1  &  -e_2 & \, & e_2 & \,e_1 &\, & \,  & e_3 \\ \hline
\end{array}
\end{equation*}

{\bf Таблица 6.}
Разрешимые алгебры Ли с двумя линейно ниль-независимыми элементами $ e_4, e_5 $.

\

{\bf Замечание.} Параметр $ \beta $ алгебры $ \mathfrak g_{5,35} $ в [5] обозначен через $ h $. Мы заменили
это обозначение на $ \beta $, т.к. везде в статье символ $ h $ используется для одной из компонент голоморфного
отображения пространства $ \Bbb C^3 $.

   Напомним, что (как и в предыдущих разделах статьи) мы интересуемся, в первую очередь, Леви-невырожденными однородными
гиперповерхностями и зафиксируем два свойства, общие для всех алгебр из блока $ \mathfrak g_{5,33} - \mathfrak g_{5,39} $.

  Во-первых, поля $ e_1 $ и $ e_2 $ (так же, как и $ e_1 $ и $ e_3 $) коммутируют друг с другом в любой из таких алгебр.
Голоморфной заменой координат такую пару полей, касательных к невырожденной гиперповерхности, можно выпрямить
 (в произвольной голоморфной реализации любой из алгебр блока, допускающей Леви-невырожденные орбиты);

  во-вторых, линейная оболочка $ \mathfrak{h} = <e_1,e_2,e_3> $ является подалгеброй в любой алгебре $ \mathfrak g $ из этого блока.
Отметим, что для пяти типов алгебр блока эта подалгебра будет абелевой (для алгебр $ \mathfrak g_{34} $,
$ \mathfrak g_{35} $, $ \mathfrak g_{38} $, $ \mathfrak g_{39} $ даже абелевым идеалом), а для двух <<исключительных>> алгебр
$ \mathfrak g_{5,36}$ и $ \mathfrak g_{5,37} $
это свойство подалгебры $ <e_1,e_2,e_3> $ (важное в наших дальнейших обсуждениях) не выполняется.

  Теперь начнем конкретные обсуждения алгебр из блока, опираясь на эти свойства. Для алгебр с абелевой подалгеброй
$ <e_1,e_2,e_3> $ наши обсуждения аналогичны предыдущему разделу статьи.
Они заключаются в рассмотрении {\bf двух случаев, в первом из которых} тройка полей $ e_1,e_2,e_3 $ считается выпрямленной.

   Во втором случае, отвечающем лемме 3 из раздела 3, выпрямляются два из трех базисных полей абелевой подалгебры, а третье лишь частично упрощается. Отметим, что выбор этого третьего поля (связанный в каждом отдельном случае с коммутационными соотношениями в исследуемой
алгебре) позволяет существенно упрощать обсуждения. Поэтому вид базисной тройки во втором случае будет специально оговариваться
для каждого набора алгебр.

\

   Итак, рассмотрим сразу пять типов алгебр
$$
    \mathfrak g_{5,33}, \mathfrak g_{5,34},\mathfrak g_{5,35},\mathfrak g_{5,38},\mathfrak g_{5,39},
\eqno (5.1)
$$
считая, что для каждой из них тройка базисных векторных полей имеет вид
$$
\begin{array}{ll}
  \, & e_1 = \left(1, 0, 0\right), \\
  \, & e_2 = \left(0, 1, 0\right), \\
  \, & e_3 = \left(0, 0, 1\right).
\end{array}
\eqno (5.2)
$$

   В последующих рассмотрениях эти пять типов алгебр мы разделим на еще более мелкие блоки. Отдельно будут
рассмотрены пара алгебр $ \mathfrak g_{5,38}$ и $ \mathfrak g_{5,39} $ и блок из трех семейств
$ \mathfrak g_{5,33}, \mathfrak g_{5,34}, \mathfrak g_{5,35}$ (при этом семейство
$ \mathfrak g_{5,35}$ из последнего блока потребует дополнительных обсуждений).

\

{\bf 5.1. Орбиты алгебр $ \mathfrak g_{5,38}$ и $ \mathfrak g_{5,39} $}

\

  {\bf Предложение 5.1.} Алгебры $ \mathfrak g_{5,38}$ и $ \mathfrak g_{5,39} $ не допускают голоморфных реализаций, для которых
  выполняются условия (5.2).

  Доказательство.
  В предположении справедливости формул (5.2) найдем координатные представления полей $ e_4 $ и $ e_5 $. Подробные выкладки
приведем для второй алгебры из этой пары, т.е. для $ \mathfrak g_{5,39} $. Для первой алгебры аналогичные вычисления оказываются несколько проще чем для второй.

  Согласно коммутационным соотношениям для алгебры $ \mathfrak g_{5,39} $, имеем
$$
   [e_1, e_4] = e_1, \  [e_2, e_4] = e_2, \  [e_3, e_4] = 0.
$$

   Отсюда легко получить представление поля
$$
   e_4 = (z_1 + A_4, z_2 + B_4,  C_4)
$$
с некоторыми комплексными константами $ A_4, B_4, C_4 $.

   Аналогично, из соотношений
$$
  [e_1, e_5] = -e_2, \  [e_2, e_5] = e_1, \  [e_3, e_5] = 0, \
$$
получаем представление поля
$$
   e_5 = (z_2 + A_5, -z_1 + B_5, C_5).
$$

   Но тогда коммутатор
$$
    [e_4,e_5] =  (z_1 + A_4)(0,-1,0) +(z_2 + B_4)(1,0,0) - (z_2 + A_5)(1,0,0) - (-z_1 + B_5)(0,1,0)
$$
имеет нулевую третью компоненту в противоречии с коммутационным соотношением
$$
    [e_4,e_5] = e_3 = (0,0,1).
$$

   Для алгебры $ \mathfrak g_{5,38} $ третьи компоненты полей $ e_4, e_5 $ по аналогичным
соображениям также могут быть только константами, а значит их коммутатор также не может быть равным $ e_3 $.

   Два этих противоречия доказывают предложение 5.1.
\hfill{$\Box$}

   Далее для пары алгебр $ \mathfrak g_{5,38}$ и $ \mathfrak g_{5,39} $ мы обсудим {\bf второй случай}, в рамках которого первая тройка базисных полей приводится к виду
$$
\begin{array}{ll}
  \, & e_1 = \left(1, 0, 0\right), \\
  \, & e_2 = \left(\hat f_2, 0, \hat h_2\right), \\
  \, & e_3 = \left(0, 0, 1\right).
\end{array}
\eqno (5.3)
$$

{\bf Предложение 5.2.} {\it Голоморфные  реализации алгебр $ \mathfrak g_{5,38}$ и $ \mathfrak g_{5,39} $, для которых
выполняются условия (5.3), могут иметь только вырожденные по Леви орбиты.
}

 Доказательство.
  Предполагая существование голоморфных реализаций с условиями (5.3) у
каждой из двух обсуждаемых алгебр, обсудим информацию о двух оставшихся базисных полях.

  Так, для алгебры $ \mathfrak g_{5,38}$ из соотношений
$ [e_1,e_5] = 0, \ [e_3,e_5] = 0 $
получаем независимость поля $ e_5 $ от переменных $ z_1, w $. В стандартных для этой статьи обозначениях имеем
$
   e_5 = (\hat f_5, \hat g_5, \hat h_5).
$

  Повторяя рассуждения Леммы 5, получим (в случае существования невырожденных орбит у этой алгебры) тройку выпрямленных полей
$ e_1, e_3, e_5 $ (составляющих базис еще одного абелева идеала в алгебре  $ \mathfrak g_{5,38}$).
  Рассматривая коммутаторы поля $ e_4 $ с полями $ e_1, e_3, e_5 $, получаем
$$
   e_4 = (z_1 + A_4, B_4, -w + C_4).
$$

  В силу независимости полученных компонент поля $ e_4 $ от переменной $ z_2 $, оно коммутирует с полем $ e_5 = (0,1,0) $
в противоречии с коммутационным соотношением $ [e_4,e_5] = e_3  $.
Тем самым, для алгебры $ \mathfrak g_{5,38}$ утверждение предложения 5.2 доказано.

  В случае голоморфной реализации алгебры $ \mathfrak g_{5,39} $ условия (5.3) и коммутационные соотношения
$$
   [e_1,e_4] = e_1, \ [e_3,e_4] = 0
$$
позволяют представить поле $ e_4 $ в виде
$$
   e_4 = (z_1 + \hat f_4, \hat g_4, \hat h_4).
$$

   Для алгебры, допускающей невырожденные орбиты и имеющей базис с такими свойствами, функция $ \hat g_4(z_2) $ не может быть
 тождественно нулевой (в силу леммы 4 о четырех нулях). Тогда, считая $ \hat g_4(0) \ne 0 $, по Лемме 5 поле $ e_4 $
 можно привести к упрощенному виду
$$
   e_4 = (z_1, 1, 0),
$$
сохранив при этом оговоренный выше вид полей $ e_1, e_2, e_3 $.

   При этом коммутатор
$$
   [e_2,e_4] = \hat f_2 (1,0,0) -  (\hat f_2', 0, \hat h_2')
$$
равен $ e_2 $. Это означает, что
$$
   e_2 = (A_2, 0, C_2 e^{-z_2})
$$
 с некоторыми комплексными константами $ A_2, C_2 $.

    Рассмотрим теперь коммутаторы $ [e_1, e_5] = -e_2 $ и $ [e_3, e_5] = 0 $. Из этих двух равенств следует, что
$$
   e_5 = (-z_1 A_2 + A_5, \hat g_5, -z_1 C_2 e^{-z_2} + C_5)
$$
 со своими комплексными константами $ A_5, C_5 $.

   Но тогда коммутационное соотношение
$$
   [e_4,e_5] = z_1(-A_2, 0, -C_2 e^{-z_2}) + (0, \hat g_5', z_1 C_2 e^{-z_2}) - (-z_1 A_2 + A_5)(1,0,0) = e_3 = (0,0,1)
$$
не может выполняться в третьей компоненте, т.к.
$$
    -z_1 C_2 e^{-z_2}  + z_1 C_2 e^{-z_2} = 0 \ne 1.
$$

   Предложение 5.2 доказано полностью.
\hfill{$\Box$}

\

{\bf 5.2. Семейства алгебр $ \mathfrak g_{5,33}, \mathfrak g_{5,34}, \mathfrak g_{5,35}$: первый случай}

\

  Рассмотрим теперь первые три семейства алгебр из набора (5.1). Как обычно, начинаем с {\bf первого случая}, в
котором тройка базисных полей абелевой подалгебры
$ I = < e_1, e_2, e_3 > $ имеет выпрямленный вид (5.2).

{\bf Предложение 5.3.} {\it
Пусть $ M $ - голоморфно однородная невырожденная по Леви вещественная гиперповерхность в $ \Bbb C^3 $, на которой имеется алгебра
$ \mathfrak g(M) $
голоморфных векторных полей со структурой одной из алгебр
семейства
$$
  \mathfrak g_{5,33}, \mathfrak g_{5,34}, \mathfrak g_{5,35}.
\eqno (5.4)
$$
Если для базиса $ \mathfrak g(M) $
выполняются условия (5.2), то $ M $ голоморфно эквивалентна трубке над аффинно однородным основанием.
}

 Доказательство.
 Как неоднократно отмечалось в предыдущих разделах,
 5-мерные орбиты алгебр, удовлетворяющих условиям (5.2), являются трубчатыми
гиперповерхностями. Остается показать аффинную однородность оснований таких трубок для алгебр (5.4).

    Рассмотрение коммутаторов пары оставшихся базисных полей $ e_4, e_5 $ с базисными полями подалгебры $ I $ приводит
во всех трех семействах алгебр к аффинным векторным полям, описываемым очень простыми формулами:
$$
   \mathfrak g_{5,33} \ : \ e_4 = (z_1 + A_4,  B_4, \beta w + C_4), \quad
                            e_5 = (A_5, z_2 +  B_5, \gamma w + C_5);
\eqno (5.5)
$$

$$
   \mathfrak g_{5,34} \ : \ e_4 = (\alpha z_1 + A_4, z_2 + B_4, w + C_4), \quad
                            e_5 = (\alpha z_1 + A_5, w + B_5,  C_5);
\eqno (5.6)
$$

$$
   \mathfrak g_{5,35} \ : \ e_4 = (h z_1 + A_4, z_2 + B_4, w + C_4), \quad
                            e_5 = (\alpha z_1 + A_5, w + B_5,  -z_2 + C_5).
\eqno (5.7)
$$

   В этих формулах все коэффициенты при переменных являются вещественными. С учетом этого <<овеществим>> поля
(5.5) - (5.7) и переведем их из $ \Bbb C^3 $ в $ \Bbb R^3 $.
Полю
$$
   e_k = (\lambda_k (z_1,z_2,w) + A_k)\frac{\partial}{\partial z_1} + (\mu_k (z_1, z_2, w)  + B_k) \frac{\partial}{\partial z_2} +
          (\nu_k (z_1,z_2,w)  + C_k) \frac{\partial}{\partial w} \ (k = 4,5),
$$
где $ \lambda_k, \mu_k, \nu_k $ --- линейные формы в $ \Bbb C^3 $ с вещественными коэффициентами,
при таком овеществлении сопоставляется поле
$$
   E_k = (\lambda_k (y_1, y_2,v) + Im\,A_k)\frac{\partial}{\partial y_1} + (\mu_k (y_1,y_2, v) + Im\,B_k) \frac{\partial}{\partial y_2} +
          (\nu_k (y_1, y_2,v) + Im\,C_k) \frac{\partial}{\partial v}.
$$

   Легко видеть, что основание 5-мерной трубки $ M $ из $ \Bbb C^3 $, описываемое уравнением $ \Phi(y_1,y_2,v) = 0 $ и системой двух соотношений
$$
    Re\left(e_k(\Phi)_{|M}\right) \equiv 0, \ k = 4,5
$$
удовлетворяет при таком пересчете системе
$$
    E_k(\Phi)_{| \Phi = 0} \equiv 0, \ k = 4,5.
$$

  Тем самым, основания всех интересующих нас в этом предложении трубок являются орбитами двумерных коммутативных
($ [e_4,e_5] = [E_4,E_5] = 0 $) алгебр, т.е. аффинно-однородными поверхностями в $ \Bbb R^3 $.
  Предложение 5.3 доказано.
\hfill{$\Box$}

{\bf Замечание.} В качестве дополнения к доказательству предложения 5.3 уточним, что орбитами
двухпараметрических семейств алгебр $ \mathfrak{g}_{33} $
и $ \mathfrak{g}_{35} $ в этом случае являются, соответственно, поверхности
$$
   v = y_1^{\gamma} y_2^{\beta} \quad \mbox{и} \quad
   v = |y_1 + i y_2|^{\beta} e^{\alpha \arg(y_1 + i y_2)}.
\eqno (5.8)
$$

  Однопараметрическое семейство однородных поверхностей
$$
   v = y_2\, (\ln y_1 - \alpha \ln y_2)
\eqno (5.9)
$$
соответствует (при $\alpha \ne 0 $) алгебрам семейства $ \mathfrak{g}_{34} $.

\

{\bf 5.3. Семейства алгебр $ \mathfrak g_{5,33}, \mathfrak g_{5,34}, \mathfrak g_{5,35}$: второй случай}

\

   Для набора алгебр (5.4) остается рассмотреть {\bf второй случай}, в рамках которого выпрямляется не вся
тройка базисных векторных полей $ e_1, e_2, e_3 $, а только два поля из трех.
Отметим при этом, что для рассматриваемой тройки алгебр удобно считать выпрямленными поля $ e_2 $ и $ e_3 $. При выпрямлении
любой другой пары из этой тройки обсуждение полного набора коммутационных соотношений в алгебрах
$ \mathfrak g_{5,33}, \mathfrak g_{5,34}, \mathfrak g_{5,35} $
сталкивается с излишними техническими трудностями.

  Итак, рассмотрим набор (5.4), считая, что тройка базисных полей абелева идеала любой из этих алгебр упрощена до вида
$$
\begin{array}{ll}
  \, & e_1 = \left(\hat f_1, 0, \hat h_1\right), \\
  \, & e_2 = \left(1, 0, 0\right), \\
  \, & e_3 = \left(0, 0, 1\right).
\end{array}
\eqno (5.10)
$$
  Рассматривая в таком случае коммутаторы выпрямленных полей $ e_2, e_3 $ с полями $ e_4, e_5 $, получаем следующие предварительные
упрощения полей $ e_4, e_5 $:
$$
   g_{5,33} \  : \ e_4 = (\hat f_4, \hat g_4, \beta w + \hat h_4), \quad
                   e_5 = (z_1 + \hat f_5, \hat g_5,  \gamma w + \hat h_5),
$$
$$
   g_{5,34} \  : \ e_4 = (z_1 + \hat f_4, \hat g_4, w + \hat h_4), \quad
                   e_5 = (w + \hat f_5, \hat g_5, \hat h_5),
\eqno (5.11)
$$
$$
   g_{5,35} \  : \ e_4 = (z_1 + \hat f_4, \hat g_4, w + \hat h_4), \quad
                   e_5 = (w + \hat f_5, \hat g_5, -z_1 + \hat h_5).
$$

   Поле $ e_4 $ можно подвергнуть дополнительному упрощению, используя лемму 4 (<<лемму о четырех нулях>>). Т.к., в силу этой леммы, на невырожденной поверхности $ g_4(z_2) $ не может быть тождественно нулевой, можно превратить эту функцию в единицу, а затем <<линеаризовать>> все поле
$ e_ 4 $. В итоге вместо (5.11) получаем еще более простой вид этого поля
$$
     e_4 = (0, 1, \beta w)\quad \mbox{для семейства} \quad \mathfrak g_{5,33} \quad \mbox{или} \quad
     \ e_4 = (z_1, 1, w ) \quad  \mbox{для} \quad \mathfrak g_{5,34}, \mathfrak g_{5,35}.
$$

  Выпишем упрощенный в этом случае вид базисов всех трех семейств $ \mathfrak g_{5,33}, \mathfrak g_{5,34}, \mathfrak g_{5,35}$.
$$
\begin{array}{ll}
  \, & e_1 = \left(\hat f_1, 0, \hat h_1\right), \\
  \, & e_2 = \left(1, 0, 0\right), \\
  \, & e_3 = \left(0, 0, 1\right), \\
  \, & e_4 = (0, 1, \beta w), \\
  \, & e_5 = (z_1 + \hat f_5, \hat g_5,  \gamma w + \hat h_5),
\end{array}
\
\begin{array}{ll}
  \, & e_1 = \left(\hat f_1, 0, \hat h_1\right), \\
  \, & e_2 = \left(1, 0, 0\right), \\
  \, & e_3 = \left(0, 0, 1\right), \\
  \, & e_4 = (z_1, 1, w), \\
  \, & e_5 = (w + \hat f_5, \hat g_5, \hat h_5),
\end{array}
\
\begin{array}{ll}
  \, & e_1 = \left(\hat f_1, 0, \hat h_1\right), \\
  \, & e_2 = \left(1, 0, 0\right), \\
  \, & e_3 = \left(0, 0, 1\right), \\
  \, & e_4 = (z_1, 1,  w), \\
  \, & e_5 = (w + \hat f_5, \hat g_5,  -z_1 + \hat h_5),
\end{array}
\eqno (5.12)
$$

  Для получения итоговых выводов о голоморфных реализациях алгебр из этих семейств в случае (5.10) остается рассмотреть
по три коммутационных соотношения, связанных с $ [e_1,e_4], [e_1,e_5], [e_4,e_5] $.

   Сначала мы обсудим два семейства алгебр из набора (5.12).

{\bf Предложение 5.4.} {\it Голоморфные реализации в пространстве $ \Bbb C^3 $ алгебр $ \mathfrak g_{5,33}, \mathfrak g_{5,34} $,
удовлетворяющие условиям (5.10), могут иметь только вырожденные по Леви 5-мерные орбиты.
}

  Доказательство.
  Для алгебр семейства $ \mathfrak g_{5,33} $ из соотношения
$ [e_1,e_4] = e_1 $
имеем:
$$
    \hat h_1 (0,0,\beta) - (\hat f_1',0,\hat h_1') = (\hat f_1,0,\hat h_1).
$$

   Две нетривиальные компоненты этого векторного равенства приводят к формуле
$$
     e_1 =(A_1 e^{-z_2},0, C_1e^{(\beta-1) z_2})
\eqno (5.13)
$$
с некоторыми комплексными константами $ A_1, C_1 $
   Аналогично, расписывая покомпонентно соотношение
$$
 [e_4,e_5] = \left( (\hat f_5', \hat g_5', \hat h_5') + \beta w (0,0,\gamma)\right) - (\gamma w + \hat h_5)(0,0,\beta) = 0,
$$
получим формулу
$$
    e_5 = (z_1 + A_5, B_5, \gamma w + C_5 e^{\beta z_2} )
$$

  Наконец, вычислим для таких полей коммутатор
$$
   [e_4, e_5] = \left(A_1 e^{-z_2}(1,0,0) + C_1e^{(\beta-1) z_2}(0,0,\gamma)\right) -
                B_5\left(-A_1 e^{-z_2},0, (\beta-1) C_1e^{(\beta-1) z_2}\right) = 0.
$$

   Отсюда получаем
$$
A_1 e^{-z_2} + B_5 A_1 e^{-z_2} = 0,
\quad 0 = 0, \quad
C_1e^{(\beta-1) z_2}\gamma - B_5 (\beta-1) C_1e^{(\beta-1) z_2} = 0
$$
или, после освобождения от экспонент,
$$
    A_1(B_5 + 1) = 0, \quad C_1 (\gamma - B_5 (\beta-1)) = 0.
\eqno (5.14)
$$

  Первое из этих равенств приводит к вырождению всех орбит алгебры $ \mathfrak g_{5,33} $ с базисом, реализованным по
схеме (5.12). В самом деле, при $ A_1 = 0 $ оказываются линейно зависимыми над $ \Bbb C $ поля $ e_1, e_2 $.
 Если же $ B_5 = -1 $, то для четверки полей
$ e_1, e_2, e_3, e_4 + e_5 $ выполняется условие леммы о четырех нулях.
Поскольку других возможностей для выполнения первого равенства (5.12) нет, для алгебр $ \mathfrak g_{5,33} $ утверждение предложения 5.4 верно.

  Аналогично, в случае алгебр $ \mathfrak g_{5,34} $ рассмотрим три последних коммутационных соотношения.
$$
   [e_1,e_4] = \left( \hat f_1 (1,0,0) + \hat h_1(0,0,1) \right) - (\hat f_1',0,\hat h_1') = \alpha e_1.
$$

  Отсюда получаем покомпонентные равенства
$$
   \hat f_1' = (1- \alpha) \hat f_1, \quad 0 = 0, \quad \hat h_1' = (1- \alpha) \hat h_1,
$$
и как следствие,
$$
    e_1 = e^{(1- \alpha)z_2}(A_1, 0, C_1).
$$

  Из второго соотношения
$$
   [e_4, e_5] = \left( (\hat f_5', \hat g_5', \hat h_5') + w (1,0,0) \right) -
                \left( (w + \hat f_5)(1,0,0) + \hat h_5 (0,0,1) \right) = 0
$$
получаем формулу для
$$
    e_5 = (w + A_5 e^{z_2}, B_5, C_5e^{z_2})
$$
с комплексными константами $ A_5, B_5, C_5 $.

  Последнее коммутационное соотношение тогда имеет вид
$$
   [e_1,e_5] = C_1  e^{(1 - \alpha)z_2}(1,0,0) - B_5 (1- \alpha) e^{(1 - \alpha)z_2}(A_1, 0, C_1) = e_1.
$$

   Из двух компонент этого равенства получаем (после освобождения от экспоненциальных множителей):
$$
   C_1 - B_5 (1- \alpha) A_1 = A_1, \quad
   - B_5 (1- \alpha) C_1 = C_1.
$$

  Переписывая эти равенства в виде
$$
   C_1 + A_1 (1 + B_5 (1- \alpha)), \quad
   C_1 (1 + B_5 (1- \alpha)) = 0,
$$
приходим к выводу $ C_1 = 0 $. Это означает, как и в предыдущем обсуждении, линейную зависимость над $ \Bbb C $ полей $ e_1, e_2 $ и, следовательно, вырождение всех орбит рассматриваемых реализаций алгебр $ \mathfrak g_{5,34} $.
  Предложение 5.4 полностью доказано.
\hfill{$\Box$}

 Наиболее интересным с точки зрения наличия невырожденных орбит является третье семейство из (5.12).

{\bf Предложение 5.5.} {\it Базис любой голоморфной реализации алгебры из семейства $ \mathfrak g_{5,35} $,
удовлетворяющей условиям (5.10)
и имеющей невырожденные орбиты, можно привести голоморфной заменой координат к следующему виду
($ A_5, C_5 $ - произвольные комплексные константы, $ \beta \ne 1 $):
$$
   e_1 = A_1 e^{(1-\beta) z_2} \left(1, 0, i\varepsilon\right), \
  \,  e_2 = \left(1, 0, 0\right), \ \
  \,  e_3 = \left(0, 0, 1\right), \ \
  \,  e_4 = (z_1, 1,  w), \ \
\eqno (5.15)
$$$$
  \,  e_5 =\left(w + A_5 e^{z_2}, \frac{\alpha - i \varepsilon}{\beta - 1},  -z_1 + C_5e^{z_2}\right).
$$
}

{\bf Замечание.} При $ \beta = 1 $ реализация алгебр из этого семейства с выполнением
условий (5.10) невозможна.

  Доказательство.
  Как и в предыдущем предложении, рассмотрим для третьего типа базиса (5.12) коммутатор
$$
   [e_1,e_4] = \left(\hat f_1 (1,0,0) + \hat h_1(0,0,1)\right) - (\hat f_1',0, \hat h_1') = \beta e_1.
$$

   Две содержательных скалярных компоненты этого векторного равенства имеют вид системы двух ОДУ
$$
    \hat f_1'  = (1-\beta) \hat f_1, \quad   \hat h_1'  = (1-\beta) \hat h_1,
$$
т. что поле $ e_1 $ можно записать в виде
$$
    e_1 = e^{(1-\beta)z_2} (A_1, 0, C_1)
\eqno (5.16)
$$
с некоторыми комплексными константами $ A_1, C_1 $. Следующее соотношение $ [e_4, e_5] = 0 $ имеет в развернутой форме вид
$$
   \left( z_1(0,0,-1) +  (\hat f_5',\hat g_5', \hat h_5') + w (1,0,0) \right) -
   \left( (w + \hat f_5 )(1,0,0) + (-z_1 + \hat h_5)(0,0,1) \right) = 0.
$$

  В трех  компонентах этого равенства имеем
$$
   \hat f_5' + w - (w + \hat f_5)  = 0, \quad  \hat g_5' = 0, \quad - z_1 + \hat h_5' - (-z_1 + \hat h_5) = 0.
$$

  Из этой системы ОДУ получаем формулу
$$
    e_5 = (w+ A_5 e^{z_2}, B_5, -z_1 + C_5e^{z_2} )
$$
с некоторыми комплексными константами $ A_5, B_5, C_5 $. Тогда завершающее коммутационное соотношение $ [e_1, e_5] = 0 $ с полученными
выражениями для полей $ e_1, e_5 $ примет вид
$$
   e^{(1-\beta)z_2} ( A_1(0,0,-1) + C_1 (1,0,0) ) - B_5 (1-\beta)e^{(1-\beta)z_2}(A_1,0,C_1) = 0.
\eqno (5.17)
$$

  В виде, свободном от экспонент, получаем здесь два соотношения на коэффициенты, входящие в (5.17):
$$
    C_1 = (B_5(1-\beta)+\alpha) A_1, \qquad A_1 + (B_5(1-\beta)+\alpha) C_1 = 0.
$$

  Подставляя первое из них во второе, приходим к равенству
$$
    A_1 \left( 1 + (B_5(1-\beta)+\alpha)^2\right) = 0.
$$

  Заметим, что равенство $ A_1 = 0 $ здесь невозможно, т.к. при таком условии базисное поле $ e_1 $ обсуждаемой алгебры оказывается нулевым.
Следовательно,
$$
    (B_5(1-\beta)+\alpha)^2 = -1 \quad \mbox{или} \quad B_5(1-\beta)+\alpha =  i\varepsilon, \ \varepsilon = \pm 1.
$$

   Ясно, что при $ \beta = 1 $ эти равенства выполняться не могут, и это уточнение доказывает справедливость
замечания к формулировке доказываемого предложения.
   Итоговый же позитивный вывод в этом случае состоит в том, что алгебра полей с базисом (5.12) возможна тогда и только
тогда, когда базисные поля имеют уточненный вид
$$
\begin{array}{ll}
  \, & e_1 = \left(A_1 e^{\lambda z_2}, 0, i\varepsilon A_1 e^{\lambda z_2}\right), \\
  \, & e_2 = \left(1, 0, 0\right), \\
  \, & e_3 = \left(0, 0, 1\right), \\
  \, & e_4 = (z_1, 1, w), \\
  \, & e_5 = (w + A_5 e^{z_2}, B_5,  -z_1 + C_5e^{z_2}),
\end{array}
\eqno (5.18)
$$
а параметры и коэффициенты, участвующие в этих формулах, удовлетворяют условиям
$$
   A_1 \ne 0, \quad \beta \ne 1, \quad \lambda = 1- \beta, \quad B_5 =  \frac{i\varepsilon - \alpha}{\lambda}.
\eqno (5.19)
$$

  Для завершения доказательства предложения 5.5 остается заметить, что ненулевая константа $ A_1 $ может быть записана в виде
$
   A_1 = \exp(\ln A_1),
$
т. что после замены $ z_2^* = z_2 + (1/\lambda) \ln A_1 $
поле $ e_1 $ примет вид
$$
    e_1 = \left( e^{\lambda z_2}, 0, i\varepsilon  e^{\lambda z_2}\right),
$$
заявленный в предложении 5.5.
\hfill{$\Box$}

{\bf Замечание 1.} Строго говоря, вид базиса, полученный в этом предложении, не гарантирует невырожденности орбит соответствующих алгебр. Такой вид является {\it необходимым}, но не достаточным условием невырожденности.

{\bf Замечание 2.} По аналогии с удалением из формул (5.18) константы $ A_1 $ можно упростить и остающийся в этих формулах
набор  $ (A_5, C_5) $. Это будет сделано (по-разному в нескольких случаях) в процессе интегрирования
полученных алгебр, отложенного до следующего раздела.

\

  В завершение же этого раздела мы обсудим две <<исключительные>>
  алгебры $ \mathfrak g_{5,36} $ и $ \mathfrak g_{5,37} $ из рассматриваемого блока VI. В частности,
будут построены упрощенные базисы возможных реализаций в $ \Bbb C^3 $ алгебры $ \mathfrak g_{5,37} $. Трудоемкое интегрирование
получаемых при этом алгебр (а также семейства алгебр $ \mathfrak g_{5,35} $) будет осуществлено в следующем разделе статьи.

Отметим, что алгебры $ \mathfrak g_{5,36} $ и $ \mathfrak g_{5,37} $, как и пара $ \mathfrak g_{5,25}, \mathfrak g_{5,26} $ рассмотрены в работе [33]
с точки зрения наличия у них <<просто однороднодных>> строго псевдо-выпуклых орбит. Обсуждения настоящей работы включают помимо СПВ-случая и поверхности с индефинитной формой Леви.

\

{\bf 5.4. Орбиты алгебры $ \mathfrak g_{5,36} $}

\

{\bf Предложение 5.6.} {\it Пятимерными орбитами в $ \Bbb C^3 $ голоморфных реализаций алгебр $ \mathfrak g_{5,36} $ являются
либо вырожденные по Леви гиперповерхности, либо голоморфные образы индефинитной квадрики $ v = |z_1|^2 - |z_2|^2 $.
}

   Доказательство.
   Предположим, как обычно, что $ M \subset \Bbb C^3 $ - однородная невырожденная по Леви вещественная гиперповерхность, на
которой имеется 5-мерная алгебра Ли $ \mathfrak g(M) $голоморфных векторных полей со структурой $ \mathfrak g_{5,36} $.

   Два коммутирующих базисных поля $ e_1 $ и $ e_2 $ этой алгебры выпрямим до состояния
$$
   e_1 = (0,0,1), \ e_2 = (1,0,0).
$$

    Тогда коммутационные соотношения
$$
   [e_1,e_3] =0, \ [e_2, e_3] = e_1,
$$
совпадающие в обеих алгебрах, приводят к следующему виду поля $ e_3 $:
$$
    e_3 = (\hat f_3, \hat g_3, z_1 + \hat h_3).
\eqno  (5.20)
$$

 Аналогичным образом из рассмотрения коммутаторов первой пары полей $ e_1, e_2 $ с полем $ e_4 $, а затем --- той же
 пары с полем $ e_5 $, приходим к формулам
$$
    e_4 = (z_1 + \hat f_4, \hat g_4, w + \hat h_4).
\eqno  (5.21)
$$$$
    e_5 = (-z_1 + \hat f_5, \hat g_5, \hat h_5).
$$

  Далее рассмотрим два случая в зависимости от свойств функции $ \hat g_3 $.

\

  {\bf В первом случае}, при $ g_3(0) \ne 0 $, стандартными голоморфными преобразованиями, описанными в начальных разделах этой работы, приводим поле $ e_3 $ к виду
$$
   e_3 = (0,1,z_1).
$$

  Далее из тройки коммутационных соотношений
$$
   [e_1,e_4] = e_1,\  [e_2,e_4] = e_2,\ [e_3,e_4] = 0,
$$
связанных с полем $ e_4 $, получаем (уточненную по сравнению с (5.21)) формулу
$$
   e_4 = (z_1 + A_4, B_4, w + A_4 z_2 + C_4)
\eqno (5.22)
$$
с некоторыми комплексными константами $ A_4, B_4, C_4 $.

   А из аналогичной тройки
$$
   [e_1,e_5] = 0,\  [e_2,e_5] =- e_2,\ [e_3,e_5] = e_3,
$$
связанной с $ e_5 $, --- аналогичную формулу
$$
   e_5 = (- z_1 + A_5, z_2 + B_5,  A_5 z_2 + C_5).
\eqno (5.23)
$$
с некоторыми комплексными константами $ A_5, B_5, C_5 $.

  Тогда с учетом формул (5.22) и (5.23) последнее коммутационное соотношение $ [e_4,e_5] = 0 $ примет в координатной форме вид
следующей системы трех уравнений:
$$
  -(z_1+ A_4) - (- z_1 + A_5) =0, \quad B_4 = 0, \quad   B_4 A_5 - A_4 (z_2 + B_5) - ( A_5 z_2 + C_5) = 0.
$$

  Из этих уравнений следует, что $ A_4 + A_5 =0, \ B_4 = 0, \ C_5 = - A_4 B_5 $, т. что пара полей $ e_4, e_5 $ примет вид
$$
    e_4 = (z_1 + A_4, 0,  w + A_4 z_2 + C_4), \quad
       e_5 = (- (z_1 + A_4), z_2 + B_5,  -A_4 (z_2 + B_5)).
$$

   Сдвиг переменных
$$
   z_1 + A_4 \rightarrow z_1, \ z_2 + B_5 \rightarrow z_2, \ w + C_4 - A_4 B_5)
$$
позволяет записать набор базисных полей изучаемой алгебры в следующей (аффинной) форме
$$
\begin{array}{ll}
  \, & e_1 = \left(0, 0, 1\right), \\
  \, & e_2 = \left(1, 0, 0\right), \\
  \, & e_3 = \left(0, 1, z_1 - A_4\right)\\
  \, & e_4 = \left(z_1, 0, w + A_4 z_2\right)\\
  \, & e_5 = \left(-z_1, z_2, -A_4 z_2\right).
\end{array}
\eqno (5.24)
$$

   Интегрирование системы уравнений, отвечающих такой пятерке полей, приводит к уравнениям орбит
$$
   v = (y_1 - b) x_2 - a y_2 + D y_1 y_2,
\eqno (5.25)
$$
где $ \ a = Re\,A_4, \ b = Im\, A_4, \ D $ - {произвольная константа}.

  Несмотря на наличие трех вещественных параметров в формуле (5.25), легко понять, что аффинным преобразованием
с подходящим набором коэффициентов каждое такое уравнение сводится к
$
     v = y_1 y_2.
$
  В свою очередь, трубчатая поверхность с таким уравнением голоморфно эквивалентна стандартной индефинитной
  квадрике $ v = |z_1|^2 - |z_2|^2 $.

\

  Во {\bf втором случае}, при $ \hat g_3(z_2) \equiv 0 $, воспользуемся леммой 4 о четырех нулях. Интересуясь
только невырожденными орбитами, можно в силу этой леммы считать, что $ \hat g_5 (z_2) $ и $ \hat g_4 (z_2) $ отличны от
тождественно нулевых функций (на любой такой поверхности).

Сместимся в близкую к началу координат точку $ p $
обсуждаемой однородной поверхности $ M $, в которой выполняется неравенство $ \hat g_5(z_2)_{|p} \ne 0 $.
Стандартными голоморфными преобразованиями поле $ e_5 $ из (5.21) можно упростить до состояния
$$
   e_5 = (-z_1, 1, 0)
\eqno (5.26)
$$
с сохранением вида (5.21) поля $ e_4 $ и выпрямленных полей $ e_1 $ и $ e_2 $.

С учетом
 проведенных упрощений базисных полей и, в частности, формулы (5.26),
рассмотрим теперь коммутационные соотношения $ [e_3,e_5] = e_3 $ и $ [e_4,e_5] = 0 $. Из них получаются очередные
уточнения вида базисных полей
$$
    e_3 = (A_3 e^{-2 z_2}, 0, C_3 e^{-z_2}), \quad e_4 = (A_4 e^{-z_2}, B_4, w + C_4)
\eqno (5.27)
$$
c некоторыми комплексными константами $ A_k, B_k, C_k $.

   Наконец, для существования голоморфной реализации алгебры $ g_{5,36} $ (с невырожденными орбитами) в этом случае
   необходимо выполнение последнего коммутационного соотношения
$  [e_3,e_4] = 0 $. Вычисляя коммутатор для полей (5.27), имеем
$$
  [e_3,e_4] = C_3 e^{-z_2}(0,0,1) - B_4 \left(-2 A_3 e^{-2 z_2}, 0, -C_3 e^{-z_2}\right) = ((2 B_4 A_3 e^{-2 z_2}, 0, C_3(1+B_4) e^{-z_2}).
$$

   Необходимыми условиями равенства нулю такого поля являются два соотношения на коэффициенты:
$$
       B_4 A_3 = 0, \quad C_3(1+B_4) = 0.
$$

   Множество формальных решений этой системы распадается на 3 случая:
$$
   1) \ A_3 = 0, \ C_3 = 0,  \quad 2) \ B_4 = 0, \ C_3 = 0, \quad 3) \ B_4 = -1, \  A_3 = 0.
$$

   При этом во всех трех случаях у базисного векторного поля $ e_3 $ из (5.27) оказывается не более
одной ненулевой компоненты. В итоге в первом случае базисное (!) векторное поле $ e_3 $ оказывается тождественно нулевым, что невозможно. А в двух других  случаях оказывается линейно зависимой над $ \Bbb C $ либо пара полей $ e_1, e_3 $ (третий случай), либо
$ e_2, e_3 $ (второй случай). Но такая зависимость
означает вырождение по Леви орбиты соответствующей алгебры.

  Тем самым, допущение о существовании невырожденных орбит у голоморфной реализации алгебры $ g_{5,36} $ по схеме второго случая
противоречиво, а предложение 5.6 полностью доказано.
\hfill{$\Box$}

\

{\bf 5.5. Упрощение базисов голоморфных реализаций алгебры $ \mathfrak g_{5,37} $}

\

{\bf Предложение 5.7.} {\it Базис любой голоморфной реализации в $ \Bbb C^3 $ алгебры $ \mathfrak g_{5,37} $,
имеющей невырожденные орбиты, можно привести голоморфной заменой координат к одному из двух видов:
$$
\begin{array}{ccccccc}
  \,  e_1 : (0, & 0, & \ 1 \ \  \ ), \\
  \,  e_2 : (1, & 0, & \ 0 \ \ \ ), \\
  \,  e_3 : (0, & 1, & z_1-B \ \ ),\\
  \,  e_4 : (z_1, & z_2, & 2 w + B z_2 \ ),\\
  \,  e_5 : (z_2, & -z_1, & (z_2^2 - z_1^2)/2 + B z_1 )
\end{array}
\ \mbox{ или } \
 \begin{array}{ccccccc}
  e_1 :   ( \ \ \ 0\ , & 0\ , &\ 1 \ \  ) \\
  e_2 : (\ \ 1 \ , &  0 \ , & 0 \ ) \\
  e_3: ( \ A \ , & 0 \ , & z_1 ) \\
  e_4 : ( z_1  , &  z_2 \ , & 2 w ) \\
  e_5: ( -A z_1  \ , &  B z_2 \ , &  - z_1^2/2 ),
\end{array}
\eqno (5.28)
$$
где $ A = \pm i, B $ - некоторая комплексная константа.
}

  Доказательство.
   Рассмотрим невырожденную однородную гиперповерхность $ M $, являющуюся орбитой голоморфной
реализации алгебры $ \mathfrak g_{5,37} $. Выпрямим два базисных поля
$$
   e_1 = (0,0,1), \ e_2 = (1,0,0)
$$
этой алгебры и учтем два коммутационных соотношения
$
   [e_1,e_3] =0, \ [e_2, e_3] = e_1,
$
совпадающих в алгебрах $ \mathfrak g_{5,36} $ и $ \mathfrak g_{5,37} $.
  Тогда получаем поле
$$
    e_3 = (\hat f_3, \hat g_3, z_1 + \hat h_3)
$$
в той же форме (5.20) и традиционно рассматриваем два случая.

  {\bf В первом случае}, при $ g_3(0) \ne 0 $, снова получаем упрощение поля
$$
   e_3 = (0,1,z_1).
$$

  Рассматривая теперь коммутаторы тройки базисных полей $ e_1, e_2, e_3 $ сначала с полем $ e_4 $, а затем с $ e_5 $,
несложно получить упрощенные представления последней пары полей:
$$
   e_4 = (z_1 + A_4, z_2 + B_4, 2 w + A_4 z_2 + C_4),
\eqno (5.29)
$$
$$
   e_5 = (z_2 + A_5, -z_1 + B_5, \frac12(z_2 + A_5)^2 -\frac12 z_1^2 + C_5)
\eqno (5.30)
$$
c произвольными комплексными константами
$ A_k, B_k, C_k $.

   Следующий этап упрощения базиса алгебры $ \mathfrak g_{5,37} $ в первом случае связан с соотношением
$ [e_4, e_5] = 0 $ для полей (5.29) и (5.30). Покомпонентное рассмотрение этого соотношения дает ряд ограничений на комплексные константы:
$$
   A_5 = B_4, \ B_5 = - A_4, \ C_5 = \frac12 A_4^2.
$$

  В итоге, после сдвигов переменных
$$
    z_1 + A_4 \rightarrow z_1, \quad z_2 + B_4 \rightarrow z_2, \quad w + \frac12(C_4 - A_4 B_4) \rightarrow w
$$
и переобозначения константы $ A_4 = B \in \Bbb C  $ набор базисных полей алгебры $ \mathfrak g_{5,37} $
в этом случае принимает первый из двух видов (5.28) доказываемого предложения.

{\bf Замечание.} После квадратичной замены переменой $ w^* = w + (1/2)z_1 z_2 $ все поля из такого базиса
становятся аффинными:
$$
\begin{array}{ll}
  \, & e_1 = \left(0, 0, 1\right), \\
  \, & e_2 = \left(1, 0, z_2/2\right), \\
  \, & e_3 = \left(0, 1, 3 z_1/2 -A \right),\\
  \, & e_4 = \left(z_1, z_2, 2 w + A z_2 \right),\\
  \, & e_5 = \left(z_2, -z_1,  A z_1 \right).
\end{array}
\eqno (5.31)
$$

   Однако при этом в полученном базисе вместо двух выпрямленных полей остается лишь одно.

\

    Во {\bf втором случае}, при $ \hat g_3(z_2) \equiv 0 $ (и наличии невырожденной орбиты $ M $ у обсуждаемой алгебры),
можно считать по лемме 4 (о четырех нулях), что $ \hat g_4(z_2) \ne 0 $.
Тогда описанная выше техника позволяет привести
базис алгебры $ \mathfrak g_{5,37} $ в некоторых голоморфных координатах к виду
$$
 \begin{array}{ccccccc}
  e_1 :   ( \ \ \ 0\ , & 0\ , &\ 1 \ \  ) \\
  e_2 : (\ \ 1 \ , &  0 \ , & 0 \ ) \\
  e_3: ( \hat f_3(z_2) \ , & 0 \ , & z_1 + \hat h_3(z_2)\ ) \\
  e_4 : (\ \  z_1 \ , &  1 \ , & 2 w \ \ ) \\
  e_5: ( -z_1 \hat f_3(z_2) + \hat f_5(z_2), & \hat g_5 (z_2), & -(1/2) z_1^2 - z_1 \hat h_3(z_2) + \hat h_5(z_2)\ ).
\end{array}
\eqno (5.32)
$$

  Рассмотрим далее три коммутационных соотношения
$
  [e_3, e_4] =e_3
$,
$
  [e_3, e_5] = e_2
$, $
  [e_4, e_5] = 0.
$

   Первое из них примет вид
$$
  \hat f_3 (1,0,0) + (z_1 + \hat h_3) (0,0,2) - z_1 (0,0,1) - (\hat f_3',0,\hat h_3') = (\hat f_3(z_2),  0, z_1 + \hat h_3(z_2)).
$$

  Две содержательные компоненты этого равенства упрощаются до
$
   \hat f_3' = 0, \ \hat h_3' - \hat h_3 = 0.
$
  Отсюда
$$
   e_3 = (A_3, 0, z_1 + C_3 e^{z_2}), \quad
\eqno (5.33)
$$
$$
   e_5 = \left(-A_3 z_1 + \hat f_5(z_2), \hat g_5 (z_2), -(1/2) z_1^2 - C_3 z_1 e^{z_2} + \hat h_5(z_2) \right)
$$
с произвольными комплексными константами $ A_3, C_3 $.

   Рассматривая теперь с учетом формул (5.33) равенство
$
  [e_3, e_5] = e_2
$,
получим
$$
  A_3 (-A_3, 0, -z_1 -C_3e^{z_2}) - (-A_3 z_1 +\hat f_5) (0,0,1) - \hat g_5 (0,0, C_3e^{z_2}) = (1,0,0).
$$

  В этом равенстве также имеется лишь две содержательных компоненты
$$
   -A_3^2 = 1, \quad -A_3 C_3 e^{z_2} - \hat f_5 - C_3 \hat g_5 e^{z_2} = 0,
$$
следствиями которых являются формулы
$$
   A_3=\pm i, \quad \hat f_5 = -  C_3 e^{z_2} (A_3 + \hat g_5).
\eqno (5.34)
$$

  Наконец, третье равенство
$
  [e_4, e_5] = 0
$
является в развернутой форме самым объемным:
$$
   z_1 (-A_3, 0, -z_1 - C_3 e^{z_2})+
   \left(-C_3 e^{z_2}(A_3+ \hat g_5 + \hat g_5'), \hat g_5', - C_3 z_1 e^{z_2} + \hat h_5' \right)+
$$
$$
  + \left(A_3 z_1 + C_3 e^{z_2} (A_3 + \hat g_5)\right) (1,0,0)
   -\left( (-1/2) z_1^2 - C_3 z_1 e^{z_2} + \hat h_5\right)(0,0,2) = 0.
$$

  После упрощения трех его отдельных компонент получим систему уравнений
$$
   C_3 e^{z_2} (A_3 + \hat g_5 + \hat g_5') + C_3 e^{z_2} (A_3 + \hat g_5) = 0,
\quad
   \hat g_5' = 0,
\quad
   \hat h_5' - 2 \hat h_5  = 0.
$$

  Очевидными следствиями этих уравнений являются соотношения
$$
   C_3 (A_3 + \hat g_5) = 0,
\quad
   \hat g_5 = B_5,
\quad
   \hat h_5  = C_5 e^{2z_2}
\eqno (5.35)
$$
с произвольными комплексными константами $ B_5 $ и $ C_5 $.

   С учетом формул (5.33) - (5.35) мы можем записать базис (5.32) в виде ($ A_3 = \pm i $)
$$
 \begin{array}{ccccccc}
  e_1 :   ( \ \ \ 0\ , & 0\ , &\ 1 \ \  ) \\
  e_2 : (\ \ 1 \ , &  0 \ , & 0 \ ) \\
  e_3: ( A_3 \ , & 0 \ , & z_1 + C_3 e^{z_2}\ ) \\
  e_4 : (\ \  z_1 \ , &  1 \ , & 2 w \ \ ) \\
  e_5: ( - A_3 z_1, & B_5, & -(1/2) z_1^2 - z_1 C_3 e^{z_2} + C_5 e^{2z_2}\ ).
\end{array}
\eqno (5.36)
$$

   Сделаем еще несколько несложных голоморфных замен координат для дальнейшего упрощения  формул (5.36).

   Во-первых, обозначим $ z_2^* = e^{z_2} $. Тогда дифференцирование $ {\partial}/{\partial z_2} $ заменится на
$ z_2^*\,{\partial}/{\partial z_2^*} $, а весь базис (5.36) запишется (после отбрасывания звездочек) в виде
$$
 \begin{array}{ccccccc}
  e_1 :   ( \ \ \ 0\ , & 0\ , &\ 1 \ \  ) \\
  e_2 : (\ \ 1 \ , &  0 \ , & 0 \ ) \\
  e_3: ( A_3 \ , & 0 \ , & z_1 + C_3 {z_2}\ ) \\
  e_4 : (\ \  z_1 \ , &  z_2 \ , & 2 w \ \ ) \\
  e_5: ( - A_3 z_1, & B_5 z_2, & -(1/2) z_1^2 - C_3 z_1{z_2} + C_5 z_2^2 \ ).
\end{array}
\eqno (5.37)
$$

     Полученное выше первое из условий (5.35) разделяет обсуждение на два подслучая:

     1) $ C_3 = 0 $,

     2) $ C_3 \ne 0, \   B_5 = - A_3 $.

  В первом подслучае, при $ C_3 = 0 $, весь базис обсуждаемой алгебры примет вид ($ A = A_3 = \pm i, B= B_5, C= C_5 $)
$$
 \begin{array}{ccccccc}
  e_1 :   ( \ \ \ 0\ , & 0\ , &\ 1 \ \  ) \\
  e_2 : (\ \ 1 \ , &  0 \ , & 0 \ ) \\
  e_3: ( \ A \ , & 0 \ , & z_1 ) \\
  e_4 : ( z_1  , &  z_2 \ , & 2 w ) \\
  e_5: ( -A z_1  \ , &  B z_2 \ , &  - z_1^2/2 + C z_2^2 ),
\end{array}
\eqno (5.38)
$$
совпадающий со вторым видом (5.28).

  {\bf Замечание.} При любом вещественном значении коэффициента $ B $ линейная комбинация полей $ -B e_4 + e_5 $ имеет
тождественно нулевую вторую компоненту. По лемме о четырех нулях такая алгебра Ли может иметь только вырожденные орбиты. Поэтому,
выделяя далее вещественную и мнимую части коэффициента $ B = b_1 + i b_2 $, можно всегда считать (для обсуждаемых
невырожденных орбит), что $ b_2 \ne 0 $.

  Остается еще второй подслучай, связанный с (5.35). Здесь, при $ C_3 \ne 0 $, мы воспользуемся еще одним голоморфным преобразованием координат.
   Линейная замена $ z_1^* = z_1 + C_3 z_2, \  z_2^* = z_2 $ и связанное с ней преобразование
$$
    \frac{\partial}{\partial z_2} = C_3 \frac{\partial}{\partial z_1^*} +  \frac{\partial}{\partial z_2^*}
$$
сохраняют вид полей $ e_1, e_2, e_4 $ и изменяют поля $ e_3, e_5 $. Так, после этой замены получим
$$
   e_3 = (A_1,0, z_1), \quad
   e_5 = \left( - A_3 z_1 + C_3(A_3 + B_5) z_2, B_5 z_2, -\frac12 z_1^2 +( \frac12 C_3^2 + C_5) z_2^2 \right).
$$

  В силу равенства $ B_5 = - A_3 $, выполняющегося в этом подслучае, с учетом обозначения
$ C =  \frac12 C_3^2 + C_5 $
весь базис обсуждаемой алгебры принимает вид
$$
 \begin{array}{ccccccc}
  e_1 :   ( \ \ \ 0\ , & 0\ , &\ 1 \ \  ) \\
  e_2 : (\ \ 1 \ , &  0 \ , & 0 \ ) \\
  e_3: ( A_3 \ , & 0 \ , & z_1 ) \\
  e_4 : (\ \  z_1 \ , &  z_2 \ , & 2 w ) \\
  e_5: ( - A_3 z_1, & - A_3 z_2, & -\frac12 z_1^2 +C z_2^2  ).
\end{array}
\eqno (5.39)
$$

   Уточним, что полученный вид является частным случаем
базиса (5.38) с дополнительным ограничением $ B = - A = -A_3 $.
  Наконец, заметим, что голоморфная замена
$$
   w^* = w - \frac{C}{2B}z_2^2
$$
сохраняет неизменными поля $ e_1, e_2, e_3, e_4 $ базиса (5.38) и
удаляет слагаемое $ C z_2^2 $ из записи поля $ e_5 $. Тем самым, базис (5.38) приводится ко второму случаю (5.28), и
предложение 5.7 доказано.
\hfill{$\Box$}

\

{\bf 6. Интегрирование <<трудных>> алгебр и итоговые выводы
}

\

   В этом разделе будут проинтегрированы два типа алгебр Ли $ \mathfrak g_{5,37} $ и $ \mathfrak g_{5,35} $, упрощенные базисы
которых в <<трудных>> случаях (5.28) и (5.15) были получены в предыдущем разделе. Орбиты семейства $ \mathfrak g_{5,35} $
из такого случая оставались некоторое время
гипотетическими претендентами на статус новых однородных поверхностей. Проверка этой гипотезы и подтверждение окончательных выводов об однородности составляют идейное содержание раздела.

\

{\bf 6.1. Сферичность Леви-невырожденных орбит алгебры $ \mathfrak g_{5,37} $}

\

{\bf Предложение 6.1.} {\it Все невырожденные по Леви пятимерные орбиты голоморфных реализаций алгебры $ \mathfrak g_{5,37} $ являются
сферическими поверхностями.
}

Доказательство.
   Рассмотрим невырожденную интегральную поверхность $ M $ алгебры Ли с фиксированным базисом любого из двух типов (5.28).
Все такие алгебры содержат касательное к поверхности $ M $ поле
$ e_1 = {\partial }/{\partial w} $.
Поэтому
для определяющей функции
$
   \Phi(z_1,z_2, w)
$
поверхности
$
  M = \{ \Phi(z_1,z_2, w) = 0 \}
$
выполняется равенство
$$
     \frac{\partial \Phi}{\partial u} = 2\, Re \left(\frac{\partial \Phi}{\partial w}\right) \equiv 0.
$$

   Так как $ M $ невырожденная, то в некоторой точке этой поверхности (которую мы теперь будем считать новым началом координат)
выполняется неравенство
$$
     \frac{\partial \Phi}{\partial v} = - 2\, Im\left(\frac{\partial \Phi}{\partial w}\right) \ne 0.
$$
  Тем самым, функцию
$
   \Phi(z_1,z_2, w)
$
можно считать имеющей (вблизи начала координат) вид
$$
   \Phi = - v + F(z, \bar z,u).
\eqno (6.1)
$$

    Более того, наличие в алгебре двух выпрямленных полей $ e_1 = {\partial }/{\partial w} $ и $ e_2 = {\partial }/{\partial z_1} $
означает, что в действительности функция $ F $ из (6.1) может зависеть только от трех вещественных переменных
$
   y_1 = Im\, z_1, x_2 = Re\, z_2, y_2 = Im\, z_2.
$

   С учетом этого выпишем условия касания обсуждаемой поверхности $ M $ тремя оставшимися базисными полями рассматриваемой алгебры Ли. Сделаем это отдельно для каждого из двух типов базисов (5.28).

\

{\bf Случай алгебр с базисами (5.28) первого типа}

   Выделим вещественные и мнимые части параметра $ B $ из формул для таких базисов, обозначая через $ a, b $, соответственно, $ Re\, B $ и $ Im\, B $.
   Тогда интересующие нас условия примут вид системы трех уравнений
$$
   \frac{\partial F}{\partial x_2} - (y_1 - b) = 0,
$$
$$
   y_1 \frac{\partial F}{\partial y_1} +
   x_2 \frac{\partial F}{\partial x_2} +
   y_2 \frac{\partial F}{\partial y_2}
-
   (2 F + a y_2 + b x_2) = 0,
\eqno (6.2)
$$
$$
   y_2 \frac{\partial F}{\partial y_1} -
   x_1 \frac{\partial F}{\partial x_2} -
   y_1 \frac{\partial F}{\partial y_2} -
    ((x_2 y_2 - x_1 y_1)+ (b x_1 + a y_1))
= 0.
$$

  Заметим, что в последнем уравнении имеется три слагаемых, содержащих множитель $ x_1 $, а именно
$$
   - x_1 \frac{\partial F}{\partial x_2} + x_1 y_1 - bx_1 = - x_1 \left(\frac{\partial F}{\partial x_2} - (y_1 - b)\right).
$$

    В силу первого уравнения системы (6.2) вся эта сумма равна нулю, т. что третье уравнение в (6.2) упрощается до
$$
   y_2 \frac{\partial F}{\partial y_1} -
   y_1 \frac{\partial F}{\partial y_2} -
    (x_2 y_2 + a y_1)
= 0.
\eqno (6.3)
$$

    Решение первого уравнения из тройки (6.2) имеет вид
$$
   F(y_1,x_2,y_2) =  x_2 (y_1 - b) + G(y_1, y_2),
$$
где $ G $ - произвольная аналитическая функция двух переменных.

  Подставляя частные производные этой функции в оставшиеся два уравнения системы (5.41), получим
$$
   y_1 \frac{\partial G}{\partial y_1} +
   y_2 \frac{\partial G}{\partial y_2}
=
   (2 G + a y_2 ) = 0,
\eqno (6.4)
$$
$$
   y_2 \frac{\partial G}{\partial y_1} -
   y_1 \frac{\partial G}{\partial y_2}
= a y_1.
\eqno (6.5)
$$

  Общее решение
$$
    G = - a y_2 + \varphi(y_1^2 + y_2^2)
$$
уравнения (6.4), получаемое из интегралов соответствующей характеристической системы ОДУ, подставим в уравнение (6.5). Так получается завершающее ОДУ
$$
    t \varphi'(t) = \varphi(t)
$$
относительно неизвестной функции $ \varphi(t) $.  Возвращаясь от его решения
$
   \varphi(t) =  N t, \ ( N $ --- произвольное вещественное)
к функции $ F $, получаем уравнение поверхности $ M $ в виде
$$
   v = x_2(y_1 -b) - ay_2  + N (y_1^2 + y_2^2).
\eqno (6.6)
$$

      Ясно, что любая невырожденная поверхность из этого семейства голоморфно эквивалентна одной из
стандартных квадрик $ v = |z_1|^2 \pm |z_2|^2 $. Тем самым, предложение 6.1 в этом случае верно.

{\bf Случай алгебр с базисами (5.28) второго типа}

  Базисы в этом случае содержат два комплексных параметра $ A, B $. Выделим у них вещественные и мнимые части, полагая
$$
  A = i a_2  \ (a_2 = \pm 1), \ B= b_1 + i b_2.
$$

   Тогда система трех уравнений на определяющую функцию $ \Phi = - v + F(y_1, x_2, y_2) $ имеет вид
$$
   a_2 \frac{\partial F}{\partial y_1} - y_1 = 0,
$$
$$
   y_1 \frac{\partial F}{\partial y_1} +
       x_2 \frac{\partial F}{\partial x_2}  + y_2 \frac{\partial F}{\partial y_2}- 2 F
 = 0,
\eqno (6.7)
$$
$$
  - a_2 x_1 \frac{\partial F}{\partial y_1} +
   (b_1 x_2 - b_2 y_2 )\frac{\partial F}{\partial x_2} +
   (b_2 x_2 + b_1 y_2 )\frac{\partial F}{\partial y_2} + x_1 y_1
= 0.
$$

\

   Решением первого из этих уравнений является функция
$$
   F = \frac1{2 a_2} y_1^2 + G(x_2,y_2)
$$
с произвольной аналитической функцией $ G(x_2,y_2) $. Два оставшихся уравнения системы (6.7) преобразуются после этого в
$$
       x_2 \frac{\partial G}{\partial x_2}  + y_2 \frac{\partial G}{\partial y_2} - 2 G
 = 0,
\quad
   (b_1 x_2 - b_2 y_2 )\frac{\partial G}{\partial x_2} +
   (b_2 x_2 + b_1 y_2 )\frac{\partial G}{\partial y_2}
= 0.
\eqno (6.8)
$$

  Первое из этих уравнений легко решается, так что (при $ x_2 \ne 0 $)
$$
   G = x_2^2\, \varphi\left(\frac{y_2}{x_2} \right)
$$
 с произвольной аналитической функцией $ \varphi(t) $ одного переменного. Последнее уравнение системы (6.8) примет
после подстановки в него полученных формул для $ F $ и $ G $ следующий вид:
$$
   b_2 (1+t^2) \varphi' + (b_1 - 2 b_2 t) \varphi = 0.
\eqno (6.9)
$$

  Как отмечалось выше, для невырожденных орбит обсуждаемой алгебры выполняется неравенство $ b_2 \ne 0 $. Тогда общее решение
ОДУ (6.9) можно описать формулой
$$
   \varphi(t) =  N (1+t^2) e^{-(b_1/b_2) \arctg t},
\eqno (6.10)
$$
где
$ N $ -- произвольная вещественная константа.

   Объединяя формулу (6.10) с предыдущими обсуждениями, получаем уравнения однородных поверхностей, отвечающие основному
подслучаю  $ b_2 \ne 0 $
$$
   v = -\frac{y_1^2}{2a_2} +  N |z_2|^2 {e}^{m \arg z_2}.
\eqno (6.11)
$$

   Ясно, что при $ N = 0 $ полученное уравнение описывает вырожденную по Леви поверхность. При $ N \ne 0 $
остается заметить, что
$$
  |z_2|^2 {e}^{m \arg z_2} = \left|e^{(1-im/2)\ln z_2} \right|^2,
$$
а потому соответствующая голоморфная замена (дополненная растяжением одной из координат) превратит уравнение (6.11) в
одну из квадрик (1.4).
  Предложение 6.1 доказано полностью.
\hfill{$\Box$}

\

{\bf 6.2. Орбиты семейства алгебр $ \mathfrak g_{5,35} $: второй случай}

\

   Согласно обсуждениям раздела 5, нам нужно рассмотреть алгебры векторных полей в $ \Bbb C^3 $ с базисами вида (5.15), т.е.
$$
\begin{array}{ll}
  \, & e_1 =  e^{\lambda z_2} \left(1, 0, i\varepsilon\right), \\
  \, & e_2 = \left(1, 0, 0\right), \\
  \, & e_3 = \left(0, 0, 1\right), \\
  \, & e_4 = (z_1, 1,  w), \\
  \, & e_5 = (w + A_5 e^{z_2}, B_5,  -z_1 + C_5e^{z_2}), \quad B_5 = ({- \alpha + i \varepsilon})/{\lambda}
\end{array}
\eqno (6.12)
$$
и проинтегрировать их.

{\bf Предложение 6.2.} {\it Если параметры $ \alpha, \lambda $ удовлетворяют неравенству
$$
   (i\varepsilon - \alpha)^2 + \lambda^2 \ne 0,
$$
то голоморфной заменой координат, сохраняющей вид полей $ e_1,e_2,e_3,e_4 $ базиса (6.12), можно привести поле $ e_5 $ к виду
$$
   e_5 = (w , B_5,  -z_1),
$$
свободному от экспоненциальных слагаемых.

  В случае равенства
$
   (i\varepsilon - \alpha)^2 + \lambda^2 =  0,
$
поле $ e_5 $ можно привести к виду
$$
   e_5 = (w , B_5,  -z_1 + C_5^* e^{z_2})
$$
(с некоторым комплексным $ C_5^* $), не изменяя остальных полей базиса (6.12).
}

{\it Доказательство.}
  Рассмотрим изменение базисных полей
(6.12) при голоморфной замене переменных
$$
   z_1^* = z_1 + M e^{z_2}, \quad z_2^* = z_2, \quad  w^* = w + N e^{z_2}
\eqno (6.13)
$$
с некоторыми комплексными константами $ M, N $.  При такой замене дифференцирования по переменным $ z_1, w $
переходят в аналогичные  дифференцирования по новым переменным $ z_1^*, w^* $, а
$$
   \frac{\partial}{\partial z_2} = \frac{\partial}{\partial z_2^*} + M e^{z_2} \frac{\partial}{\partial z_1^*}  +
                                                                     N e^{z_2} \frac{\partial}{\partial w^*}.
\eqno (6.14)
$$

   Поэтому три первых базисных поля очевидным образом сохранят свой вид в новых координатах. Поле $ e_4 $ также
сохранит свой вид, т.к. в силу (6.14),
$$
   e_4 = (z_1, 1, w) = z_1 \frac{\partial}{\partial z_1} + \frac{\partial}{\partial z_2} + w \frac{\partial}{\partial w} =
                     (z_1 + M e^{z_2}) \frac{\partial}{\partial z_1^*} +
\frac{\partial}{\partial z_2^*} +
                     (w + N e^{z_2}) \frac{\partial}{\partial w^*}  = (z_1^*, 1, w^*).
$$

  А поле $ e_5 $ преобразуется по следующим формулам:
$$
   e_5 = \left( (w + A_5 e^{z_2}) + B_5 M e^{z_2}\right) \frac{\partial}{\partial z_1^*} + B_5 \frac{\partial}{\partial z_2^*} +
         \left( (-z_1 + C_5 e^{z_2}) + B_5 N e^{z_2}\right) \frac{\partial}{\partial w^*} =
\quad \quad
$$$$
\quad \quad
 =  \left( w^* + (- N + A_5 + B_5 M) e^{z_2}\right) \frac{\partial}{\partial z_1^*} + B_5 \frac{\partial}{\partial z_2^*} +
         \left( -z_1^* +(M + C_5 + B_5 N) e^{z_2}\right) \frac{\partial}{\partial w^*}.
$$

   Следовательно, при выполнении двух равенств
$$
   - N + A_5 + B_5 M =0, \quad M + C_5 + B_5 N = 0
\eqno (6.15)
$$
поле $ e_5 $ освободится в новых координатах от слагаемых  $ A_5 e^{z_2}, C_5 e^{z_2} $.
Заметим, однако, что подбор констант $ M, N $, обеспечивающих выполнение условий (6.15), возможен не всегда и связан
с определителем системы (6.15)
$$
   \Delta =
\left| \begin{array}{ll}
B_5 & -1 \\
1 & B_5
\end{array}
\right|
= B_5^2 + 1.
\eqno (6.16)
$$

   Если этот определитель не равен нулю, то при наборе $ (M,N) $, являющемся решением системы (6.15), замена (6.13) делает поле $ e_5 $ линейным.

   Если же определитель (6.16) равен нулю, но константа $ A_5 \ne 0 $, мы можем воспользоваться упрощенной заменой (6.13), полагая
$ M = 0, N = A_5 $. Тогда изменение только переменной $ w $ позволяет, в силу формул (6.14) удалить экспоненциальное слагаемое в первой компоненте поля $ e_5 $ и заменить коэффициент $ C_5 $ поля $ e_5 $ на
$ C_5^* = C_5 + B_5 N = C_5 + B_5 A_5 $.

  Остается заметить, что с учетом выражения для
$
  B_5 = ({- \alpha + i \varepsilon})/{\lambda}
$
определитель (6.16) можно записать в виде
$$
    \Delta = \frac{(- \alpha + i \varepsilon)^2 + \lambda^2}{\lambda^2}.
\eqno (6.17)
$$

   Условия предложения 6.2 связаны с равенством или неравенством нулю числителя последней формулы. Тем самым, предложение 6.2 доказано.
\hfill{$\Box$}

{\bf Замечание.} С учетом условия $ \varepsilon = \pm 1 $ числитель формулы (6.17) обращается в нуль лишь при двух наборах
$ (\alpha, \lambda) $, а именно,
$$
    \alpha = 0, \ \lambda = \pm 1.
$$

   Дальнейшее исследование орбит алгебр Ли с базисами вида (6.12) мы разделяем на два подслучая.
  В первом, общем, подслучае считаем базис любой интересующей нас алгебры имеющим вид
$$
\begin{array}{ll}
  \, & e_1 =  e^{\lambda z_2} \left(1, 0, i\varepsilon\right), \\
  \, & e_2 = \left(1, 0, 0\right), \\
  \, & e_3 = \left(0, 0, 1\right), \\
  \, & e_4 = (z_1, 1,  w), \\
  \, & e_5 = (w , B_5,  -z_1 ), \quad B_5 = ({- \alpha + i \varepsilon})/{\lambda}, \ (\alpha,\lambda) \ne (0,\pm1).
\end{array}
\eqno (6.18)
$$

  Во втором, частном, случае имеем базис
$$
\begin{array}{ll}
  \, & e_1 =  e^{\lambda z_2} \left(1, 0, i\varepsilon\right), \\
  \, & e_2 = \left(1, 0, 0\right), \\
  \, & e_3 = \left(0, 0, 1\right), \\
  \, & e_4 = (z_1, 1,  w), \\
  \, & e_5 = (w , i\lambda\varepsilon, -z_1 + C e^{z_2} ), \quad  \varepsilon = \pm 1, \ \lambda = \pm 1, \  C \in \Bbb C.
\end{array}
\eqno ((6.19)
$$

  Нам необходимо проинтегрировать такие алгебры.

{\bf Предложение 6.3.} {\it  С точностью до (локальной) голоморфной эквивалентности орбита любой алгебры с базисом (6.18) в
пространстве $ \Bbb C^3 $ описывается уравнением
$$
    v \sin y_2 - y_1 \cos y_2 = e^{mx_2 + n y_2}, \quad
m\in \Bbb R, \ n \ge 0, \ (m,n)\ne (\pm1,0).
\eqno (6.20)
$$
}

{Доказательство.} Система уравнений в частных производных, отвечающая тройке нетривиальных полей базиса (6.18)
и определяющим функциям $ v = F(y_1, x_2, y_2) $ искомых орбит, имеет вид
$$
   e^{-\lambda x_2} \sin{\lambda y_2}\, \frac{\partial F}{\partial y_1} -
   \varepsilon e^{-\lambda x_2} \cos{\lambda y_2} = 0,
\eqno (6.21)
$$$$
    \frac{\partial F}{\partial y_1} +  \frac{\partial F}{\partial x_2} - F = 0,
\qquad
  F \frac{\partial F}{\partial y_1} + \frac1{\lambda} \left(-\alpha\, \frac{\partial F}{\partial x_2}
         + \varepsilon \,\frac{\partial F}{\partial y_2} \right) + y_1 = 0.
$$

  Освобождаясь в первом уравнении этой системы от экспоненты, получаем его решение в виде
$$
   F = \varepsilon y_1\, \ctg \lambda y_2 + G(x_2, y_2)
\eqno (6.22)
$$
с произвольной аналитической функцией $ G(x_2, y_2) $.

  Подстановка такой функции во второе и третье уравнения (6.21) превращает их в
$$
    \frac{\partial G}{\partial x_2} = G, \quad
        -\alpha\, \frac{\partial G}{\partial x_2}
         + \varepsilon \,\frac{\partial G}{\partial y_2}  = - \varepsilon \lambda\, \ctg \lambda y_2\, G.
$$

   Очередная подстановка решения
$$
   G (x_2, y_2) = e^{x_2} H(y_2)
\eqno (6.23)
$$
первого из этих уравнений во второе приводит (после сокращения на экспоненту) к ОДУ
$$
    \varepsilon H' = (\alpha - \varepsilon \lambda\, \ctg \lambda y_2 )H.
$$

  Связывая его решение
$$
   \ln H = \alpha \varepsilon y_2 - \ln (\sin \lambda y_2) + D,
$$
где $ D $ -- произвольная вещественная константа, с предыдущими шагами, получаем итоговое решение системы (6.21) в виде
$$
   F =  \varepsilon y_1\, \ctg \lambda y_2 - e^{x_2} e^D e^{\alpha \varepsilon y_2}\, \frac1{\sin \lambda y_2}.
$$

   Тогда искомые орбиты алгебр с базисами вида (6.18) задаются уравнениями
$$
   v \sin \lambda y_2 - \varepsilon y_1\, \cos \lambda y_2 = e^D e^{x_2 + \alpha \varepsilon y_2}.
$$

  Сжатие в $ e^D $ раз переменных $ z_1, w $ (и, в т.ч., их мнимых частей) позволяет записать последнее уравнение в более простой форме
$$
     v \sin \lambda y_2 - \varepsilon y_1\, \cos \lambda y_2 = e^{x_2 + \alpha \varepsilon y_2}.
$$

   А растяжение переменной $ z_2 $ в $ \lambda $ раз и введение новых параметров
$$
    m = \frac1{\lambda}, \quad n = \frac{\alpha \varepsilon}{\lambda}
$$
превращает его в
$$
       v \sin y_2 - y_1 \cos y_2 = e^{mx_2 + n y_2}.
$$

   Остается отметить, что этот результат по форме совпадает с доказываемым утверждением. Дополнительное ограничение  $ n \ge 0 $
достигается в силу возможности еще одной несложной замены
$$
    z_1^* = -z_1, \ z_2^* = -z_2,
$$
сохраняющей вид (6.20) и изменяющей знак параметра $ n $.
   Предложение 6.3 доказано.
\hfill{$\Box$}

{\bf Замечание 1.} При $ m = 0 $ уравнение (6.20) описывает известное [4] однопараметрическое семейство однородных
гиперповерхностей в $ \Bbb C^3 $ с 6-мерными алгебрами симметрий.

{\bf Замечание 2.} Домножая уравнение (6.20) на $ e^{x_2} $, можно записать его в виде
$$
    v \left(e^{x_2}\sin y_2\right) - y_1 \left(e^{x_2}\cos y_2\right) = e^{(m+1)x_2 + n y_2}.
$$
С учетом введения новой переменной $ z_2^* = e^{z_2} $ можно далее преобразовать (6.20) к виду
$$
    v y_2 - y_1 x_2 =  |z_2|^A e^{B\arg z_2} \ (A= m+1, \ B = n).
\eqno (6.24)
$$

{\bf Предложение 6.4.} {\it  С точностью до (локальной) голоморфной эквивалентности орбита любой алгебры с базисом (6.19) в
пространстве $ \Bbb C^3 $ описывается одним из двух уравнений:
$$
    v y_2 + y_1 x_2 = 0 \quad \mbox{или} \quad
v y_2 + y_1 x_2 = |z_2|^2 \arg z_2.
\eqno (6.25)
$$
}

{\bf Замечание.} Преобразование $ z_1^* = -z_1 $ позволяет заменить сумму разностью в левых частях уравнений (6.24) и (6.25).

{Доказательство.} Система трех содержательных уравнений, отвечающая алгебре с базисом (6.19), лишь в последнем уравнении
незначительно отличается от системы (6.21). Это уравнение имеет вид
$$
  F \frac{\partial F}{\partial y_1} - {\lambda} \varepsilon \,\frac{\partial F}{\partial y_2}  + y_1 -
        ( c_1 \sin y_2 + c_2\cos y_2)e^{x_2} = 0,
\eqno (6.26)
$$
где $ c_1 = Re\,C, \  c_2 = Im\, C $.

  Подстановка в него формул (6.22) и (6.23) приводит (так же, как и при доказательстве предложения 6.3, после сокращения
на $ e^{x_2} $ ) к обыкновенному дифференциальному уравнению
$$
   H'+ \lambda \ctg({\lambda y_2}) H= - \lambda \varepsilon ( c_1 \sin y_2 + c_2\cos y_2)
\eqno (6.27)
$$
или, с учетом условия $ \lambda = \pm 1 $ и нечетности котангенса, к более простому
$$
    H'+ \ctg({y_2}) H= - \lambda \varepsilon ( c_1 \sin y_2 + c_2\cos y_2).
$$

   Общим решением этого уравнения является функция
$$
   H(y_2) = \frac{\lambda \varepsilon}{4 \sin y_2}\left( c_1 (\sin 2y_2 -2y_2) + c_2\cos2 y_2  \right) + \frac{D}{\sin y_2},
$$
где $ D $ - произвольная константа.

  Итоговое уравнение любой орбиты обсуждаемой в этом случае алгебры имеет вид
$$
    v = - \varepsilon y_1 \ctg({\lambda}y_2) +  e^{x_2}  \frac{\lambda \varepsilon}{4 \sin y_2}\left( c_1 (\sin 2y_2 -2y_2) + c_2\cos2 y_2  \right) + e^{x_2}\frac{D}{\sin y_2}.
$$

  Домножим его на $ e^{x_2} \sin y_2 $ и перегруппируем слагаемые:
$$
    v (e^{x_2}\sin y_2) + \varepsilon y_1 (e^{x_2}\cos({\lambda}y_2)) =
       e^{2x_2}\frac{\lambda \varepsilon}{4}\left( (c_1 \sin 2y_2  + c_2\cos2 y_2) - 2 c_1 y_2 \right) + e^{2x_2}{D}.
$$

  Вводя теперь новые переменные
$
   z_1^* = \varepsilon z_1, \quad z_2^* = e^{z_2}
$
и переобозначая константы,
запишем последнее уравнение в форме
$$
    v y_2^* + y_1 x_2^* = A_1 (e^{2x_2} \sin 2y_2) + A_2 (e^{2x_2} \cos 2y_2) + A_3 e^{2x_2} - 2 A_1 e^{2x_2} y_2.
\eqno (6.28)
$$

  Здесь, например,
$$
    A_1 = \frac{\lambda \varepsilon}{4} c_1,
\eqno (6.29)
$$
и аналогично пересчитываются другие константы.

 Слагаемые из правой части также можно выразить через переменную $ z_2^* $, т.что
$$
   e^{2x_2} \sin 2y_2 = Im\left((z^*)^2\right) = 2 x_2^* y_2^*,\quad
   e^{2x_2} \cos 2y_2 = Re\left((z^*)^2\right) = (x_2^*)^2 - (y_2^*)^2,
$$$$
   e^{2x_2} = |e^{z_2}|^2 = |z_2^*|^2 =  (x_2^*)^2 + (y_2^*)^2, \quad
   e^{2x_2} y_2 = |z_2^*|^2( \arg z_2^*)
$$

  С учетом этого уравнение (6.28) можно теперь записать в виде (звездочки в обозначениях отбрасываем)
$$
    v y_2 + y_1 x_2 = 2 A_1 x_2 y_2 + A_2 (x_2^2 - y_2^2) + A_3 (x_2^2 + y_2^2) - 2 A_1  |z_2|^2 \arg z_2.
\eqno (6.30)
$$

   Все квадратичные слагаемые из правой части этого уравнения можно перенести в левую, после чего это уравнение перепишется как
$$
   (v - 2 A_1 x_2 + A_2 y_2 - A_3 y_2 ) y_2 + ( y_1 - A_2 x_2 - A_3 x_2) x_2 =  - 2 A_1  |z_2|^2 \arg z_2.
\eqno (6.31)
$$

  Линейные формы
$$
    (v - 2 A_1 x_2 + A_2 y_2 - A_3 y_2), \quad (y_1 - A_2 x_2 - A_3 x_2)
$$
из левой части (6.31) можно теперь принять за новые переменные $ v^* = Im\, w^*, \ y_1^* = Im\, z_1^* $, соответственно.

   Тогда после очередного отбрасывания звездочек получим уравнение
$$
   v y_2 + y_1 x_2  = - 2 A_1  |z_2|^2 \arg z_2,
\eqno (6.32)
$$
в котором единственный коэффициент равен, согласно (6.29),
$
   - 2 A_1 = - {\lambda \varepsilon}c_1/{2}.
$

   Если при этом $ c_1 = 0 $, то уравнение (6.32) совпадает с первой формулой (6.25). Если же $ c_1 \ne 0 $, то
после еще одного растяжения координат
$$
    z_1^* = \frac{z_1}{-2 A_1},  \quad w^* = \frac{w}{-2 A_1}
$$
мы получим вторую из формул (6.25).
  Предложение 6.4 доказано.
\hfill{$\Box$}

  {\bf Замечание.} При $ \lambda = 1 $ все элементы базиса (6.19) можно превратить голоморфной заменой $ z_2^* = e^{z_2} $
в аффинные векторные поля.

\

{\bf 6.3. Проверка гипотезы о новизне орбит алгебры $ \mathfrak g_{5,35} $}

\

  Нам необходимо теперь обсудить вопросы невырожденности по Леви,
сферичности и новизны полученных однородных гиперповерхностей (9.9) и (9.13).

Отметим, что
в целом задача проверки (локальной) голоморфной эквивалентности двух произвольных заданных поверхностей, а тем более,
построения биголоморфного отображения одной из них на другую, является достаточно сложной.
Вместе с тем, проверка сферичности конкретных гиперповерхностей в $ \Bbb C^3 $ является вполне конструктивной задачей.
Полезным инструментом при этом (как и при проверке возможной голоморфной эквивалентности пары гиперповерхностей в $ \Bbb C^n $ и, в частности, в $ \Bbb C^3 $), является, как отмечалось выше, нормальная форма Мозера уравнений таких поверхностей.

   Например, в соответствии с известным фактом многомерного комплексного анализа (см. [23]), омбилическая
в каждой своей точке вещественно-аналитическая гиперповерхность --- сферична. В свою очередь, омбиличность аналитической
гиперповерхности в $ \Bbb C^3 $,
означающая обращение в нуль многочлена $ N_{220} $ из нормального уравнения такой поверхности,
 может быть проверена чисто техническими средствами. Достаточно подробно такая техника описана в [23] и [30],
а потому здесь мы лишь кратко прокомментируем необходимые действия и шаги нормализации.

  Отметим также, что имеются работы (см. [21,44-46]), в которых
свойство сферичности поверхностей (как в СПВ-случае, так и в индефинитной ситуации), устанавливается за счет их других геометрических
или аналитических характеристик. Укажем здесь, в первую очередь, статью [45], ссылка на которую позволяет достаточно легко
решить часть интересующих нас конкретных вопросов.

{\bf Предложение 6.5.} {\it Обе поверхности (6.25) являются индефинитными сферическими.
}

{ Доказательство.}
  В списке [45] имеется сферическая (индефинитная) трубка
$$
    v = y_1 \tg(y_2).
\eqno (6.33)
$$
   Первая из поверхностей (6.25), т.е.
$$
   v y_2 + y_1 x_2  = 0 \ \mbox {или} \quad y_1 = -v \frac{y_2}{x_2}
$$
сводится к названной трубке Исаева-Мищенко голоморфным преобразованием  $ z_2^* = \ln {z_2} $. В самом деле,
$$
    \frac{y_2}{x_2} = \tg(\arg (x_2 + i y_2)) = \tg(Im (z_2^*)) = \tg(y_2^*),
$$
а потому при таком преобразовании мы получаем новое уравнение проверяемой поверхности (9.13) в виде
$
   y_1 = - v \tg(y_2^*).
$
  Еще одной заменой координат $ z_1 = w^*, \ w = -z_1^* $ последнее уравнение переводится в (6.33).

  Переходя к обсуждению второй поверхности (6.25), мы воспользуемся техникой нормальных форм и соответствующими
формулами работы [30].

  Удобно для этого записать исходное уравнение в измененной форме
$$
   v x_2 - y_1 y_2 = |z_2|^2 \arg z_2,
\eqno (6.34)
$$
а в качестве точки, в которой будет осуществляться нормализация уравнения, возьмем точку $ Q(0,1,0) $ пространства $ \Bbb C^3 $.
После сдвига в эту точку $ x_2 \rightarrow 1+ x_2 $
получаем уравнение поверхности в виде
$$
    v = \frac1{1+x_2}\left(y_1 y_2 +\left((1+x_2)^2+y_2^2\right)\arctg{\frac{y_2}{1+x_2}} \right)
\eqno (6.35)
$$

  Выписывая разложение правой части уравнения (6.35) в ряд Тейлора до 4-й степени включительно (это можно сделать <<руками>>, но надежнее
воспользоваться пакетом символьной математики типа Maple), получаем
$$
   v = \sum_{k=1}^{\infty} F_k (y_1, x_2, y_2),
\ \mbox{
где}
\ \
   F_1 = y_2, \
   F_2 = y_1 y_2, \
    F_3 = - y_1 x_2 y_2 + \frac23y_2^3, \
    F_4 = y_1 x_2^2 y_2 - \frac43 x_2 y_2^3.
\eqno (6.36)
$$

   Далее от вещественных переменных нужно перейти к комплексным. Чтобы не загромождать выкладки дробями, удобно сделать это по следующим формулам
$$
    x_2 \rightarrow (z_2 + \bar z_2), \quad y_k \rightarrow i(\bar z_k - z_k), \ k = 1,2,
$$
имея в виду дополнительное растяжение координат
$$
    z_k = 2 z_k^* = 2 ( x_k^* + i y_k^*), \quad k = 1,2.
$$

  Эрмитову часть получающегося после такой замены слагаемого
$$
   F_2 = (z_1\bar z_2 + z_2 \bar z_1) - (z_1 z_2 + \bar z_1 \bar z_2),
$$
т.е. (индефинитную) форму Леви обсуждаемой поверхности, обозначим через
$$
   <z,z> = z_1\bar z_2 + z_2 \bar z_1.
$$

   Все другие слагаемые из уравнения (6.36), выраженные через комплексные переменные, распадаются на группы мономов одинаковых
бистепеней по переменным $ z=(z_1,z_2) $ и $ \bar z = (\bar z_1, \bar z_2) $. Например,
$$
   F_3(z,\bar z) = F_{3\bar 0}+ F_{2\bar 1} + F_{1\bar 2}+ F_{0\bar 3},
$$
где $ F_{3\bar 0} $ означает группу голоморфных слагаемых, $ F_{2\bar 1} $ -- группу слагаемых бистепени $(2,\bar 1) $ и т.п.

   При нормализации уравнения (6.36) оно примет вид
$$
     v = <z,z> + N_{2\bar20} + ...,
$$
где многоточиями обозначены слагаемые более высоких степеней по переменным $ z, \bar z, u = Re\,w $.
 При этом многочлен $ N_{2\bar20} $ равен проекции в $ \mathfrak N $-пространство Мозера многочлена
$$
    H_{2\bar2} = F_{2\bar2} - <f_2,f_2>.
$$

  Здесь вектор-функция $ f_2 = (f_2^{(1)}, f_2^{(2)}) $ определяется (однозначно) из соотношения
$$
     F_{2\bar 1}(z,\bar z) = < f_2, z > = f_2^{(1)} \bar z_2 + f_2^{(2)} \bar z_1.
$$

   Разворачивая $(2,\bar 1) $- и $(2,\bar 2) $-компоненты слагаемых $ F_3, F_4 $ из уравнения (6.36), получим
$$
   F_{2\bar 1} = -z_2^2 \bar z_1 - 2i z_2^2 \bar z_2, \quad
   F_{2\bar 2} = z_2^2 \bar z_1 \bar z_1 + z_1 z_2 \bar z_2^2.
\eqno (6.37)
$$

  Это означает, что вектор-функция $ f_2 = (- 2i z_2^2, -z_2^2) $, а выражение
$
   < f_2, f_2 > = f_2^{(1)} \overline{f_2^{(2)}} + f_2^{(2)} \overline{f_2^{(1)}}
$
равно нулю.
   Остается заметить, что в соответствии с формулами из [30, Предложение 1.2],  проекция многочлена
$$
    (\bar C z_2^2 \bar z_1 \bar z_1 +  C z_1 z_2 \bar z_2^2)
$$
в $ \mathfrak N $-пространство Мозера равна
$$
    i\, Im\, C ( z_1 z_2 \bar z_2^2 - z_2^2 \bar z_1 \bar z_1 ).
$$

  Так как коэффициенты многочлена $ F_{2\bar 2} $ из формулы (6.37) вещественны, соответствующий многочлен $ N_{2\bar 20} $ из
нормального уравнения поверхности (6.34) равен нулю. Тем самым, эта поверхность омбилична, а в силу ее однородности -- сферична.
Предложение 6.5 доказано.
\hfill{$\Box$}

   Заключительным результатом этого подраздела является утверждение о
размерности алгебр симметрии невырожденных
несферических орбит 5-мерной алгебры $ \mathfrak g_{5,35} $, описываемых уравнениями (6.20), т.е.
$$
    v \sin y_2 - y_1 \cos y_2 = e^{mx_2 + n y_2}, \quad
m\in \Bbb R, \ n \ge 0, \ (m,n)\ne (\pm1,0).
$$

{\bf Предложение 6.6.} {\it На любой гиперповерхности (6.20)
имеется 6-мерные алгебры голоморфных векторных полей. Все эти поверхности имеют индефинитную форму Леви
}

   Доказательство.
   Знаконеопределенный характер формы Леви для всех поверхностей (6.20) устанавливается простым вычислением.

   Далее, по аналогии с доказательством предложения 3.7 мы дополним базис (6.18)
$$
   e_1 =  e^{\lambda z_2} \left(1, 0, i\varepsilon\right), \ \
   e_2 = \left(1, 0, 0\right), \ \
   e_3 = \left(0, 0, 1\right), \ \
   e_4 = (z_1, 1,  w), \ \
  e_5 = (w , B_5,  -z_1 ),
$$$$
 \quad B_5 = ({- \alpha + i \varepsilon})/{\lambda}, \quad  \lambda = 1 - \beta, \ (\alpha,\lambda) \ne (0,\pm1)
$$
(возникающий при реализации в рамках второго случая любой из алгебр
$ \mathfrak g_{5,35} $) шестым голоморфным векторным полем
$$
    e_6 = e^{-\lambda z_2} (1,0,-i).
$$

   Непосредственными вычислениями проверяется, что это поле является касательным к поверхности (6.20) при
любых $ m, n $. Кроме того, для поля $ e_6 $ и выписанного базисного набора полей (6.12) при при любых $ \alpha, \beta $ вещественная линейная оболочка
$
   <e_1, e_2,e_3,e_4,e_5,e_6>
$
образует 6-мерную алгебру векторных полей. Структура этой алгебры следующая:
$$
   [e_1,e_4] = \beta e_1, \quad [e_2,e_4] = e_2, \quad [e_3,e_4] = e_3; \quad [e_4,e_6] = (\beta -2) e_6,
\eqno (6.38)
$$
$$
   [e_1,e_5] = \alpha e_1, \quad [e_2,e_5] = - e_3, \quad [e_3,e_5] = e_2; \quad [e_5,e_6] = \alpha e_6.
$$

  Предложение 6.6 доказано.
\hfill{$\Box$}

{\it Следствие.} Поверхности (6.20) НЕ являются новыми в задаче описания голоморфно однородных гиперповерхностей в $ \Bbb C^3 $.

   Вопрос об известных однородных гиперповерхностях, которым эквивалентны поверхности (6.20), мы здесь не обсуждаем, как и
аналогичный вопрос о поверхностях (3.10), затронутый в разделе 3.
Отметим лишь, что алгебры (6.38) изоморфны 6-мерным алгебрам, ассоциированным с семейством однородных гиперповерхностей (см. [4])
$$
   v = {i}(z_1\bar z_2 - z_2 \bar z_1) + (z_1^A)\overline{(z_1^A)}, \ A \in \Bbb C \setminus \{-1,0,1,2\}.
\eqno (6.39)
$$

  Вместе с тем для получения точного ответа на обозначенный вопрос необходимо, вообще говоря, проверить возможность вложения
алгебры (6.38) в 6-мерные и 7-мерные алгебры векторных полей, отвечающие другим известным однородным гиперповерхностям.

\

{\bf 6.4. Итоговые выводы о <<просто однородных>> гиперповерхностях}

\

  После отклонения гипотезы о новизне поверхностей (6.20) остается подвести окончательные итоги
в <<просто однородной>> части обсуждаемой в статье задачи. Приведем здесь список просто однородных
Леви-невырожденных гиперповерхностей
пространства $ \Bbb C^3 $ (из Теоремы 2.5)
c указанием структур соответствующих алгебр Ли
голоморфных векторных полей (из классификации [5]).

$$
\begin{array}{|c|c|c|c|c|c|c|c|c|c|c|c|}
\hline
N & \mbox{ Поверхности } & \mbox{ Алгебры Ли } & \mbox{Ограничения } \\ \hline
1 & v = x_1^{\alpha} x_2^{\beta} & \mathfrak{g}_{5,33} & 0 < |\alpha|\le |\beta|\le 1, \alpha + \beta \ne 1  \\ \hline
2 & v = (x_1^2 + x_2^2)^{\alpha} \exp(\beta \arg(x_1+ix_2)) & \mathfrak{g}_{5,35} & \alpha \ne 1/2, \beta \ge 0, (\alpha,\beta) \ne (1,0)
 \\ \hline
3 & v = x_1({\alpha}\ln x_1 + \ln x_2) & \mathfrak{g}_{5,34} & {\alpha} \notin \{-1,0\}  \\ \hline
4 & ( v - x_1 x_2 + x_1^3/3 )^2 = \beta (x_2 - x_1^2/2)^3 & \mathfrak{g}_{5,30} (\alpha = 0) & \beta \notin \{0,4\}  \\ \hline
5 & x_1 v = x_2^2 \pm x_1^{\alpha} & \mathfrak{g}_{5,30}  & \alpha \notin\{0,1,2\}   \\ \hline
6 & x_1 v = x_2^2 \pm x_1^2 \ln x_1 & \mathfrak{g}_{5,32}(\alpha = \pm 1) & \,  \\ \hline
7 & v (1\pm x_2 y_2) = y_1 y_2 & \mathfrak{g}_{5,32}(\alpha = 0) & \, \\ \hline
8 & (v - x_2 y_1)^2+ y_1^2 y_2^2 = y_1 & \mathfrak{g}_{5} & \, \\ \hline
\end{array}
$$
\

{\bf Таблица 7.} <<Просто однородные>> невырожденные  гиперповерхности в $ \Bbb C^3 $
\

{\bf Замечание.}   Приведем несколько простых комментариев к этой таблице:

--- при $ \alpha = 0 $, как и при $ \alpha + \beta = 1 $, поверхности из строки 1 вырождены по Леви (см. [3]);

--- при $ \alpha = 1/2 $ поверхности из строки 2 вырождены по Леви (см. [3]), а при $ \alpha = 1, \ \beta = 0 $ получаем здесь
трубчатую СПВ-квадрику
 $ v = x_1^2 + x_2^2 $ c 15-мерной алгеброй симметрий в $ \Bbb C^3 $;

--- при значении $ \alpha = -1 $ поверхность из строки 3 является вырожденной, а при $ \alpha = 0 $ соответствующая поверхность
эквивалентна индефинитной квадрике;

--- поверхность из строки 4 вырождается при $ \beta = 4 $, а при $ \beta = 0 $ соответствующая поверхность
эквивалентна индефинитной квадрике;

--- поверхности из строки 5 превращаются в вырожденные трубки над конусами при выделенных в таблице значениях $ \alpha = 1, 2 $, а
при $ \alpha = 0 $ --- в гиперболоиды.

 Теперь мы прокомментируем эту таблицу более подробно и, тем самым, приведем доказательство Теоремы 2.5.

  Так, доказательства простой однородности и новизны двух последних типов поверхностей
$$
     v\, ( 1 \pm y_2 x_2) = y_1 y_2
\quad \mbox{и} \quad
   \ (v - x_2 y_1)^2 + y_1^2 y_2^2 = y_1
$$
из таблицы 7,
использующие технику нормальных форм, приведены, соответственно, в статьях [7] и [8]. Уточним здесь, что
индефинитные поверхности
$
     v\, ( 1 \pm y_2 x_2) = y_1 y_2
$
имеют нормальные формы антитрубчатого типа 6) в терминологии Предложения 1.2 настоящей статьи. Поверхность
$
   \ (v - x_2 y_1)^2 + y_1^2 y_2^2 = y_1
$
относится к общему типу 1) этого же предложения.

   Оба эти типа не сводятся голоморфными преобразованиями к трубчатым многообразиям, в отличие от всех остальных
поверхностей из Теоремы 2.5.

  Факт такой сводимости обеспечивается (голоморфно инвариантным) свойством наличия 3-мерного абелева идеала
в алгебре симметрий обсуждаемой поверхности $ M $. Используя этот факт и обсуждение коммутационных свойств в алгебрах
из блока VI, можно получить вывод о голоморфно разной структуре поверхностей из строк 1 и 2 таблицы 7, а также аналогичный вывод
о строках 1 и 3.

   Для получения более общих утверждений о поверхностях из таблицы 7, можно использовать технику нормальных форм, а также
следующее утверждение о трубчатых гиперповерхностях в $ \Bbb C^3 $ из работы  [33].

{\bf Предложение 6.7 ([33, Proposition 4.1]).}  {\it Пусть $ M_1, M_2 \subset \Bbb C^ 3 $ - две трубчатые гиперповерхности над аффинно однородными
основаниями $ B_1, B_2 $, соответственно. Предположим, кроме того, что $ M_1, M_2 $ просто однородны и что абелев
идеал, порождаемый сдвигами $ i {\partial}/{\partial z_k} \ (k = 1,2,3) $ вдоль мнимых осей,
является единственным 3-мерным абелевым идеалом в обеих алгебрах $ g_1(M_1) $ и $ g_2(M_2) $
голоморфных векторных полей на $ M_1, M_2 $.
Тогда биголоморфная эквивалентность $ M_1 $ и $ M_2 $ в некоторых их точках равносильна аффинной эквивалентности их оснований.
}

  В работе [33] это утверждение использовано при описании однородных СПВ-гиперповерхностей, однако
формулировка и сам результат остаются справедливыми и для индефинитных Леви-невырожденных гиперповерхностей.

   В силу этого, для построения полного списка голоморфно не эквивалентных <<просто однородных>>
трубок в $ \Bbb C^3 $ над аффинно-однородными основаниями достаточно иметь полный список
{\it аффинно не эквивалентных} аффинно однородных поверхностей в $ \Bbb R^3 $,
порождающих эти <<просто однородные>> трубки. Список однородных поверхностей из работы [11], приведенный
в теореме 2.1 настоящей статьи, является полным (по крайней мере, признан таковым многими специалистами).

Вместе с этим,
он является <<переполненным>>: в нем имеются формально разные поверхности (например, отвечающие разным значениям параметров,
описывающих семейства однородных поверхностей), являющиеся, тем не менее, аффинно эквивалентными.
   Очевидный пример такого типа представляет первое семейство из списка [11]
$
    x_3 = x_1^{\alpha} x_2^{\beta},
$
для которого, например, пары
 $ (\alpha, \beta) $, $ (\beta, \alpha) $, $ ( 1/{\alpha}, -\beta/\alpha ) $ задают аффинно эквивалентные поверхности.

   Ограничения на параметры подобных семейств, которые нужно ввести при работе с ними, авторы [11] назвали <<очевидными>>.
Частично соглашаясь с этим, отметим, что в
теореме 2.5 выписаны (см. также [33]) условия на параметры $ {\alpha}, {\beta} $ этого семейства, гарантирующие взаимно-однозначное
описание его аффинно различных представителей в параметрической форме. Аналогичные уточнения, связанные с этим
семейством (и с группой 6-го порядка, действующей на множестве параметров и приводящей к аффинно эквивалентным поверхностям),
имеются в [47].

   Еще одно уточнение о единственности набора параметров, определяющих аффинно однородные поверхности в $ \Bbb R^3 $,
относится ко второму семейству
$$
    v = (x_1^2 + x_2^2)^{\alpha} e^{\beta arg (x_1 + ix_2)}
\eqno (6.40)
$$
из Теоремы 2.5, приведенному в [11] без каких-либо комментариев о параметрах.
  За счет замены $ x_2^* = - x_2 $ поверхности из $ \Bbb R^3 $, отвечающие парам $ (\alpha, \beta) $  и $ (\alpha, -\beta) $, являются
аффинно эквивалентными. Как следствие, голоморфно эквивалентны трубки над поверхностями с такими наборами параметров.

Отметим, что 5-мерные алгебры голоморфных
векторных полей над трубками (6.40) изоморфны $ \mathfrak g_{5,35} $. При этом замена $ e_3^* = -e_3, \ e_5^* = -e_5 $
в коммутационных соотношениях
$$
 [ e_1, e_4] = \beta e_1, \ [ e_2, e_4] =  e_2, \ [ e_3, e_4] = e_3,
\eqno (6.41)
$$
$$
  [ e_1, e_5] = \alpha e_1, [ e_2, e_5] = - e_3, \ [ e_3, e_5] = e_2,
$$
для этих алгебр сохраняет вид этих соотношений, но изменяет знак у параметра $ \alpha $.

  Такая возможность изменения знака у одного из параметров алгебр $ \mathfrak g_{5,35} $ в описании [5] не учтена.
Остался не замеченным этот факт и в книге [38], авторы которой исправили другую неточность [5] в описании семейства алгебр
$ \mathfrak g_{5,32} $.

   Других значимых замечаний и уточнений к спискам [11] и [5] автором настоящей статьи (с соавторами последних
публикаций об однородности) не обнаружено.
   Аффинные алгебры Ли, обеспечивающие аффинную однородность оснований всех остальных поверхностей из таблицы 6, как показано в [4], отличаются друг от друга. Тем самым, все поверхности из этой таблицы попарно голоморфно не эквивалентны.
К формулировке Теоремы 2.5, несколько отличающейся от непосредственного цитирования работы [11], привел учет обозначенных выше неточностей.

\newpage

\centerline{\bf Приложение: общий список однородных гиперповерхностей в $ \Bbb C^3 $}

\

{\bf I. <<Кратно-транзитивные>> Леви-невырожденные поверхности}

\

\quad {\bf I.1. Трубки с аффинно-однородными основаниями}

\

1) $ x_3 = \ln x_1 + \alpha \ln x_2,  \ \alpha \in [-1,1]\setminus \{0\} $,

2) $ x_3 = \alpha \arg (x_1 + i x_2 ) + \ln (x_1^2 + x_2^2), \ \alpha \ge 0 $,

3) $  x_3 = x_2^2 \pm x_1^{\alpha}, \ \alpha\notin\{ 0,1\}  $,

4) $  x_3 = x_2^2 \pm \ln x_1 $,

5) $  x_3 = y^2 \pm x \ln x $,

6) $  x_3 = x_1 x_2 + e^{x_1} $,

7) $  x_3 = x_1 x_2 +  x_1^{\alpha} $,

8) $  x_3 = x_1 x_2 + \ln x_1 $,

9) $  x_3 = x_1 x_2 + x_1 \ln x_1 $,

10) $  x_3 = x_1 x_2 + x_1^2 \ln x_1 $,

11) $ x_1 x_3 = x_2^2 \pm x_1 \ln x_1 $,

12) $ \varepsilon_1 x_1^2 +\varepsilon_2 x_2^2 + x_3^2 = 1, (\varepsilon_1, \varepsilon_2 \in \{\pm 1\}) $.

\

\quad {\bf I.2. Трубки, основания которых не являются аффинно однородными}

\

$
  1) \  x_3 = (1 + e^{2x_1}) (x_2 + \ln(1 + e^{2x_1})),
$

$
  2) \  x_3 = x_1 x_2 + x_1^3 \ln x_1,
$

$
  3) \ x_3 = x_2 e^{x_1} + e^{\alpha x_1}, \ \alpha \notin\{-1,0,1,2 \}, \ \alpha \sim (1-\alpha),
$

$
  4) \ x_3 \cos x_1 + x_2 \sin x_1 = e^{\alpha x_1}, \ \alpha \ge 0,
$

$
  5) \ x_3 = x_1 e^{x_2} + x_1^2,
$

$
   6) \  x_3 = \alpha \ln(1 + e^{2x_1}) + \ln x_2, \ \alpha \in [-1,1]\setminus \{0\},
$

$
   7) \ x_3 = \alpha \ln(1 + e^{2x_1}) + \ln(1 + e^{2x_2}), \ \alpha \in [-1,1]\setminus \{0\},
$

$
   8) \ x_3 = x_2^2 + \varepsilon \ln(1 + e^{2x_1}).
$

\

\quad {\bf I.3. Поверхности Картанова типа}

\

$
  1)  \ 1 + |z_1|^2 + |z_2|^2 + |z_3|^2 = a | 1 + z_1^2 + z_2^2 + z_3^2 | \ \ ( a > 1 ),
$

$
  2) \ 1 + |z_1|^2 + |z_2|^2  - |z_3|^2 = a | 1 + z_1^2 + z_2^2 - z_3^2 | \ \ (0 < a < 1, \ a > 1 ),
$

$
  3)  \ 1 - |z_1|^2 - |z_2|^2  - |z_3|^2  = a | 1 - z_1^2 - z_2^2 - z_3^2| \ \ (0 < a < 1 ),
$

$
  4)  \ 1 + |z_1|^2 - |z_2|^2  - |z_3|^2 = a | 1 + z_1^2 - z_2^2 - z_3^2 | \ \ (0 < a < 1, \ a > 1 ),
$

\

\quad {\bf I.3. Кватернионные модели однородных поверхностей}

$$
   Im( \bar \xi_1 \xi_3 + \bar \xi_2 \xi_4) = \gamma \sqrt{ |Re(\bar \xi_1 \xi_3 + \bar \xi_2 \xi_4)|^2 + |\xi_1 \xi_4 - \xi_2 \xi_3|^2},
\quad \gamma \in \Bbb R\setminus\{0\}.
$$

\

\quad {\bf I.4. Двухпараметрическое семейство кратно-транзитивных поверхностей Винкельманнова типа}

$ v = (z_1 \bar z_2 + z_2 \bar z_1) + z_1^A\overline{z_1^A}, \ A \in \Bbb C \setminus \{-1,0,1,2\}  $

\

{\bf II. <<Просто однородные>> Леви-невырожденные поверхности}

\

{\bf II.1. <<Просто однородные>> трубки с аффинно однородными основаниями}

\

1) $ x_3 = x_1^{\alpha} x_2^{\beta}, \ 0 < |\alpha|\le |\beta|\le 1, \alpha + \beta \ne 1 $,

2) $  x_3 = (x_1^2 + x_2^2)^{\alpha} e^{\beta arg (x_1 + ix_2)}, \ \alpha \ne 1/2, \beta \ge 0, (\alpha,\beta) \ne (1,0) $,

3) $  x_3 = x_1( \alpha \ln x_1 + \ln x_2 ), \ {\alpha} \notin \{-1,0\}  $,

4) $ (x_3 - 3 x_1 x_2 + 2 x_1^3)^2 = \alpha(x_1^2 - x_2)^3, \ \alpha \notin \{0,4\}, $

5) $ x_1 x_3 = x_2^2 \pm x_1^{\alpha}, \ \alpha \notin\{0,1,2\} $,

6) $ x_1 x_3 = x_2^2 \pm x_1^2 \ln x_1 $.

\

{\bf II.2. Не сводимые к трубкам <<просто однородные>> гиперповерхности}

\

$
   1) \  v\, ( 1 \pm y_2 x_2) = y_1 y_2.
$

$
  2) \ (v - x_2 y_1)^2 + y_1^2 y_2^2 = y_1.
$

\

{\bf III. Вырожденные по Леви однородные поверхности}

\

\quad {\bf III.1. Трубки над поверхностями нулевой гауссовой кривизны  }

\quad \quad \quad \quad
($ \ x_1 = Re\, z_1, x_2 = Re\,z_2, x_3 = Re\,z_3 $)

\

$
 1) \
 x_1^2 + x_2^2 = x_3^2 \ (x_3 > 0),
$

$
 2) \ x_3 = |x_1+ix_2| e^{ \omega \arg (x_1 + i x_2)} \ (\omega > 0)
$

$  3) \ x_3 = x_1 (\ln x_2 - \ln x_1), \
$

$
  4) \ x_3 = x^{1 - \theta} x_2^{\theta} \ (x_1 > 0, x_2 > 0),
$

$
   5) \ (x_3 - 3 x_1 x_2 + 2 x_1^3)^2 = 4(x_1^2 - x_2)^3.
$

\

\quad {\bf III.2. Произведения комплексной прямой $ \Bbb C $ на Леви-невырожденные однородные
вещественные гиперповерхности из $ \Bbb C^2_{z_1, z_2} $} ($ \ x_1 = Re\, z_1, x_2 = Re\,z_2 $)

\

$
   1) \ x_2 = x_1^s \ ( s \in (-1,0) \cup (0,1/2) \cup(1/2,1) ),
$

$
  2) \ x_2 = \ln x_1,
$

$
   3) \ x_2 = x_1 \ln x_1,
$

$
    4) \ r = e^{a\varphi} \
   (  a \ge 0, \ r-\mbox{полярный радиус},\ \varphi - \mbox{полярный угол в плоскости} \ \Bbb R^2_{x_1,x_2}),
$

$
  5)  \ 1 + |z_1|^2 + |z_2|^2 = a | 1 + z_1^2 + z_2^2 | \ \ ( a > 1 ),
$

$
  6) \ 1 + |z_1|^2 - |z_2|^2 = a | 1 + z_1^2 - z_2^2 | \ \ ( a > 1 ),
$

$
   7) \  |z_1|^2 + |z_2|^2 - 1 = a | z_1^2 + z_2^2 - 1 | \ \ ( 0 < | a | < 1 ).
$

\

\quad {\bf III.3. Вещественная гиперплоскость в $ \Bbb C^3 $}

$
   \ v = 0.
$

{\bf Примечание 1.} Всего в приведенном списке имеется 47 различных типов однородных гиперповерхностей.

{\bf Примечание 2.} Ограничения на параметры в разделах I и II приведенного списка объясняются двумя
причинами: 1) вырождение поверхностей при некоторых значениях <<формальных>> параметров,  2) голоморфная эквивалентность поверхностей, отвечающих различным значениям <<формальных>> параметров.

   Так, уравнение $ x_3 = \ln x_1 + \alpha \ln x_2 $ задает однородные гиперповерхности при
любых (а не только отвечающих ограничениям семейства I.1.1) вещественных значениях параметра $ \alpha $. Однако
при $ \alpha = 0 $
соответствующая поверхность вырождена по Леви (см. п. III.2.2.)  Кроме того, любая поверхность из семейства I.1.1
голоморфно эквивалентна поверхности с аналогичным уравнением при замене параметра $  \alpha $ на значение $ 1/\alpha $,
не удовлетворяющее ограничениям этого семейства.

{\bf Примечание 3.}
 <<Стандартные>> квадрики $ v = |z_1|^2 \pm |z_2|^2 $ с 15 мерными алгебрами симметрий голоморфно эквивалентны многим поверхностям, упомянутым в этих списках. Например, каждая из двух трубчатых поверхностей $ x_3 = x_2^2 \pm x_1^2 $ (из п.I.1.9) эквивалентна, в зависимости от знака, одной из этих стандартных квадрик. Аналогично, трубка $ x_3 = x_1 x_2 $ (п. II.1.1, $ \alpha = \beta = 1 $) также
голоморфно эквивалентна индефинитной квадрике  $ v = |z_1|^2 - |z_2|^2 $ и т.п. По этой же причине трубки над поверхностями
$  x_3 = y^2 \pm e^x $ и $ x_3 = \arg(x_1 +i x_2) $ (п. 8 и 4 классификации [11]), также эквивалентные стандартным квадрикам, не вошли в итоговые списки однородных гиперповерхностей.

\newpage

\begin{center}
СПИСОК ЛИТЕРАТУРЫ
\end{center}

\

1. \textit{Poincare H.} Les fonctions analytiques de deux variables et la
         representation conforme // Rend. Circ. Math. Palermo (1907),
         P. 185 - 220.

2.   \textit{Cartan E.} Sur la geometrie pseudoconforme des hypersurfaces de deux variables complexes // Ann. Math. Pura Appl.
      1932. V. 11, N 4. P. 17--90.

3.  \textit{Fels G., Kaup W.} Classification of Levi degenerate homogeneous CR-manifolds in dimension 5 // Acta Math. 2008. V. 201. P. 1--82.

4. \textit{ B. Doubrov, A. Medvedev, D.The.} Homogeneous Levi non-dege\-nerate hypersurfaces in $ \Bbb C^3 $
//arXiv:1711.02389v1 [math.DG] 7 Nov. 2017.

5. \textit{Мубаракзянов~Г.~М.} Классификация вещественных структур алгебр Ли пятого порядка// Изв. вузов. Матем. --- 1963. --- \No~3. --- С.~99--106.

6.  \textit{Акопян Р.С., Лобода А.В.} “О голоморфных реализациях нильпотентных алгебр Ли”,
 Функц. анализ и его прил., 53:2 (2019),  59–63

7. \textit{Р. С. Акопян, А. В. Лобода.} “О голоморфных реализациях пятимерных алгебр Ли”, Алгебра и анализ, 31:6 (2019),  1–37

8. \textit{Атанов А.В., Лобода А.В.} Об орбитах одной неразрешимой 5-мерной алгебры Ли // Математическая физика и компьютерное моделирование, ВолГУ. Т.23, N 2. с.5-20.

9. \textit{Атанов А.В., Коссовский И.Г., Лобода А.В.}
Об орбитах действий 5-мерных неразрешимых алгебр Ли
в трехмерном комплексном пространстве. ДАН, т. 487, No 6, с. 607–610.

10. \textit{Атанов А. В., Лобода  А. В.} Разложимые пятимерные алгебры Ли в задаче о голоморфной однородности в $ \Bbb C^3 $, Итоги науки и техн. Сер. Соврем. мат. и ее прил. Темат. обз., 173 (2019),  86–115

11.   \textit{Doubrov B., Komrakov B., Rabinovich M.} Homogeneous surfaces in the three-dimensional affine geometry // Geometry and topology of submanifolds. VIII / Eds. F. Dillen et al. River Edge (NJ): World Sci., 1996. P. 168--178.

12.  {Beloshapka V.K., Kossovskiy I.G.} Homogeneous hypersurfaces in $\mathbb{C}^3$, associated with a model CR-cubic // J. Geom. Anal. 2010. V. 20, N 3. P. 538--564.

13. \textit{Дубровин Б.А., Новиков С.П., Фоменко А.Т.} Современная геометрия. М.: Наука. 1986, 760 с.

14. \textit{Zaitsev D.} On Different Notions of Homogeneity for CR-Manifolds. Asian J. Math. V. 11, Number 2 (2007), 331-340.

15. \textit{Azad H., Huckleberry A., Richthoffer W.} Homogeneous CR manifolds. Journal Reine Angew. Math. 1985, V. 358. P. 125 - 154.

16. \textit{Широков А.П., Широков П.А.} Аффинная дифференциальная геометрия. М.
       Физматгиз, 1959, 319 с.

17.  \textit{Guggenheimer H.} Differential geometry, McGraw-Hill, New York, 1963.

18. \textit{Nomizu K., Sasaki T.} A new model of unimodular-affinely
         homo\-geneous surfaces // Manuscr. Math. 1991. 73, N 1. P. 39 - 44.

19. \textit{Tanaka N.} On the pseudo-conformal geometry of hypersurfaces of the space of n complex variables.
     J. Math. Soc. Japan 14 (1962), no. 4, 397--429.

20.  \textit{Morimoto A, Nagano T.} On pseudo-conformal transformations of hypersurfaces // J. Math. Soc. Japan. 1963. V. 15, N 3. P. 289--300.

21.  \textit{ Burns D., Shneider S.} Spherical hypersurfaces in complex space //
         Inv. Math. 1976, V. 33, N 3, P. 283 - 289.

22.  \textit{ Baouendi M.S., Rotshild L.P.} Geometric properties of mappings between
       hypersurfaces in complex space // J. Diff. Geom. 1990, V. 31. N 2.
       P. 473 - 499.

23.   \textit{Chern S.S. Moser J.K.} Real hypersurfaces in complex manifolds // Acta Math. 1974. V. 133. P. 219--271.

24. \textit{А. Г. Витушкин, В. В. Ежов, Н. Г. Кружилин} Продолжение локальных отображений псевдовыпуклых поверхностей,
    Докл. АН СССР, 270:2 (1983),  271–274

25. \textit{Витушкин А.Г.} Голоморфные отображения и геометрия поверхностей
        // Современные проблемы математики. Т.7.М.:ВИНИТИ, 1985.
        С. 167-226.

26. \textit{Витушкин А.Г.} Вещественно-аналитические гиперповерхности комплексных
        многообразий // Успехи матем. наук, 1985. Т. 40, N 2. С. 3 - 31.

27. \textit{Пинчук С.И.} Голоморфные отображения в $ \Bbb C^n $ и проблема
        голоморфной эквивалентности //
В книге
        "Современные проблемы математики", Фундаментальные направления,
        т.9, М.:ВИНИТИ, 1986, С. 195 - 223.

28. \textit{Winkelmann J.} The classification of 3-dimensional homoge\-neous
        complex manifolds // Lecture Notes in Math. Springer, N 1602
        (1995). P. 230.

29.\textit{Шабат~Б.~В.} Введение в комплексный анализ: Ч.2. Функции нескольких переменных. --- М.:~Наука,~1985.

30. \textit{Лобода~А.~В.} Однородные вещественные гиперповерхности в $\mathbb{C}^3$ с двумерными группами изотропии// Тр. МИАН. --- 2001. --- \textit{235}. --- С.~114--142.

31.  \textit{Лобода А.В.} Локальное описание однородных вещественных гиперповерхностей двумерного комплексного пространства в терминах их нормальных уравнений // Функц. анализ. 2000. Т. 34, №2. С. 33--42.

32.\textit{Лобода~А.~В.} О размерности группы, транзитивно действующей на гиперповерхности в $\mathbb{C}^3$// Функц. анализ и его прил. --- 1999. --- \textit{33}, \No~1. --- С.~68--71.

33. \textit{Kossovskiy I., Loboda A.} Classification of Homogeneous Strictly Pseudoconvex Hypersurfaces in $ \Bbb C^3 $//  https://arxiv.org/abs/1906.11345

34. \textit{Лобода~А.~В.} Однородные строго псевдо-выпуклые гиперповерхности в $\mathbb{C}^3$ с двумерными группами изотропии// Мат. сб. --- 2001. --- \textit{192}, \No~12. --- С.~3--24.

35. \textit{Лобода~А.~В.} Об определении однородной строго псевдо-выпуклой гиперповерхности  по коэффициентам ее нормального уравнения// Мат. заметки. --- 2003. --- \textit{73}, \No~3. --- С.~453--456.

36. \textit{Ежов~В.~В., Лобода~А.~В., Шмальц~Г.} Каноническая форма многочлена четвертой степени в нормальном уравнении вещественной гиперповерхности в $\mathbb{C}^3$// Мат. заметки. --- 1999. --- \textit{66}, \No~4. --- С.~624--626.

37. \textit{Repovs D., Skopenkov A.B., Schepin E.V.} $ C^1 $- homogeneous compacta
         in $ \Bbb R^n $ are $ C^1 $-submanifolds of $ \Bbb R^n $ //
         Proc. Amer. Math. Soc. V. 124 (1996). P. 1219 - 1226.

38. \textit{Shnobl L., Winternitz P.} Classification and Identification of Lie Algebras. CRM Monograph Series. Volume 33. 2014.
 Centre de Recherches Mathematiques, Montreal, QC, Canada. 306 pp;

39. \textit{Vrancken, L.} Degenerate Homogeneous Affine Surfaces in $\Bbb R^3 $ // Geometriae Dedicata, (1994). 53 (3), 333-351.

40. \textit{Мубаракзянов~Г.~М.} О разрешимых алгебрах Ли// Изв. вузов. Матем. --- 1963. --- \No~1. --- С.~114--123.

41. \textit{Лобода~А.~В.} Всякая голоморфно-однородная трубка в $\mathbb{C}^2$ имеет аффинно-однородное основание// Сиб. матем. журнал. --- 2001. --- \textit{42}, \No~6. --- С.~1335--1339.

42 \textit{Sabzevari M., Hashemi A., Alizadeh B. M., Merker J.} Applications of differential algebra for computing Lie algebras of infinitesimal CR-automorphisms //\\ http://arxiv.org/abs/1212.3070

43. \textit{Лобода А.В., Суковых В.И.} Использование компьютерных алгоритмов в задаче коэффициентного описания однородных поверхностей // Вестник ВГУ. Cистемный анализ. 2015. №1. C. 14--22.

44. \textit{Dadok~J., Yang~P.} Automorphisms of tube domains and spherical hypersurfaces// Amer. J. Math. --- 1985. --- \textit{107}, \No~4. --- С.~999--1013.

45. \textit{Исаев~А.~В., Мищенко~М.~А.}	Классификация сферических трубчатых гиперповерхностей, имеющих в сигнатуре формы Леви один минус// Изв. АН СССР. Сер. матем. --- 1988. --- \textit{52}, \No~6. --- С.~1123--1153.

46. \textit{Isaev A.V.} Rigid spherical hypersurfaces // Complex Variables,
         V. 31 (1996), P. 141 - 163.

47. \textit{Eastwood M., Ezhov V.} On Affine Normal Forms and a Classification of Homogeneous Surfaces in Affine Three-Space.
     Geometriae Dedicata v. 77, 11–69 (1999)

\end{document}